\documentclass{article}
\usepackage{row}
\usepackage{graphicx}
\begin{document}
\title{The Origins of Complex Geometry in the 19th Century
\thanks{The author would like to thank Howard Resnikoff for his careful reading and comments on an earlier draft of this paper}
}
\author{Raymond O. Wells, Jr.
\footnote{Jacobs University Bremen; University of Colorado at Boulder; row@raw.com}}

\maketitle
\tableofcontents

\section{Introduction}
One of the most beautiful and profound developments in the 19th century is {\it complex geometry}. We mean by this a constellation of discoveries that led to the modern theory of complex manifolds (and more generally complex spaces: complex manifolds with specified types of singularities) and modern algebraic geometry, both of which have played an important role in the 20th century. The primary aspects to the theory of complex manifolds are the geometric structure itself, its topological structure, coordinate systems, etc., and holomorphic functions and mappings and their properties. Algebraic geometry over the complex number field uses polynomial and rational functions of complex variables as the primary tools, but the underlying topological structures are similar to those that appear in complex manifold theory, and the nature of singularities in both the analytic and algebraic settings are also structurally very similar.  Algebraic geometry uses the geometric intuition which arises from looking at varieties over the the complex and real case to deduce important results in arithmetic algebraic geometry where the complex number field is replaced by the field of rational numbers or various finite number fields.  This has led to such important results in the latter half of the twentieth century, most notably, Wiles's solution of the Fermat Last Theorem problem.  

Complex geometry includes such diverse topics as Hermitian differential geometry, which plays an important role in Chern classes of holomorphic vector bundles, for instance, Hermitian symmetric domains or more generally homogeneous spaces with complex structure, or real differentiable manifolds with some complex structure in the tangent bundle such as almost complex manifolds and CR (Cauchy-Riemann) manifolds, and many other examples.  Of course a domain%
\footnote{We will use the generic word {\it domain} to mean a connected open set.}
 in the complex plane \(\BC\) was an initial example of a complex manifold, much studied in the 19th century, and that will be an important part part of the story.

During the 17th and 18th centuries mathematics experienced major developments in geometry and analysis, specifically the geometry of curves and surfaces in \(\BR^3\), following the pioneering work of Descartes and Fermat, and the flourishing of analysis after the creation of differential and integral calculus by Leibniz, Newton and others. In all of this work, geometry was restricted to real geometric objects in the Euclidean plane and  three-space.  Complex numbers, on the other hand, were developed and referred to as {\it imaginary numbers}, as they were called for several centuries, and they arose as solutions of specific polynomial algebraic equations.  In the 18th century they became part of the standard tools of analysis, especially in the development of  fundamental elementary functions, developed by Euler,  which is epitomized in his famous formulation of the complex exponential function
\[
e^z := e^x(\cos y+i\sin y),
\]
for a complex variable \(z= x+i y\). However, the study of the geometry of curves and surfaces in \(\BR^3\)  did not include complex numbers in any substantive manner, whereas in the 20th century, complex geometry has become one of the main themes of 20th century mathematical research.

The purpose of this paper is to highlight  key ideas developed in the 19th century which became the basis for 20th century complex geometry.  We shall do this by looking in some detail at some of the innovators and their initial publications on a selection of research topics that, in the end, contributed in various ways to what we now call complex geometry.

In Section \ref{sec:complex-plane} we discuss the work of the Norwegian surveyor Wessel, the French mathematician Argand and the German astronomer-mathematician Gauss, all of whom contributed to our understanding of the complex plane as the usual Euclidean plane with complex coordinates \(z=x+iy\), including its polar coordinate representation as well as expressing the distance between points in terms of complex coordinates (modulus of a complex number). Over the course of the century this understanding became universally adopted, but at the beginning of the century, it was quite unknown.

In Section \ref{sec:abel-theorem} we look in some detail at Abel's fundamental paper concerning what is now known as Abel's Theorem concerning his generalization of the addition theorem for elliptic integrals due to Euler, which was itself a generalization of the addition theorem for trigonometric functions.  This  paper became a major motivation for major work by Riemann, Weierstrass and many others in the second half of the 19th century, as we discuss in the paragraphs below.  

In the next section (Section \ref{sec:elliptic-functions}) we discuss two fundamental papers by Abel and Jacobi which created the theory of elliptic functions, the 19th century generalization of trigonometric functions (which were periodic), and these new functions were doubly periodic in two independent directions in the complex plane. Elliptic functions utilized the geometry of the complex plane in a fundamental manner, for instance in the role of the period parallelogram, whose translates covered the complex plane. This theory was developed further in the work of Cauchy, Liouville and Weierstrass, among many others, and we trace this development in some detail, as it became quite standard material in the texts at the end of the century.

A key development in the 19th century was the creation of a theory of complex-valued functions that were intrinsically defined on domains in the complex plane and this is the theory of holomorphic and meromorphic functions.  The major steps in this theory were taken by Cauchy, in his theory of the Cauchy integral theorem and its consequences, by Riemann, in his use of partial differential equations, in particular, the study of harmonic functions, and by Weierstrass with his powerful use of power series (pun intended!).  The unification of all three points of view towards the end of the 19th century had created what is now called function theory, and has been ever since a standard tool for almost all mathematical studies today.  In Section \ref{sec:holomorphic} we shall look at some of the initial papers by these innovators and see how the point of view for this important topic evolved over time.

 Finally, in Section \ref{sec:riemann-surfaces}, we come to a pivotal development in complex geometry, namely Riemann's creation of Riemann surfaces.  Riemann's paper of 1857, which we discuss in some detail in this section, takes some of the main ideas from Abel's paper on Abel's theorem concerning multivalued functions of one real variable, and creates a theory of single-valued holomorphic functions on an abstractly defined surface with complex coordinates. These surfaces are looked at from the point of view of analysis, from algebraic geometry as the solution of algebraic equations of two complex variables, and from the point of view of topology, including the important notion of connectivity of a surface, which led to later developments in algebraic topology.

The conclusion of this paper outlines some topics which are today important for complex geometry and which were also developed during the latter part of the 19th century. These include the theory of transformation groups of Lie and Klein, the development of set theory by Cantor and the subsequent developments of topological spaces by Hausdorff and Kuratowski, and the fundamental work on foundations of algebraic topology. Two other topics that are important in complex geometry are projective geometry (in particular, the development of projective space) and differential geometry.  These two topics are covered from an historical perspective (with a number of references) in an earlier paper of ours which concerned itself with important developments in geometry (not necessarily complex) in the 19th century \cite{wells2013}.  We conclude this paper by discussing briefly the creation of abstract topological, differentiable and complex manifolds in the definitive book by Hermann Weyl in 1913, who used all of the topics discussed above, and which became the cornerstone of what became complex geometry in the 20th century.

\section{The Complex Plane}
\label{sec:complex-plane}
The  well known quadratic formula
\[x= \frac{-b \pm \sqrt{b^2-4ac}}{2a},
\]
as a solution to the quadratic equation
\[
ax^2 +bx+c=0,
\]
is attributed to the Babylonians  during their very creative period of mathematical discovery (circa 1800BC to 300BC) (see \cite{vanderwaerden1963}, \cite{resnikoff-wells1984} for discussions of the splendid mathematical accomplishments of the Babylonians, mostly preceding and greatly influencing the Greek mathematicians and astronomers).  Of course they used different notation, but their understanding was clear.   This formula led to the problem of understanding what one means by the  the square root in the cases where \(b^2-4ac\) happens to be negative.  This problem has been a part of mathematical culture ever since.  By the 18th century numbers involving \(\sqrt{-1}\) were used by numerous mathematicians in the solutions of a variety of problems, and Leonhard Euler (1707-1783) introduced the well known notation
 of \(i\) to represent%
\footnote{However, we note that in the work of a number of several 19th century mathematicians, the notation \(\sqrt{-1}\) was used for emphasis, for instance in the well known dissertation of Bernhard Riemann from 1851 \cite{riemann1851}, which we will discuss in the paragraphs below.}
 \(\sqrt{-1}\)
 and gave us his famous formula involving our basic mathematical constants
\[
e^{\pi i}=-1.
\]
These numbers became known as {\it imaginary numbers}, indicating clearly that they were figments of the imagination, in some sense, but weren't real mathematical objects.  The mathematicians of the 18th century, many of whom were very interested in questions of geometry, including Euler, missed the opportunity to come up with a geometric interpretation of these imaginary numbers.  This opportunity was not missed at the beginning of the 19th century.  

There were three independent contributions to the creation of the {\it complex plane} at the beginning of the 19th century, namely by Caspar Wessel (1745--1818) in 1797 \cite{wessel1797}, Jean-Robert Argand (1768--1822) in 1806 \cite{argand1806}, and Carl Friedrich Gauss (1777--1855) in 1831 \cite{gauss1831}.  We can cite this creation of the geometric complex plane as having been the birth of complex geometry, and it took some time for this new perspective to become an ordinary part of mathematical discourse.

Wessel and Argand both wrote definitive papers on the geometric representation of complex numbers in the Cartesian plane \(\BR^2\), and neither paper was recognized at the time of publication for the great breakthrough they both represented. In the extreme case, Wessel's paper was not recognized until a century later when it was translated into French (from the original Danish). Today this paper is available in a beautiful book \cite{wessel1797} (translated into English), along with a personal and mathematical biography of Wessel.%
\footnote{In addition the book contains a detailed excellent article by Kirsti Andersen entitled {\it Wessel's Work on Complex Numbers and its Place in History} which concerns the history of the use of the plane to represent complex numbers from Wessel to Hamilton, including the contributions of numerous other mathematicians including Argand and Gauss.}   

\begin{figure}
\vspace{6pt}
\centerline{
	\includegraphics[width=6in]{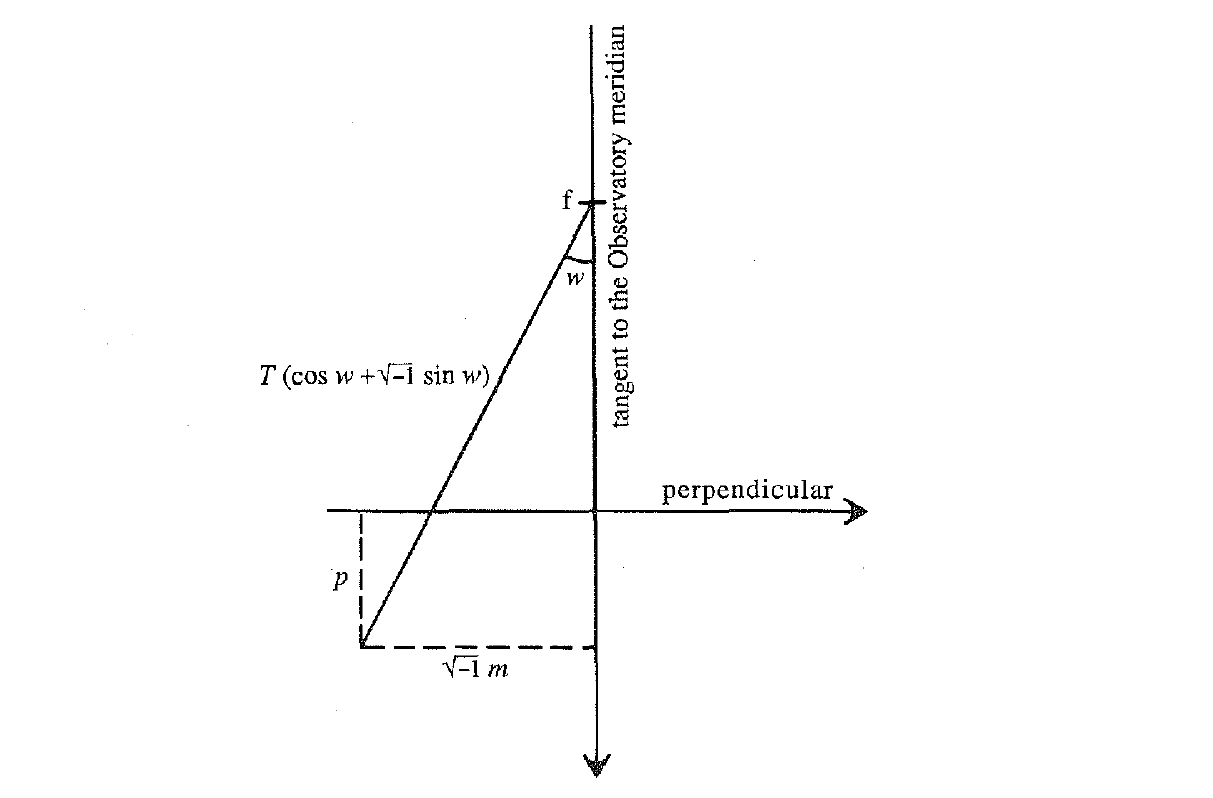}}
 \caption{A figure taken from a manuscript of Wessel called {\it Trigonometric Calculations} from 1779 (on p. 46 of \cite{wessel1797}) illustrating the use of complex coordinates.}
 \label{fig:wessel-figure}
\end{figure}

\begin{figure}
\vspace{6pt}
\centerline{
	\includegraphics[width=12cm]{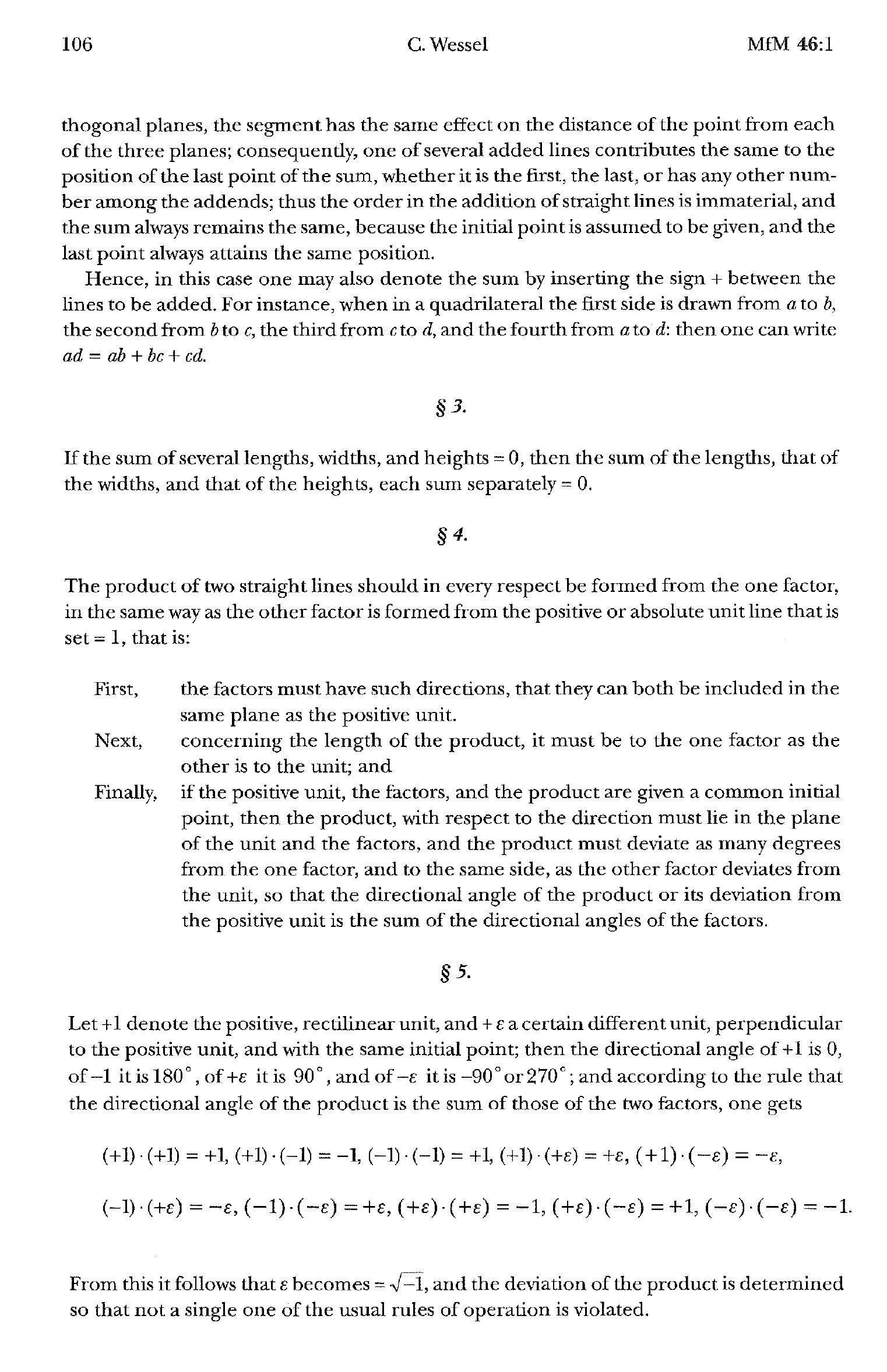}}
 \caption{Wessel's notion of sum and product of complex numbers from his paper of 1797 (English translation \cite{wessel1797}, p. 106) }
 \label{fig:wessel-p106}
\end{figure}

\begin{figure}
\vspace{6pt}
\centerline{
	\includegraphics[width=12cm]{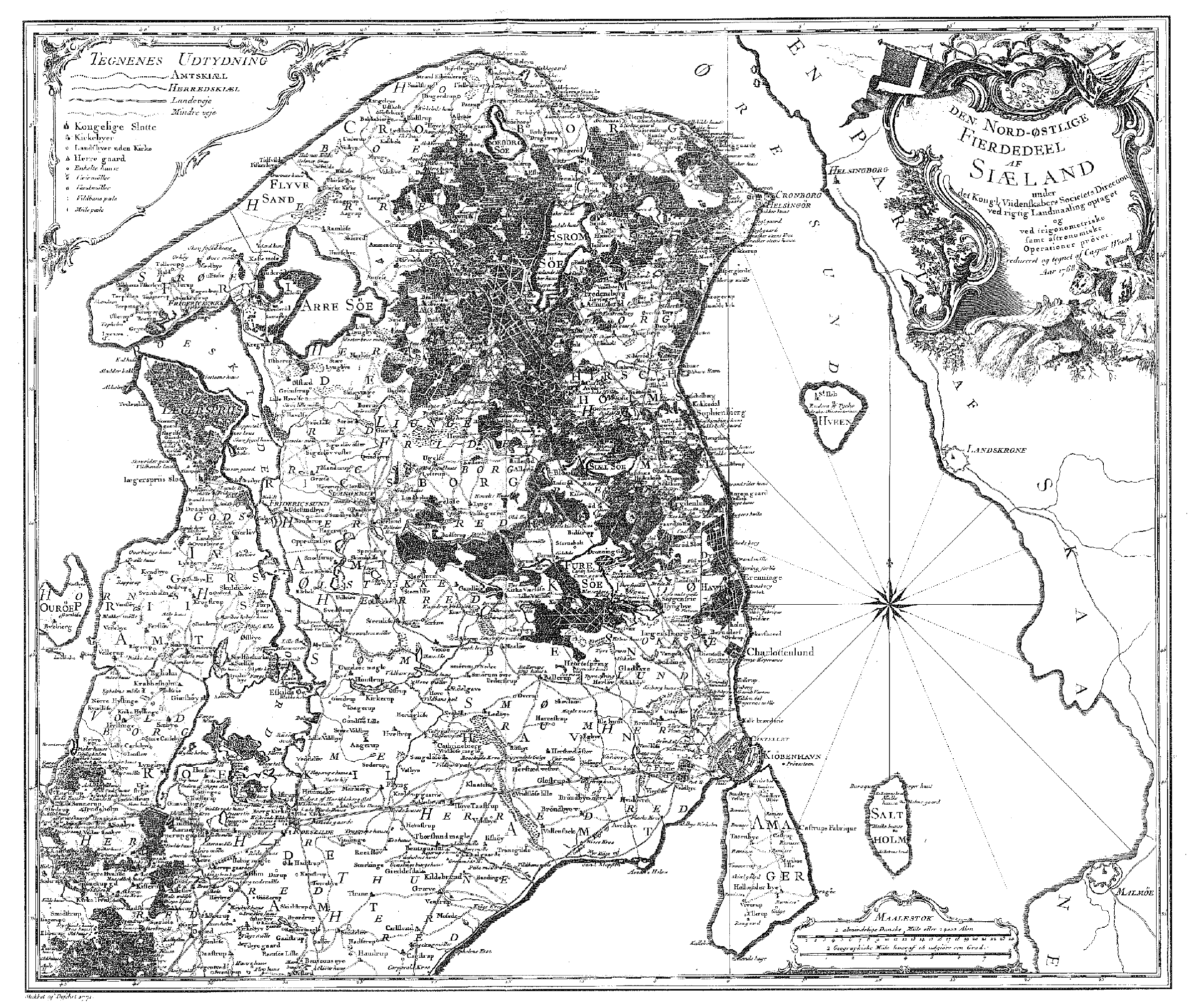}}
 \caption{Wessel's map of Denmark from 1768 (\cite{wessel1797}, plate following p. 21) }
 \label{fig:wessel-map}
\end{figure}

Wessel was a geographical and trigonometrical surveyor who surveyed large parts of Denmark and one section of Germany (Duchy of Oldenburg, northwest of Bremen, at the time under control of the Danish crown).  In fact, Gauss, in his survey of the land southeast of Oldenburg (Bremen to G\"{o}ttingen), used some of Wessel's survey data to lend accuracy to his own measurements. Wessel came upon his idea of representing complex numbers%
\footnote{The term {\it complex numbers} was introduced by Gauss in 1831 \cite{gauss1831}, although the term {\it imaginary numbers} was used till the latter half of the 19th century by many mathematicians, including, in particular, Cauchy.}
in a geometric manner as a tool for simplifying trigonometrical calculations, which were so prevalent in his surveying work. He described a complex number as a length and a direction from a given point and a given axis passing through that point, just as we do today.  More importantly he described how to add and multiply numbers using this language.  In Figure \ref{fig:wessel-figure}
we see Wessel using a polar coordinate system involving complex numbers as coordinates, and in Figure \ref{fig:wessel-p106} we show a page of the English translation from his paper of 1797 where he describes addition and multiplication of complex numbers.  Note that at the bottom of the page in Figure \ref{fig:wessel-p106} he identifies his geometric quantity \(\e\), a unit vector perpendicular to the real axis, as \(\sqrt{-1}\).
He uses this representation to give a complete description of the \(n\) roots of unity of degree \(n\) in the form:
\[
\{1, \cos(2\pi/n)+\e\sin(2\pi/n), \cos(4\pi/n) + \e\sin(4\pi/n),...\},
\] 
where \(\e\) is his notation for \(\sqrt{-1}\).. Finally, his main task in the remainder of this paper is to tackle problems of spherical geometry in three dimensions.  We note that his product of two directed line segments (see Figure \ref{fig:wessel-p106}) from a common point {\it lies on a plane spanned by the two segments}, indicating that he has been conceptualizing his ideas in three dimensions from the beginning.  We conclude this discussion of Wessel by including in Figure \ref{fig:wessel-map} a beautiful map from his earlier work, showing his skill as a cartographer.

We turn now to Argand, who published a small pamphlet \cite{argand1806} in a limited print edition in 1806 entitled {\it Essai sur un mani\`{e}re de repr\'{e}senter les quantit\'{e}s imaginaires dans les constructions g\'{e}om\'{e}trique}. This was later reprinted in an influential mathematical journal edited by Joseph Diaz Gergonne (Annales de Math\'{e}matiques pures et appliques) in 1813, which included papers by Jacques-Fr\'{e}d\'{e}ric Fran\c{c}ais (1775--1833), Fran\c{c}ois Joseph Servois (1768--1847), responses by Argand, and some commentary by Gergonne concerning the new ideas in Argand's work%
\footnote{See the article by Andersen \cite{andersen1999} for a detailed analysis of this interesting mathematical discussion in Gergonne's journal}%
. Argand also gave in this paper the first definitive proof of the fundamental theorem of algebra using his geometric representation.  Indeed, he formulated the theorem in the form that for any polynomial equation of degree \(n\) with complex numbers as coefficients,
\[
p(z)= a_0 + a_1 z + ... + a_n z^n = 0, a_j \in \BC,
\]
there exists at least one complex number \(z_0 \in \BC\) such that \(p(z_0) = 0\). This proof is not constructive and is a proof by contradiction, like many other proofs given later by others, and it utilizes substantively the notion of the {\it modulus} of a complex number,
\[
|x+ i y| := \sqrt{x^2 + y^2},
\]
which was first introduced by Argand in his paper and is the length of the directed line segment used by Wessel.

Gauss had thought about the issue of the geometric representation of complex number for some decades at the beginning of the 19th century, but didn't publish anything on the subject until his ``Second Commentary on Quadratic Residues" \cite{gauss1831} in 1831, in which he specifically defined a complex number \(z\) of the form \(x + iy\) to correspond to the point \((x,y)\) in the Euclidean two-plane \(\BR^2\), and the usual arithmetic (addition and multiplication) of complex numbers 
\bean
(a+ib) +(c+id) &= &(a+c) + i(b+d)\\
(a+ib)(c+id) &=& (ac-bd) + i ( ad+bc)
\eean
corresponded to new specific points in the plane.  In this paper, he did not consider a polar coordinate representation of complex numbers so that multiplication corresponded to multiplying moduli and adding angles of complex numbers as did Wessel and Argand, although he surely was aware of this by this time.  He was more concerned with emphasizing that this relation of arithmetic and geometry was a valid way of doing mathematics, and that had such numbers not been called``imaginary" centuries earlier, they would have been accepted much earlier.  His main purpose in this short note is to indicate that a number of his number-theoretic results from his well-known treatise on number theory from 1801 {\it Disquistiones Arithmeticae} \cite{gauss1801} could be extended to the setting of complex numbers, and specifically he discusses complex numbers of the form \(a+ib\) where \(a\) and \(b\) are integers, but emphasizing that such numbers were points in a two-dimensional plane. The only earlier reference by Gauss to a geometric representation of complex numbers was in a detailed letter from Gauss to Bessel in 1811%
\footnote{This reference and English translation is from Andersen \cite{andersen1999}.}%
, and the original letter was printed in Gauss's Werke Vol 8 in 1900 \cite{gauss1811}).  Here is what Gauss had to say:
\begin{quote}
What should one understand by \(\int \phi x\cdot dx\) for \(x=a+bi\)? Obviously, if we want to start from clear concepts,
we have to assume that \(x\) passes from the value for which the integral has to be \(0\) to \(x=a+bi\)
through infinitely small increments (each of the form \(x=a + bi)\), and then to sum all the \(\phi x \cdot dx\) .
Thereby the meaning is completely determined . However, the passage can take place in infinitely
many ways : Just like the realm of all real magnitudes can be conceived as an infinite straight line, so
can the realm of all magnitudes, real and imaginary, be made meaningful by an infinite plane, in
which every point, determined by abscissa \(=a\) and ordinate \(=b\), as it were represents the quantity \(a+bi \).
The continuous passage from one value of \(x\) to another \(a+bi\) then happens along a curve and is therefore possible in infinitely many ways. I claim now that after two different passages the integral
\(\int \phi x \cdot dx\) acquires the same value when \(\phi x\) never becomes equal to \(\infty\) in the region enclosed by the two
curves representing the two passages . This is a very beautiful theorem whose not exactly difficult proof
I shall give at a suitable occasion . 
\end{quote}
We see here in this quote also Gauss's quite specific understanding of what became known as the {\it Cauchy Integral Theorem} (Augustin-Louis Cauchy (1789-1857) )  which we will discuss in Section \ref{sec:holomorphic} below.

\section{Abel's Theorem}
\label{sec:abel-theorem}
In the 18th century trigonometric functions (often called circular functions), and the related logarithmic and exponential functions became fundamental tools of analysis.  The trigonometric functions first appeared in the work of Hipparchus of Nicaea (c.190BC-- c.120BC) in the context of spherical trigonometry for use in astronomy, and later plane trigonometry was developed and used for practical engineering and building problems. In Euler's well known text on calculus from 1748 \cite{euler1748} we see these functions used in the form we are familiar with today.  These functions and others like them were called {\it transcendental functions} in that they were a more general class of functions then the {\it rational functions}, which were ratios of polynomial functions. It's important to note that all of the important transcendental functions of the 18th century, including many of the newer transcendental functions of the 19th century (e.g., Bessel functions, Riemann zeta function, etc.) were accompanied by {\it numerical tables} of their values, so that they could be used in applied computational settings. Only with the advent of computers in the mid-20th century did the use of such tables become obsolete.

Calculus became an important tool involving calculating with symbols which could often reduce a complicated problem to a simpler one before tables of values or approximation tools (such as power series) had to be used. As was known from the beginning of the use of calculus, it was most often much simpler to differentiate a given function than to find its integral, i.e., a formula for its antiderivative.  Definite integrals of specific functions which didn't seem to have an antiderivative were studied extensively in the first half of the 19th century by the means of integration in the complex plane using Cauchy residue theory, as we will see in Section \ref{sec:holomorphic}. But toward the end of the 18th century and the first half of the 19th century a great deal of effort went into understanding specific classes of indefinite integrals.  In fact, the notation often used, \(\int f(x) dx\), meant, in more precise terms we use today, \(\int^x_0 f(x)dx\), where the lower limit (here \(0\) might be some other constant), and often a constant of integration was implied or explicitly mentioned. This notation  simply meant \(\int f(x)dx\) was a function whose derivative was \(f(x)\). 

In the 18th century it was well known that the trigonometric functions and logarithm and exponential functions were defined as integrals of specific rational or algebraic functions or inverses of such functions. For instance, 
\be
\label{eqn:trig-integrals}
\log (x) = \int \frac{dx}{x}, \;\textrm{arcsin} (x) = \int \frac{dx}{\sqrt{1-x^2}}, \;\textrm{arctan} (x) = \int \frac{dx}{1+x^2},
\ee
i.e., the derivatives of these transcendental functions were these specific rational and algebraic functions.  A function such as \(\sqrt{1-x^2}\) was often referred to in the literature of the time as an {\it irrational} function, i.e., an algebraic function (involving possible roots of rational functions)  which was not rational.%
\footnote{The notion of irrational function as used at the time didn't seem to refer to transcendental functions, which, of course, were also not rational functions.}

The question of understanding integrals of various classes of functions became an important topic in the 18th and first half of the 19th century, and this led to very important work in complex geometry, as we shall see.

First, since the creation of calculus and the fundamental theorem of calculus, it was well known how to integrate a polynomial, i.e.,  \(\int x^n dx = \frac{1}{n+1}x^{n+1}\). Moving up one step in complication, let \(R(x)\) be the field of rational functions in one real variable \(x\) and let \(r(x) \in R(x)\). Then, by basic algebra, namely, using the fact that any polynomial with real coefficients could be factored into linear and irreducible quadratic terms%
\footnote{This was well known and used regularly throughout the 18th century, but the {\it proofs} of the fundamental theorem of algebra didn't appear until the 19th century.},
and the method of partial fractions, one was able to write:
\be
\label{eqn:bernouilli}
\int r(x)dx = \int p(x)dx  + \int \sum\frac{a_j dx}{x- b_j} + \int\sum \frac{(e_k x + d_k)dx}{x^2 + e_k x +f_k},
\ee
where \(p(x)\) is a polynomial. Hence each integral of the form \(\int r(x) dx\), where \(r(x) \in R(x)\) can be reduced to a rational function and integrals of the form:
\[
\int \frac{dx}{x} = \log (x), \int \frac{dx}{1+x^2} = \textrm{arctan} (x),
\]
two transcendental functions. This general principle was formulated by Johann Bernoulli,  who published a short paper on this topic  in 1703 \cite{bernoulli1703} in which he outlined the process described above as a general algorithm for integrals of rational functions%
\footnote{In fact, in this paper Bernoulli assumed simple complex roots of a polynomial reducing (\ref{eqn:bernouilli}) to simply logarithmic terms.}.

If we now consider a rational function \(r(x,y)\) of two real variables (again a ratio of two polynomials \(p(x,y), q(x,y)\) of the two variables \(x\) and \(y\)), and let \(x\) and \(y\) be related by the quadratic equation
\[
y^2 = a + bx + cx^2,
\]
and hence, 
\[y(x) = \pm \sqrt{a+bx+cx^2},
\]
then the question arose in the 18th century, can one reduce an integral of the form
\be
\label{eqn:irrational}
\int r(x,y(x)) dx
\ee
to a sum of rational and elementary transcendental functions (i.e., trigonometric and logarithmic functions).  Special cases of this were known for some time, as in (\ref{eqn:trig-integrals}) for \(\int \frac{dx}{\sqrt{1-x^2}}\), for instance, where \(r(x,y)\) and \(y^2=1-x^2\). These kinds of problems arose in a variety of problems in elasticity, astronomy, and other sciences, and was an important motivation for finding general solutions (see Kline \cite{kline1972}, Chapter 19, for an overview of this intertwined scientific and mathematical development in the 18th century). In 1768 Euler  published an important book on integral calculus (\cite{euler1768-e342} (this was the first of three volumes; Vol. 2 was published in 1769 and Vol. 3 was published in 1770), which solved this particular problem and also set the stage for the work of Abel and Jacobi some 50 years later. Euler proved, by making a judicious changes of variables of the form \(x = x(t)\), where \(x(t)\) was an explicit rational function of \(t\), that the integral (\ref{eqn:irrational}) became
\be
\label{eqn:euler-quadratic}
\int r(x,y(x))dx = \int g(t) dt, 
\ee
where \(g(t)\) was a rational function, and hence the problem was reduced to the older one. Thus, such an integral reduced to a sum of a rational function and elementary functions, as before.

This change of variables due to Euler became later known as the rational parametrization of an algebraic curve of degree two, which we want to illustrate here due to its simplicity.  Suppose we have an algebraic curve in \(\BR^2\) of degree 2 of the form:
\[
ax^2 +bxy +cy^2 + dx +ey +f = 0.
\]
First we make a translation in the plane to make the constant term vanish and we have (in the new coordinates)
\be
\label{eqn:degree2}
ax^2 +bxy +cy^2 + dx +ey = 0.
\ee
Then the origin \((0,0)\) is a point on the curve, and we can consider the one-parameter family of straight lines of the form 
\[
y=tx, \;\textrm{for}\; t\in\BR,
\]
which will intersect the curve at both the origin and one other point on the curve for a fixed \(t\).  Substituting \(y=tx\) into (\ref{eqn:degree2}), we obtain
\be
\label{eqn:degree2a}
ax^2 +btx^2 +ct^2x^2 +dx +etx= 0.
\ee
Solving for \(x\) in terms of \(t\) we find the parametrization of the curve in terms of \(t\) to be:
\bea
\label{eqn:parametrization-x}
x(t) &=& \frac{-(d+et)}{a + bt +ct^2},\\
\label{eqn:parametrization-y}
y(t) & = & t\left(\frac{-(d+et)}{a + bt +ct^2}\right).
\eea
and from (\ref{eqn:parametrization-x}) we see that \(dx(t)=r(t)dt\), where \(r(t)\) is a rational function of \(t\). It follows then that
\be
\label{eqn:rational-curve}
\int r(x,y)dx = \int r(x(t), y(t)) r(t)dt,
\ee
when \(x\) and \(y\) are related by (\ref{eqn:degree2}). This verifies that such an integral is computable in terms of rational functions and elementary functions, Euler's result from 1768. An algebraic curve which has a parametrization in terms of rational functions of the form (\ref{eqn:parametrization-x}) and (\ref{eqn:parametrization-y}) is called a {\it rational curve}, and there are many examples of polynomials \(f(x,y)\) of degree higher than two which are also rational curves%
\footnote{For instance, there is the {\it folium of Descartes} given by \(x^3+y^3-3axy=0,\) which is parametrized by the rational functions \(x=\frac{3at}{1+t^3},y=\frac{3at^2}{1+t^3}\).}.
In the same book from 1768 \cite{euler1768-e342-2-6} Euler discussed the more difficult problem of the form
\be
\label{eqn:quartic}
\int \frac{dx}{\sqrt{A + Bx +Cx^2 +Dx^3 + Ex^4}},
\ee
or, more generally,
\be
\label{eqn:elliptic-integral}
\int r(x,y)dx,
\ee
where
\[
y^2 = A + Bx +C x^2 +Dx^3 +Ex^4.
\]
Functions of the type (\ref{eqn:elliptic-integral}) have been known since the 18th century as {\it elliptic integrals} as they originally arose in the context of computing via integration the lengths of arcs of an ellipse, just as the classical trigonometric  functions arose in conjunction with measuring the lengths of circular arcs.  Note that elliptic integrals {\it are functions} of the variable \(x\), indeed, they are transcendental functions, just as the elementary functions are, even though they are referred to as integrals.
In his original paper in \cite{euler1768} and in the text \cite{euler1768-e342-2-6} Euler discovered algebraic relations between elliptic integrals of the same type, e.g., the differential equation
\be
\label{eqn:elliptic-eqn}
\frac{dx}{\sqrt{A+ Bx + Cx^2 + D x^3 + E x^4}} = \frac{dy}{\sqrt{A +By + Cy^2 + Dy^3 + Ey^4}}
\ee
has a solution as an algebraic complete integral (an algebraic one-parameter family of algebraic curves). First Euler makes a change of variables of the form
\[
x=\frac{\a t+ \beta}{\g t +\delta},
\]
to get rid of the linear and cubic terms, reducing the problem to
\be
\label{eqn:E1DE}
\frac{dx}{\sqrt{A+  Cx^2 + Ey^4}} = \frac{dy}{\sqrt{A +Cy^2 + Ey^4}}.
\ee
Then, by several more quite nontrivial (and nonlinear) changes of variables and integrating,  he is able to produce the integral of this equation as a very specific polynomial function of degree 4 with coefficients which depend on \(A, C,\) and \(E\) and an arbitrary constant \(f\).  His solution has the form 
\be
\label{eqn:E1}
A(x^2+y^2)= f^2(A+Ex^2y^2) +2xy\sqrt{A(A+Cf^2+Ef^4)},
\ee
See \S15 of \cite{euler1768}, and he has a number of variations of this solution in this paper; we shall see a special case of this below. This relationship became known as an {\it Euler addition formula for elliptic integrals}. 

Let's illustrate this in the simpler case of
\[
\frac{dx}{\sqrt{A+Cx^2}} = \frac{dy}{\sqrt{A+C y^2}}, 
\]
which Euler had discussed earlier in his text \cite{euler1768-e342-2-5}.  He obtained a solution of the form
\be
\label{ref:E2}
y = x\sqrt{\frac{A+Cb^2}{A}} + b\sqrt{\frac{A+Cx^2}{A}},
\ee
having solved for \(y\) in terms of the other variables (here \(b\) is the constant of integration) from his solution.  Let's assume the special case of \(A=1, C= -1\), and then we have the function 
\be
\label{eqn:E3}
y= x\sqrt{1-b^2} +b\sqrt{1-x^2}
\ee
is the solution of 
\be
\label{eqn:E3DE}
\frac{dx}{\sqrt{1-x^2}}= \frac{dx}{\sqrt{1-y^2}}.
\ee
If we integrate both sides we find that
\be
\label{eqn:AF-1}
\int_0^y \frac{dt}{\sqrt{1-t^2}} = \int^x_0\frac{dt}{\sqrt{1-t^2}} + \textrm{constant}.
\ee
But from (\ref{eqn:E3}) we see that for \(x=0\), we must have \(y=b\), and thus the constant in (\ref{eqn:AF-1}) is the form
\[
\textrm{constant} = \int^b_0 \frac{dt}{\sqrt{1-t^2}},
\]
and hence
\[
\int^y_0 \frac{dt}{\sqrt{1-t^2}} = \int_0^x \frac{dx}{\sqrt{1-t^2}} + \int^b_0 \frac{dt}{\sqrt{1-t^2}},
\]
where \(x, y\), and \(b\) are related by (\ref{eqn:E3}). By relabeling the variables, as did Euler, we find the familiar 
formula
\[
\int^x_0 \frac{dt}{\sqrt{1-t^2}} + \int_0^y \frac{d}{\sqrt{1-t^2}} = \int^b_0 \frac{dt}{\sqrt{1-t^2}},
\]
where
\[
b= x\sqrt{1-y^2} +y\sqrt{1-x^2},
\]
which is the classical addition formula for the inverse trigonometric functions
\[
\label{eqn:arcsin-addition}
\textrm{arcsin}(x) + \textrm{arcsin}(y) = \textrm{arcsin}(b),
\]
which becomes 
\be
\label{eqn:AFtrig}
\sin (x+y) = \cos (x) \sin (y) + \sin(x) \cos (y).
\ee
Thus Euler's solution of the equation (\ref{eqn:E3DE}) yields the classical addition formula (\ref{eqn:AFtrig}) for circular functions, which was known to the ancient trigonometers.  The corresponding half-angle formulas allowed the Greek astronomers to compute the trigonometric tables which were so critical for their astronomical research. 

Niels Henrik Abel (1802--1829) in his very short lifetime%
\footnote{He was not yet 27 years old when he died.}
wrote a number of quite important papers, several of which came to play an important role in the development of complex geometry.  We will discuss two of these papers: \cite{abel1826}, his paper on what is now referred to as Abel's Theorem in algebraic geometry and \cite{abel1827}, his foundational paper on elliptic functions, the doubly periodic functions in the complex plane that generalized the classical periodic functions of trigonometry.  Both of these papers were influenced by the work of Euler which was described in the paragraphs above as well as the follow-up to Euler by Adrien--Marie Legendre (1752-1833) in his several decades long study of elliptic integrals and their applications.  

Legendre's principal contributions were contained in three monographs he published in the decade before Abel's work. These three volumes were entitled {\it Exercices de Calcul Int\'{e}gral}.  Volume 1 \cite{legendre1811} in 1811 was his major theoretical work on elliptic integrals, which showed how all elliptic integrals of a general kind could be reduced, via algebra and calculus, to three specific types of integrals, which Legendre referred to as integrals of the first, second and third kind, which we shall see shortly. Volume 2 \cite{legendre1817} from 1817 contained a major survey of approximation methods, methods of creating tables and numerous applications to geometry and applied mathematics, in particular to mechanics.  Volume 3 \cite{legendre1816} (which was actually published in 1816 before Volume 2) contains detailed tables for elliptic functions of the first and second kind and their logarithms, as well as a discussion of the issues of reducing computations of some integrals of the third kind to those of the first and second kind (there were too many free parameters in these transcendental functions of the third kind to allow the creation of reasonable tables).  After the groundbreaking work of Abel and Jacobi in 1826 and 1827 he continued his surveys of the development of what has now become the theory of elliptic functions.

The first paper of Abel we want to mention was presented to the French Academy of Science in 1826 (\cite{abel1826}) and was finally published posthumously in 1841.  It gives a vast generalization of the addition formula for elliptic integrals that was due to Euler and discussed in the paragraphs above and is now called {\it Abel's Theorem} in algebraic geometry.%
\footnote{There are a number of theorems known as Abel's Theorem in different parts of mathematics, e.g., on the convergence of power series, on the unsolvability of quintic polynomial equations, etc.}
In his second major paper \cite{abel1827}, published in Volumes 2 and 3 of Crelle's journal in 1827 and 1828, Abel wrote a definitive and foundational paper on elliptic functions.  The title of this paper, {\it Recherches sur les fonctions elliptique}, is misleading, and at the same time, so very appropriate. What he meant in the title by ``elliptic functions" were the transcendental functions studied by Euler and Legendre, etc., which were defined by and known as elliptic integrals. In this paper he introduced, for the first time, the {\it inverses} of the elliptic integral functions, and these became the now familiar doubly periodic meromorphic functions on the complex plane known as elliptic functions, that we will discuss in the forthcoming paragraphs.  So the title is absolutely correct in modern times, even if Abel didn't know it at the time!

Let us preface our formulations of Abel's Theorem%
\footnote{There are more than one algebraic-geometric theorems referred to historically over the past century as {\it Abel's Theorem}.  The very informative paper by Stephen Kleiman entitled {\it What is Abel's Theorem anyway"} \cite{kleiman2004} discusses four variants of what have been called Abel's Theorem. This paper is an article in a beautiful book \cite{laudal-piene2004} representing the proceedings of a conference held in honor of the mathematical legacy of Abel in 2002, 200 years after his birth in 1802.} with a specific version of Euclid's addition formula for elliptic integrals.  Namely, in 1761 \cite{euler251} Euler studied the differential equation 
\be
\label{eqn:EIE}
\frac{dx}{\sqrt{1-x^4}} = \frac{dy}{\sqrt{1-y^4}},
\ee
a special case of  (\ref{eqn:elliptic-eqn}) discussed briefly above, and he finds the complete algebraic integral to be
\be
\label{eqn:ECI}
x^2 +y^2 +c^2x^2y^2 = c^2+2xy\sqrt{1-c^2},
\ee
where \(c\) is the constant of integration. Now consider the specific elliptic integral  
\be
E(x):=\int^x \frac{dx}{\sqrt{1-x^4}},
\ee
where there is some fixed lower limit of integration, which in this case we choose to be \(x=0\) . 
Then one finds by integrating each side of (\ref{eqn:EIE}) that
\[
E(x) = E(y) + C,
\]
where \(C\) is a constant.  From the complete integral of (\ref{eqn:EIE}) given by (\ref{eqn:ECI}) we see that when \(x=0\), then \(y=c\) (we take the positive square root in this case, for convenience), and hence
\[
0 = E(0) = E(c) +C,
\]
and hence \(C= -E(c)\), yielding
\[E(x) = E(y) -E(c),
\]
or 
\[
E(x) +E(c) = E(y).
\]
Changing the names of the variables \(x_3= y, x_1= x, x_2= c\), we obtain the addition theorem for this particular elliptic integral of the form
\be
\label{eqn:AEI}
E(x_1) + E(x_2) = E(x_3),
\ee
where 
\[
{x_1}^2 +{x_3}^2 +{x_1}^2{x_2}^2{x_3}^2 = {x_2}^2 + 2{x_1}{x_3}\sqrt{1-{x_2}^2},
\]
which gives after squaring
\[
4{x_1}^2{x_3}^2(1-{x_2}^2) = ({x_1}^2 + {x_3}^2 +{x_1}^2{x_2}^2{x_3}^2 -{x_2}^2)^2,
\]
a polynomial relation of degree 12 among the arguments of the three transcendental functions \(E(x_1), E(x_2), E(x_3)\) (see \cite{kleiman2004}, p. 20 for various references to this formula). Note that for the arcsine addition formula (\ref{eqn:arcsin-addition}), which we can write as
\[
\int_0^{x_1} \frac{dx}{\sqrt{1-x^2}} + \int_0^{x_1} \frac{dx}{\sqrt{1-x^2}} = \int_0^{x_3} \frac{dx}{\sqrt{1-x^2}},
\]
 we have the same sort of algebraic relation which takes the (familiar) form
\[
x_3= x_1\sqrt{1-{x_2}^2} +x_2\sqrt{1-{x_1}^2},
\]
which we discussed earlier, and which when squared twice yields a polynomial relation among the three arguments of these three transcendental functions of  degree six.

Let now \(r(x,y)\) be a rational function, and let \(f(x,y)\) be a polynomial, and let \(y(x)\) be the implicit (multivalued) function defined by \(f(x,y)=0.\) The general {\it Abelian integral} is defined to be
\be
\label{eqn:abelian-integral}
A(x):=\int^x r(x,y(x))dx.
\ee
What Abel originally meant by this was an antiderivative (as did Euler), i.e. \(A(x)\) is a function whose derivative is \(r(x,y(x))\), and we are expressing this as a definite integral from an initial point (unspecified) to an upper limit \(x\), using the same symbol \(x\) as the variable of integration.  A first version of Abel's Theorem asserts that if \(g(x,y)\) is an auxiliary polynomial, and if the curve \(g(x,y)= 0\) intersects the curve \(f(x,y)=0\) in the points \((x_1,y_1),..., (x_N,y_N)\), then there are rational functions \(u,v_1,...,v_r\) of the variables \(x_1,...,x_N\) and the coefficients of the polynomial \(g(x,y)\) so that 
\be
\label{eqn:EAT}
A(x_1) +A(x_2)+ ... + A(x_N) = u +k_1\log v_1 + ... + k_r\log v_r,
\ee
where \(k_1,...,k_r\) are constants.  This says that the left hand side of (\ref{eqn:EAT}), a sum of \(N\) transcendental functions, is an {\it elementary function}, i.e., in this case a sum of a rational function and logarithmic terms. Thus (\ref{eqn:EAT})  says that such a sum of Abelian integrals is an elementary function.  Note that this is a generalization of the much simpler case that the integral
\[
\int^x r(x,y(x))dx
\]
is the sum of elementary functions, when \(y^2=Ax^2 +Bx +C\), i.e., in the trigonometric case (Euler's theorem, see (\ref{eqn:rational-curve}))%
\footnote{Note that there is no auxiliary polynomial \(g(x,y)\) in this simple case.}%
. This version of Abel's theorem (\ref{eqn:EAT}) is sometimes referred to as the {\it elementary addition theorem}, i.e., a specific sum of Abelian integrals is an elementary function (see Kleiman \cite{kleiman2004}).  

The more general version of Abel's Theorem, often known as the Abel Addition Theorem (see again \cite{kleiman2004}), asserts that, given \(r(x,y)\) and \(f(x,y)\) as before, then there is an integer \(p\ge 0\), depending only on \(f\), so that, for any set of points \(\{x_1,,...,x_\a\}\), there are points \(\{y_{1},..., y_{p}\}\) so that
\be
\label{eqn:AT}
A(x_1) + A(x_2) + ... + A(x_\a) = e(x_1,...,x_\a) + A(y_1) + A(y_2) + ... + A(y_p),
\ee
where \(e\) is an elementary function of \((x_1,...,x_\a)\) and \(y_1,...,y_p\) are algebraic functions of \((x_1,..., x_\a)\). Note that in (\ref{eqn:EAT}) we have only elementary functions on the right-hand side, and in the special elliptic integral case \(r(x,y)=1/y, f(x,y) = y^2-x^4-1,\) (\ref{eqn:AEI}), there is only one elliptic integral on the right-hand side (no elementary functions). In this case we had \(\a = 2\), but we could have iterated (\ref{eqn:AEI}) and had any number of terms on the left-hand side and still one term on the right-hand side.  Thus, in this case, \(p\) for \(f=x^2-x^4-1\), seems to be equal to 1, and that indeed turns out to be the case.  We will discuss the significance of the integer \(p\) in Abel's Theorem (\ref{eqn:AT}) somewhat later in this section.

There are two major issues in understanding or interpreting Abel's two theorems here (\ref{eqn:EAT}) and (\ref{eqn:AT}). The first is the multivalued nature of \(y(x) \) as implicitly defined by the equation \(f(x,y)=0\), and the second is: what does the integral 
\be
\label{eqn:integral}
\int^x_{x_0} (r(x,y(x))dx
\ee
mean? Here we are now thinking of the integral in (\ref{eqn:integral}) as a definite integral from some fixed point \(x_0\) to some variable end point \(x\). As we mentioned earlier, Abel thought in terms of antiderivatives and differentiation, and his proofs involve differentiation, the fundamental theorem of calculus, the implicit function theorem, and, quite importantly, the general fact, apparently quite well known at the time, that a rational symmetric function of the roots of a polynomial was a rational function of the coefficients of the polynomial (a result due to Vandermonde \cite{vandermonde1771}, as pointed out by Kleiman \cite{kleiman2004}).  This was used repeatedly by Abel to reexpress various (symmetric) functions of the multivalued functions as single-valued functions.  

Abel's work in this early part of the 19th century led to vigorous work in the latter half of that same century to understand better this issue of the multivalued functions appearing in his work; the most decisive next steps were taken by Bernhard Riemann (1826--1866)  \cite{riemann1857} in 1857, as we shall see later in Section \ref{sec:riemann-surfaces}. One aspect of the integration issue that was recognized by Abel, and which was definitively pursued by Riemann was the fact that the integral \(\int_{x_0}^x r(x,y(x))dx \) could have different values depending on the path one took from the initial point \(x_0\) to the final point \(x\). On the real line this seems to be only one path, but one could specify which signs to use in any formula for \(y(x)\) involving various combinations of radicals, for instance.  

The possible ambiguities in this integral became known as {\it periods} of the integral, as differences of two such integrals were specific multiples of fixed entities. At the time of Riemann and later, the variables \((x,y)\) were interpreted as complex numbers, and the integral (\ref{eqn:integral}) was considered as a complex path integral from \(x_0\) to \(x\) along some complex path \(\g\). Whether the integral along two different paths was the same or not became a major subject of study in complex analysis (Cauchy's integral theorem and residue theory) and in what became algebraic topology (whether the two paths bounded a simply-connected domain or not).  Both topics became major research areas in second half of the 19th century.

Finally, we want to discuss the significance of the integer \(p\) in Abel's Theorem (\ref{eqn:AT}). First, let us quote from p. 172 of Abel's paper \cite{abel1826}, where he denoted the Abelian integrals \(A(x_j)\) as \(\psi_j x_j\), and \(p\) was the difference of the two integers \(\m\), the total number of variables and integrals appearing in the theorem, and the the integer \(\a\), the number of variables (and integrals) appearing on the left-hand side of the theorem:
\begin{quote}
Dan cette formule les nombre des fonctions \(\psi_{\a+1}x_{\a+1}\),\(\psi_{\a+1}x_{\a+2}\),...,\(\psi_{\m}x_{\m}\) est tr\`{e}s-remarquable. Plus il est petit, plus la formule es simple. Nous allons, dans ce qui suit, chercher la moindre valeur dont ce nombre, qui est eprim\'{e} par \(\m-\a\), est susceptible.%
\footnote{In this formula the number of functions \(\psi_{\a+1}x_{\a+1},\psi_{\a+1}x_{\a+2},...,\psi_{\m}x_{\m}\) is very remarkable. Moreover, it is small and the formula is simple. We shall, in that which follows, search for the the smallest value for  which this number, which is expressed by \(\m-\a\), can be attained.}
\end{quote}

Strangely enough, Abel never expressed this number, which we have called \(p\), by a single symbol, in spite of the significance he did attribute to this integer, which only depended on the polynomial \(f(x,y)\). Abel proceeds to derive formulas which allow him to compute this number in various special cases, and we mention  three such cases here.  The first is the most complex.  Namely, consider a polynomial \(f(x,y)\) of degree 13, i.e.,
\[
f(x,y)= p_0 +p_1y+p_2y^2 +...+p_{12}y^{12} + y^{13},
\]
where the degrees of the polynomials (in the variable \(x\)) \(p_0, p_1,...,p_{12}\) are
\[
2,3,2,3,4,5,3,4,2,3,4,1,1.
\]
In this case, after four pages of computation (pp. 181-185 of \cite{abel1826}), Abel obtains \(p= 38\).  This number \(p\) turns out to be the celebrated {\it genus} of the algebraic curve defined by \(f(x,y) = 0\), and is a topological invariant  of the Riemann surface (and topological manifold) defined by the algebraic curve. Riemann formulated the genus in the more modern sense a half-century later.  Note that the definition of genus as defined by Abel was an invariant of the analytical data he had at his disposal, and became later a topological invariant in the hands of Riemann.  

In the case that
\[
f(x,y)= y^2-\phi(x),
\]
the {\it hyperelliptic} case, which was studied extensively by Abel in \cite{abel1828}, one finds that if \(d = \deg \phi\), where we assume that \(\phi\) has distinct roots, then
\[
p = \left\{
\begin{array}{l}
(d-1)/{2},\; \textrm{if}\; d \;\textrm{is odd,}\\
(d-2)/{2},\; \textrm{if}\; d \;\textrm{is even.}\\
\end{array}\right.
\]
So, if we have an {\it elliptic curve} in this hyperelliptic case, i.e., \(d= 3\) or 4, then \(p =1\), which means topologically (as we learn later from Riemann \cite{riemann1857}) that the elliptic curve is a two-dimensional torus.  In this case any sum of Abelian integrals (these would be now {elliptic integrals}) is the sum of one such elliptic integral plus an elementary function (as in the special case of (\ref{eqn:AEI}) above).

Our final and simplest example is the case \(y^2=Ax^2+Bx+C\), which gives \(p=0\).  This means that the underlying Riemann surface is the Riemann sphere, which is, topologically, a simple two-sphere.  We mention again in this very special hyperelliptic case that since \(p\) = 0,  the right-hand side of Abel's theorem (\ref{eqn:AT}) contains no Abelian integrals, only elementary functions, as we know from the earlier work of Euler discussed above (\ref{eqn:rational-curve}) on the rational parametrization of an algebraic curve of degree two.

One final note is that an Abelian integral is called of the {\it first kind}, if the integral is finite for all \(x\).  This terminology was introduced by Legendre in the case of elliptic integrals in \cite{legendre1811}.  For instance, the following Abelian integrals in the hyperelliptic case (where \(f(x,y) = y^2-\phi(x)\) and \(\f(x)\) has distinct roots) are of the first kind, where \(p\) is again the genus of the hyperelliptic curve \(f(x,y) = 0\),
\be
\label{eqn:first-kind}
\int^x_{x_o} \frac{dx}{\sqrt{\phi(x)}}, \int^x_{x_o} \frac{xdx}{\sqrt{\phi(x)}}, ..., \int^x_{x_o} \frac{x^{p-1}dx}{\sqrt{\phi(x)}}.
\ee
In this case these \(p\) Abelian integrals of the first kind in (\ref{eqn:first-kind}) are linearly independent and they span the space of all such Abelian integrals of the first kind (see Markushevich \cite{markushevich1992}).  We will see this in greater detail later in this section. Note that the genus \(p\) appears here explicitly, and the dimension of this vector space of all Abelian integrals of the first kind can be used as a second and equivalent definition of genus in this case. 

\section{Elliptic Functions}
\label{sec:elliptic-functions}

We now turn to the second major paper by Abel \cite{abel1828} which developed the theory of elliptic functions. This was followed one year later by the equally definitive and independent work by Carl Gustav Jacob Jacobi (1804-1851) \cite{jacobi1829} on precisely the same subject (Jacobi had published a shorter introduction to his work at the end of 1827 \cite{jacobi1827}). These two  long papers laid the foundation for the rich development of the theory of doubly-periodic functions in the complex plane that was pursued by numerous mathematicians throughout the 19th century in a wide variety of forms (complex analysis, algebraic geometry, number theory, etc.).

However, before we look at Abel's and Jacobi's work, let's briefly review what functions of a complex variable meant to mathematicians at the beginning of the 19th century.  As we saw in Section \ref{sec:complex-plane},
the geometric representation of complex numbers in the complex plane had not yet been developed. Complex numbers were simply algebraic combinations of real numbers with the imaginary unit \(i = \sqrt{-1}\) of the form \(a+ib\) manipulated according to the well known rules of addition and multiplication of such numbers.  In reading through the  work of Euler from the mid-18th century \cite{euler1748} that we have cited in our earlier chapter, one sees that imaginary numbers arose from solving algebraic equations and were manipulated by the usual rules of algebra. A function \(f\) of a complex variable \(x+iy\) for rational function computed \(f(x+iy)\) by algebra, i.e.,
\[
(x+iy)^2 = x^2-y^2 +i(2xy),
\]
and a series of the form
\[
\sum_{n=0}^\infty a_n(x+iy)^n
\]
would be expressed in terms of its real and imaginary parts by term-by-term algebra.  For transcendental functions we find a pregnant remark of Euler on p. 96 of \cite{euler1796} (the 1796 French edition of his analysis book from 1748) which says (in English translation), where here \(x\) is a real number,
\begin{quote}
Since \(\sin^2x +\cos^2x =1\), in decomposing into factors, one would have
\[
\cos x + \sqrt{-1}\sin x)(\cos x - \sqrt{-1} \sin x) = 1.
\]
These factors, although imaginary, are of great usage in the combination and multiplication of arcs [radian angles].
\end{quote}
A few pages later in the same book Euler observes that (now letting \(i = \sqrt{-1}\) , for convenience),
by letting
\[
e^{ix} = \cos x +i \sin x,
\]
then
\bean
\cos x = \frac{e^{ix} + e^{-ix}}{2},\\
\sin x = \frac{e^{ix} - e^{-ix}}{2}.
\eean
Using the addition formula for exponentials he then obtains (by definition)
\[
e^{x+iy}= e^xe^{iy} = e^x(\cos x + i \sin y),
\]
with similar expressions for the transcendental functions of a complex variable \(\sin(x+iy)\), \(\cos(x+iy)\), etc.. 
These are then examples of transcendental functions of a complex variable represented as algebraic combinations (involving the imaginary unit \(i\)) of real-valued functions of a complex variable. 

This was the stage that was set for Abel and Jacobi as they set out to create their theories of elliptic functions (which would also be formulated initially as algebraic combinations of real-valued functions, just as Euler did with the trigonometric functions).  Let us now formulate what an elliptic function is in the standard language of complex analysis.  Namely, let  \(\omega_1\) and \(\omega_2\) be two fixed complex numbers such that \(\Im(\omega_1/\omega_2)\ne 0\), then an {\it elliptic function} \(f(z)\) with the two periods \(\omega_1\) and \(\omega_2\) is a meromorphic function on the complex plane \(\BC\) such that
\[
f(z+m\omega_1 +n\omega_2) = f(z), \;\textrm{for all}\; n, m \in \BZ,
\]
where \(\BZ\) denotes the ring of integers. We say that such a function is {\it doubly-periodic with the two periods \(\omega_1\) and \(\omega_2\)}.  This is completely analogous to the simply-periodic functions from trigonometry, where, for instance,  
\[\sin(x+2\pi n)=\sin(x), \;\textrm{for all}\; n\in \BZ,
\]
 with the period \(2\pi\). Abel and Jacobi gave the first examples of such doubly-periodic functions and proved many of their important properties, as well as giving a variety of ways to represent such functions (power series, infinite products, etc.). The theory of trigonometric functions was a model for both of them.

Let's start with Abel's paper \cite{abel1828}, and we will follow the notation and normalizations used in his paper, although the formalism and notation of Jacobi became the standard in the literature in the following decades.  Due to Abel's early death, he was not able to participate in the later developments.  The basic idea of both mathematicians was to study the {\it inverse} of the elliptic integral functions, that had been studied extensively by their predecessors.  In this manner the addition theorems for elliptic integrals, \`{a} la Euler, became addition formulas for the elliptic functions, which generalized the addition formulas for trigonometric functions. Let us note that if one starts with the transcendental function
\[
\arcsin (x) := \int ^x_0 \frac{dx}{\sqrt{1-x^2}},
\]
then one can define its inverse \(\sin(x)\) and obtain the full theory of trigonometric functions.  This is, in effect, what Abel and Jacobi do in the elliptic integral context.

Abel begins in \cite{abel1828} by recalling the work of Euler and Legendre that we discussed in the preceding paragraphs. He notes that every elliptic integral of the form
\[
\int \frac{R(x)dx}{\sqrt{\a +\beta x +\g x^2 +\d x^3 +\e x^4}},
\]
where \(R(x)\) is a rational function can be reduced to 
\[
\int \frac{P(y)dy}{\sqrt{a + b y^2 +c y^4}},
\]
where \(P(y)\) is a rational function of \(y^2\).  This can, in turn, be reduced to the form
\[
\int \frac{A+By^2}{C+dy^2}\frac{dy}{\sqrt{a+by^2+cy^4}},
\]
and by yet one more change of variables, this can be reduced to the trigonometric form
\[
\int\frac{A+B\sin^2\theta}{C+D\sin^2\theta}\frac{d\theta}{\sqrt{1-c^2\sin^2\theta}},
\]
where \(c\) is real and \(|c| < 1\).  Finally, Abel notes that (all of this is from Legendre's book \cite{legendre1811}), every elliptic integral, by this type of reduction, can be reduced to the three cases:
\[
\int \frac{d\theta}{\sqrt{1-c^2\sin^2\theta}},
\int d\theta\sqrt{1-c^2\sin^2\theta},
\int \frac{d\theta}{(1+n^2\sin^2\theta)\sqrt{1-c^2\sin^2\theta}},
\]
which Legendre calls elliptic integrals of the first, second and third kind.  Abel decides to concentrate on the elliptic integrals of the first kind and on p. 164 of \cite{abel1828}, after the brief introduction outlined above, he says
\begin{quote}
Je me propose, dans ce m\'{e}moire, de considerer le fonction inverse, c'est a dire la fonction \(\phi\a\), determin\'{e}e par les \'{e}quations
\[
\ba{l}
\a = \int \frac{d\theta}{\sqrt{1-c^2\sin^2\theta}},\\
\sin \theta = \phi\a = x.
\ea\]
\footnote{I propose, in this memoir, to consider the inverse function, that is to say the function \(\phi\a\) determined by the equations
\[
\ba{l}
\a = \int \frac{d\theta}{\sqrt{1-c^2\sin^2\theta}},\\
\sin \theta = \phi\a = x.
\ea
\]}
\end{quote}
Abel then considers specifically the elliptic integral of the first kind in the form
\be
\label{eqn:legendre-integral}
\a(x) = \int_0^x \frac{dtx}{\sqrt{1-t^2}\sqrt{1-c^2t^2}},
\ee
in terms of the variable \(x\), where again \(c^2 >0\)%
\footnote{Abel doesn't distinguish between the upper limit of the integral and the variable of integration, but we do to clarify the discussion.}%
.  Now Abel makes two changes in notation to suit his purposes.  He replaces \(c^2\) by \(-e^2\) and replaces the term \(\sqrt{1-x^2}\) by \(\sqrt{1-c^2x^2}\) for symmetry, and finally considers the specific elliptic integral of the first kind in the form
\be
\label{eqn:EID}
\a(x) = \int^x_0\frac{dt}{\sqrt{1-c^2t^2}\sqrt{1+e^2t^2}}.
\ee
We let  \(x(\a)\) be the inverse of \(\a(x)\) given by (\ref{eqn:EID}), which is well defined near \(x=0\), and Abel defines \(\phi(\a)\) to be this inverse \(x(\a)\) on a suitable interval containing \(x=0\).  The derivative of \(\a(x)\) is simply given by 
\[
\a'(x) = \frac{1}{\sqrt{1-c^2x^2}\sqrt{1+e^2x^2}},
\]
and, by the inverse function theorem, the derivative of \(\phi(\a)\) is given by
\be
\label{eqn:phi-derivative}
\phi'(\a) = \sqrt{1-c^2\phi(\a)^2}\sqrt{1+e^2\phi(\a)^2}.
\ee
Then Abel introduces two additional functions of \(\a\) defined by
\be
\label{eqn:E2}
f(\a) := \sqrt{1-c^2\phi(\a)^2}, F(\a) := \sqrt{1+e^2\phi(\a)^2},
\ee
which appear in (\ref{eqn:phi-derivative}), yielding \(\phi'(\a) = f(\a)F(\a)\). 
These {\it three} functions of a real variable%
\footnote{The inverse function \(\f(\a)\) and its related functions \(f(\a)\) and \(F(\a)\) are well defined locally near \(\a=0\) by the inverse function theorem.  The extension to the full real line is discussed later in this section.} \(\a\) are the generalizations of the {\it two} trigonometric functions \(\sin(\a)\) and \(\cos(\a)\), and, as Abel says on p. 265 of his paper:
\begin{quote}
Plusieurs proprie\'{e}t\'{e}s de ces fonctions se d\'{e}dusierent imm\'{e}diatement des propri\'{e}t\'{e}s connues de la fonction elliptique [elliptic integral] de la premi\`{e}re esp\`{e}ce, mais d'autres sont plus cach\'{e}es.  Par exemple on d\'{e}montre que les \'{e}quations \(\phi\a=0, f\a=0, Fa=0\) on un nombre infini de racines, qu'on peut trouver toutes. Une des les plus remarquables est qu'on peut exprimer rationellement \(\phi(m\a), f(m\a), F(m\a)\) (\(m\) un nombre entier) en \(\phi\a, f\a, Fa\). 
Aussi rien n'est plus facile que de trouver \(\phi(m\a), f(m\a), F(m\a)\), lorsqu'on conna\^{i}t \(\phi\a, f\a, F\a\); mais le probl\`{e}me inverse, savoir de d\'{e}terminer \(\f\a,f\a, F\a\) en \(\f(m\a), f(m\a), F(m\a)\), est plus difficile, parcequ'il d\'{e}pend d'une equation d'un degr\'{e} \'e{e}lev\'{e} (savoir du degr\'{e} \(m^2\)).

La r\'{e}solution de cette \'{e}quation es l'objet principal de ce m\'{e}moire. D'abord on fera voir, comment on peut trouver toutes les racines, au moyen des fonctions \(\f, f, F\). On traitera ensuite de la r\'{e}solutions alg\'{e}brique de l'\'{e}quation en question, et on parviendra \`{a} ce r\'{e}sultat remarquable, que \(\f\frac{\a}{m}, f\frac{\a}{m}, F\frac{\a}{m}\) peuvent \^{e}tre exprim\'{e}es en \(\f\a,f\a,F\a\), par une formule qui, par rapport \`{a} \(\a\), ne contient d'autre irrationalit\'{e} que des radicaux. Cela donne une classe tr\`{e}s g\'{e}n\'{e}rale d'quations qui sont r\'{e}soluble alg\'{e}briquement.%
\footnote{Several properties of these functions are deducible immediately from the known properties of the elliptic function [elliptic integral] of the first kind, but others are more hidden. For example, one can show that the equations \(\phi\a=0, f\a=0 , Fa=0\) has an infinite number of roots, where one can find all of them. One of the most remarkable properties is that one can express rationally \(\phi(m\a), f(m\a), F(m\a)\) (\(m\) an integer) in \(\phi\a, f\a, Fa\). Also, nothing is more simple than to find 
\(\phi(m\a), f(m\a), F(m\a)\), when one knows \(\phi\a, f\a, F\a\); but the inverse problem, to know how to determine 
\(\f\a,f\a, F\a\) in  \(\f(m\a), f(m\a), F(m\a)\), is more difficult, since it depends on an equation of high degree (more specifically of degree \(m^2\)). 

The solution of this equation is the principal object of this memoir.  At first  one can see how one can find all the roots, by means of the functions \(\f, f, F\). Then one treats the algebraic solution of the equation in question, and one comes to this remarkable result, that \(\f\frac{\a}{m}, f\frac{\a}{m}, F\frac{\a}{m}\) can be expressed in \(\f\a,f\a,F\a\), by a formula, which, with respect to \(\a\), doesn't contain any irrationality except radicals. This gives a very general class of equations which are solvable algebraically.
}
\end{quote}

We note that this last comment of Abel's about solvability of high degree equations by means of extracting roots relates to one of his first papers \cite{abel1824} in which he shows for the first time the unsolvability in terms of radicals of generic algebraic equations of degree 5 or higher, a problem that had been outstanding for a long time. The definitive work on whether a given polynomial equation was solvable in terms of radicals was due to \'{E}variste Galois (1811-1832) in his work which established the now well-known Galois theory.  This was published in 1846 \cite{galois1846}, 14 years after his very early death at the age of 20.

Let us now turn to Abel's construction of his version of elliptic functions and their first fundamental properties.  He first defines each of these functions for all real values of \(\a\) in a specific interval around the orign, and then proceeds to define them as functions of a complex variable \(\a+i\beta\) on the entire complex plane in a sequence of steps.  First he sets
\[
\frac{\omega}{2}:= \int^{\frac{1}{c}}_0\frac{dt}{\sqrt{1-c^2t^2}\sqrt{1+e^2t^2}},
\]
where it is simple to verify that the limiting integral at the singular point \(x=\frac{1}{c}\) is well-defined.
Thus one sees that \(\f(\a)>0\) on \((0,\omega/2)\), and \(\f(0)=0\) and \(\f(\omega/2)= 1/c\).  Also, from the definition of \(\f(\a)\), one sees that \(\f(-\a)= -\f(\a)\), and thus we have \(\f(\a)\) well defined on \([-\omega/2,\omega/2]\), and similarly for \(f(\a)\) and \(F(\a)\).  Now Abel wants to define these functions for imaginary numbers of the form \(i\beta\).  

For this he substitutes formally \(iy\) for \(x\) in (\ref{eqn:EID}) and, integrates the integrand of the elliptic integral in (\ref{eqn:EID}) on the imaginary axis from \(0\) to \(iy\), and  obtains 
\[
i\int^y_0\frac{dt}{\sqrt{1+c^2t^2}\sqrt{1-e^2t^2}},
\]
where we see that the roles of \(e\) and \(c\) have been interchanged. Let
\[
\beta(y) := \int^y_0\frac{dt}{\sqrt{1+c^2t^2}\sqrt{1-e^2t^2}},
\]
which is again a monotone increasing function on the interval \([-\frac{\tilde{\omega}}{2},\frac{\tilde{\omega}}{2}]\),
where
\[
\frac{\tilde{\omega}}{2} := \int^\frac{1}{e}_0 \frac{dx}{\sqrt{1+c^2x^2}\sqrt{1-e^2x^2}},
\]
and we let the inverse of \(\beta(y)\) on this interval be denoted by \(y(\beta)\).  

We have already defined \(\phi(\a)\) to be \(x(\a)\) on \([-\omega/2,\omega/2]\), and now we define similarly \(\phi(i\beta):= iy(\beta)\) on the interval \([-i\frac{\tilde{\omega}}{2},i\frac{\tilde{\omega}}{2}]\) on the imaginary axis.  We then define on this same interval
\[
f(i\beta) := F(\beta),\;\textrm{and}\; F(i\beta)=f(\beta),
\]
using the interchange of \(c\) and \(e\) in the definition of \(\a(x)\) and \(\beta(y)\).  We note that \(\phi(\pm\frac{\omega}{2}) = \pm\frac{1}{c}\) and \(\phi(\pm i\frac{\tilde\omega}{2}) = \pm i\frac{1}{e}\).

Thus, at this point \(\f(\a)\) and \(\f(i\beta)\) are defined for \(\omega/2 \le \a\le \omega/2\), and \(-\tilde{\omega}/2 \le \beta \le \tilde{\omega}/2\). The problems remains to define \(\f(\a)\) and \(\f(i\beta)\) for all \(\a\) and \(\beta\), and to then define \(\f(\a+i\beta)\) for all complex numbers \(\a+i\beta\). 

For both of these tasks Abel needs a tool, and that is a specific generalization of the usual addition formulas for sines and cosines.  Abel formulates these new addition formulas for the three functions \(\f(\a), f(\a), F(\a)\) as follows:
\bea
\label{eqn:AF1}
\f(\a+\beta) &=& \frac{\f(\a)f(\beta) +\f(\b)f(\beta)f(\a)F(\a)}{1+e^2c^2\f^2(\a)\f^2(\beta)},\\
\label{eqn:AF2}
f(\a+\beta) &=& \frac{f(\a)f(\beta) -c^2\f(\a)\f(\beta)F(\a)F(\beta)}{1+e^2c^2\f^2(\a)\f^2(\beta)},\\
\label{eqn:AF3}
F(\a+\beta) &=& \frac{F(\a)F(\beta) +e^2\f(\a)\f(\beta)f(\a)f(\beta)}{1+ e^2c^2\f^2(\a)\f^2(\beta)}.
\eea
We recall briefly the classical formulas for trigonometric functions (as one finds in Euler's {\it Introductio} \cite{euler1748} from 1748, for instance):
\bea
\sin(\a+\beta) &=& \sin(\a)\cos(\beta) + \cos(\a)\sin(\beta),\\
\cos(\a+\beta) &=& \cos(\a)\cos(\beta) - \sin(\a)\sin(\beta),\\
\tan(\a+\beta) &=& \frac{\tan(\a) + \tan(\beta)}{1-\tan(\a)\tan(\beta)},
\eea
which has the same type of rational expressions as in (\ref{eqn:AF1}), (\ref{eqn:AF2}), and (\ref{eqn:AF3}).

Abel points out that these addition formulas follow from Legendre's theory of elliptic integrals \cite{legendre1811}, which follows up on the Euler addition theorem for elliptic integrals that we discussed earlier. He also gives a simple and elegant proof which we can sketch here (the same proof will work for the trigonometric formulas listed above as well). First, using the fact that
\bean
f^2(\a)&=& 1-c^2\f^2(\a),\\
F^2(\a)&=& 1 + e^2\f^2(\a),
\eean
then, by differentiating, we obtain
\bea
\label{eqn:*1}
f(\a)f'(\a) &= & -c^2\f(a)\f'(\a),\\
\label{eqn:*2}
F(\a)F'(\a) &=& 1 + e^2\f(\a)\f'(\a),
\eea
and from (\ref{eqn:EID}) we have
\be
\label{eqn:*3}
\f'(\a) = \sqrt{1-c^2\f^2(\a)}\sqrt{1+e^2\f^2(\a)}= f(\a)F(\a).
\ee
Substituting (\ref{eqn:*3}) in (\ref{eqn:*1}) and (\ref{eqn:*2}), we find that
\bean
\f'(\a) &=& f(\a)F(\a),\\
f'(\a) &=& -c^2\f(\a)F(\a),\\
F'(\a) &=& c^2\f(\a)f(\a),
\eean
the elliptic function analogue to \((\sin(\a))' = \cos(\a)\), etc. 

Now for the proof of, for instance, (\ref{eqn:AF1}), we denote the right-hand side of (\ref{eqn:AF1}) by \(r(\a,\beta)\), and compute both \(\frac{\partial r}{\partial \a}\) and \(\frac{\partial r}{\partial \beta}\) using the differentiation formulas above.  It turns out that \(\a\) and \(\beta\) appear symmetrically in these expressions and that one verifies by inspection that
\be
\label{eqn:PDE}
\frac{\partial r}{\partial \a} = \frac{\partial r}{\partial \beta}.
\ee
As was  known at the time, a solution of the partial differential equation (\ref{eqn:PDE}) is a function of the sum \(\a+\beta\), and hence there is a function \(\psi\) of one variable such that
\[
r(\a,\beta) = \psi(\a+\beta).
\]
One can find \(\psi\) by looking at particular values, and for instance, for \(\beta=0\), we have \(\f(0) = 0, f(0)=1, F(0)=1\), and hence 
\[
r(\a,0)=\f(a)=\psi(\a),
\]
and hence
\[
r(\a,\beta)= \psi(\a+\beta) = \f(\a+\beta),
\]
and (\ref{eqn:AF1}) is proved.  The addition formulas (\ref{eqn:AF2}) and (\ref{eqn:AF3}) can be proved in the same manner%
\footnote{Of course this proof depends on {\it knowing} what the right hand side of such an addition formula looks like, and this knowledge stems from the work of Euler and Legendre.}. 
Abel uses these addition formulas to define in a natural manner the evaluation of the elliptic functions on the real line for \(|\a|>\omega/2|\) and on the imaginary axis for \(|i\beta|> \tilde{\omega}/2\).  Then he also invokes the addition formula to define
 for instance 
\bean
\f(\a+i\beta) &=& \frac{\f(\a)f(i\beta) +\f(i\beta)f(i\beta)f(\a)F(\a)}{1+e^2c^2\f^2(\a)\f^2(i\beta)},\\
 &=& \frac{\f(\a)F(\beta) -i\f(\beta)F(\beta)f(\a)F(\a)}{1-e^2c^2\f^2(\a)\f^2(\beta)},
\eean
and similarly for the other two elliptic functions \(f\) and \(F\).

\begin{figure}
\vspace{6pt}
\centerline{
	\includegraphics[width=15cm]{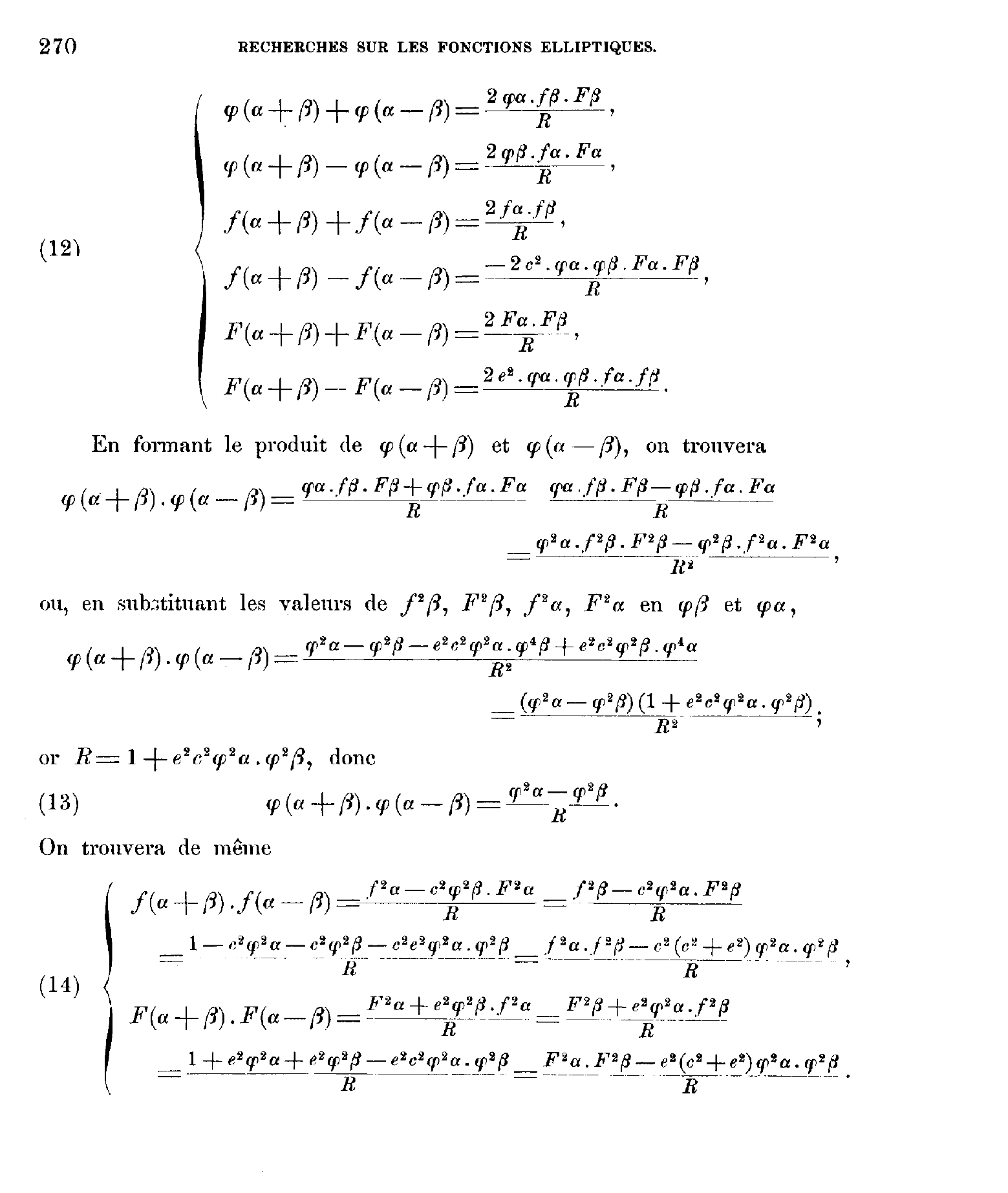}}
 \caption{Page 270 of Abel's paper on elliptic functions \cite{abel1828}}
 \label{fig:p270abel}
\end{figure}

After having used the addition formulas in this manner, then Abel says on p. 279 of \cite{abel1828}
\begin{quote}
Des formules (\ref{eqn:AF1}), (\ref{eqn:AF2}), (\ref{eqn:AF3}) on peut deduire une foule d'autres.
\footnote{From the formulas (\ref{eqn:AF1}), (\ref{eqn:AF2}), (\ref{eqn:AF3}) one can deduce a crowd of others.}
\end{quote}
In Figure \ref{fig:p270abel} 
we see a sample of his plethora of formulas that he derives from the basic addition theorems.  Here he has used the abbreviation 
\[
R= 1+e^2c^2\f^2(\a)\f^2(\beta).
\]
\begin{figure}
\vspace{6pt}
\centerline{
	\includegraphics[width=15cm]{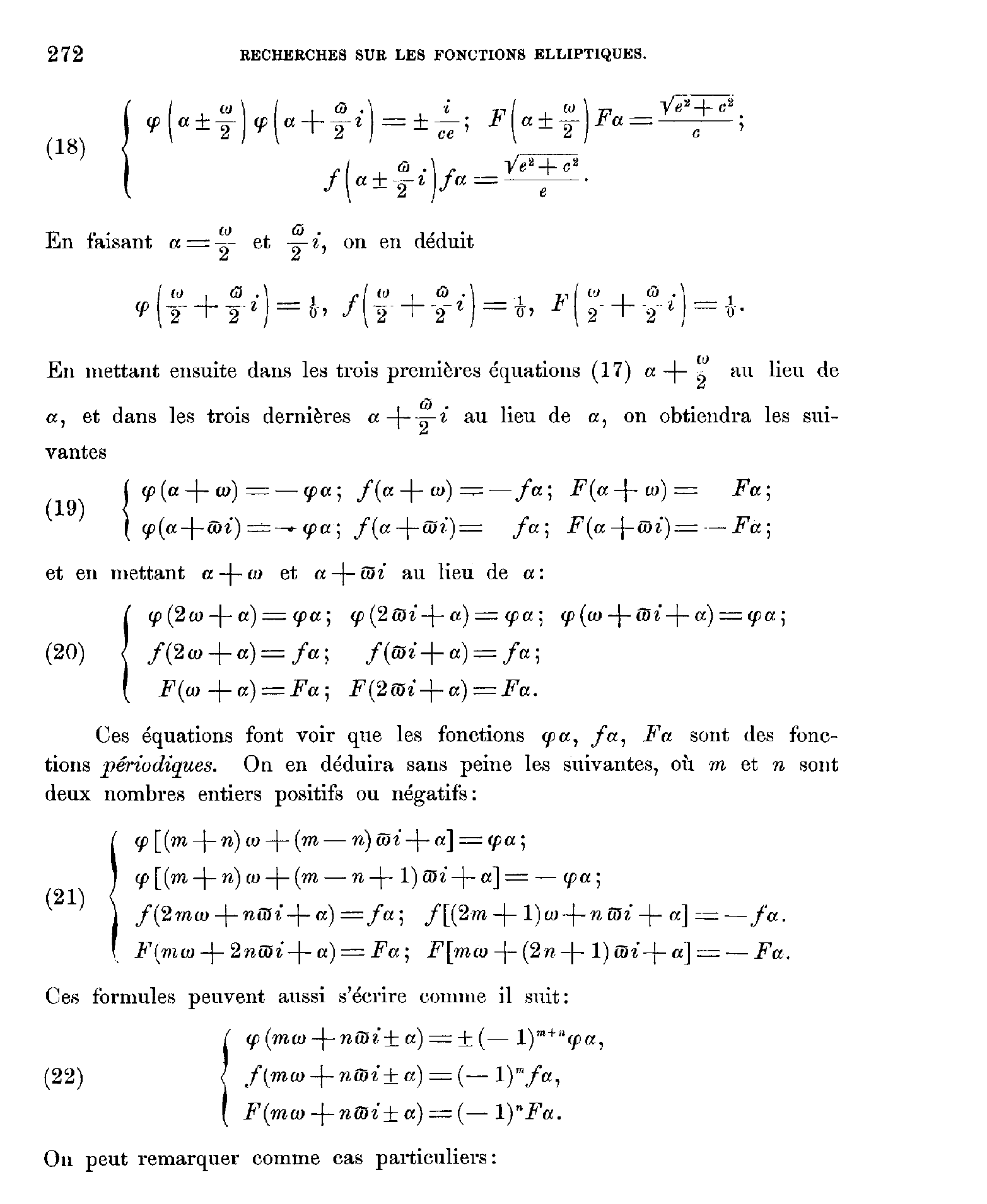}}
 \caption{Page 272 of Abel's paper on elliptic functions \cite{abel1828}}
 \label{fig:p272abel}
\end{figure}
After two more  pages of calculations we find on  p.272 of his paper (reproduced in Figure \ref{fig:p272abel}) the first formulation of the doubly-periodic nature of his elliptic functions.  This is equation no. 20 on this page in Figure \ref{fig:p272abel}.  At the top of the same page we see in the second equation that these elliptic functions all have a pole at the point \((\frac{\omega}{2},i\frac{\tilde{\omega}}{2})\) (and at the suitable translates of this point as well). This the first instance in the literature of a doubly-periodic function of a single complex variable.  

What is significant for us here is that one cannot formulate this notion of doubly-periodicity without the use of complex variables and in the decades that followed, these functions and others related to them, became important objects of study of meromorphic functions in the complex plane. Later in his paper Abel found many different kinds of representations of these functions.  An important historical point is that these functions played a role in applied mathematics as well.

In the remainder of his paper \cite{abel1828} Abel goes on to establish a variety of identities and properties for the elliptic functions he created in this paper, along with applications to the transformations of elliptic integrals and to the special case of the elliptic integral
\[
\int \frac{d}{\sqrt{1-x^4}}
\]
that Euler had studied, which describes the arc length of a leminiscate (\cite{abel1828}, pp. 361-362).  In addition he obtains a variety of representations of the elliptic functions in terms of infinite series and infinite products.

As we mentioned earlier, Jacobi had announced his discovery of elliptic functions in a short paper in December of 1827 \cite{jacobi1827} and followed up with foundational 190 page paper \cite{jacobi1829}. Interestingly, both of these papers were published in the {\it Astronomische Nachrichten}, edited by Heinrich Christian Schumacher (1780-1850), an important astronomer at that time.  Applications to astronomy of this new theory seemed to have been an important motivation for Jacobi at the time.

We will look at some of the innovations in Jacobi's paper \cite{jacobi1829}.  First, he proceeds in a similar manner to what Abel did at roughly the same time. Namely, he considers the inverse of the elliptic integral
\be
\label{EID2}
u(x) = \int_0^x \frac{dx}{\sqrt{(1-x^2)(1-k^2x^2)}},
\ee
to be
\bean
\f&=& \am u,\\
x&=& \sin \am  u.
\eean
Jacobi defines
\[
K :=\int_0^\frac{\pi}{2} \frac{d\f}{\sqrt{1-k^2\sin^2\f}},
\]
using the substitution \(x=\sin \f\). He then defines a number of other functions related to \(\sin \textrm{am} (u)\), which have now become standard in the theory of elliptic functions. 	

We discuss the most important ones here briefly.  Namely, we have the two additional functions \(\cos \am u\) and 
\[
\D \am u := \frac{d\am u}{du}= \sqrt{1-k^2\sin \am^2 u}.
\]
This was the notation of Jacobi, and towards the end of the 19th century it has become standard to write
\bean
\sn u&:= &\sin\am u,\\
\cn u &:=&\cos\am u,\\
\dn u&:=& \D\am u,
\eean
 for these three functions, which are the analogues of the three elliptic functions of Abel, \(\f, f\), and \(F\). These satisfy the properties
\bean
\sn^2 u+\cn^2 u &=&1,\\
\sn^2 u+ k^2\dn^2 u&=&1,
\eean
the analogues to \(\sin^2 x+\cos^2 x=1\) in this context, and they satisfy addition formulas, which are formulated explicitly by Jacobi for these three functions and other related functions (see \cite{whittaker-watson1962} or \cite{hurwitz-courant1964} for  proofs of these addition formulas). Jacobi does not prove these formulas, but depends on the earlier work of Legendre on elliptic integrals \cite{legendre1811}  for proofs in this elliptic functions context.  

Jacobi then defines \(\sn iv,\), \(\cn iv\), and \(\dn iv\), in the same manner as Abel, and using the addition formulas extends his elliptic functions to be functions of a complex variable, e.g., \(\sn (u+iv)\).  These functions are doubly-periodic, which follows easily from the addition formulas.  For instance, letting \(K'\) be defined by 
\[
K^2+(K')^2=1,
\]
one finds that
\[
\sn(u+iv+ 4K)=\sn(u+iv), \:\textrm{and}\;\sn(u+iv +i2K')=\sn(u+iv),
\]
which shows that \(\sn(u+iv)\) has two independent periods, \(4K\) along the real axis, and \(i2K'\) along the imaginary axis. One can find a complete set of these period relations for all of the Jacobi elliptic functions  in \cite{whittaker-watson1962}.

In his paper Jacobi obtains an extensive set of properties for the elliptic functions, many of which are similar to those derived by Abel (representation in terms of power series, infinite products, solutions of certain differential equations, etc.). Then on p. 186 of \cite{jacobi1829} he defines for the first time a new concept, which gives a new method of representing the Jacobi elliptic functions, and which becomes intrinsically very important in mathematics, independent of the theory of elliptic functions. This is Jacobi's discovery of {\it theta functions}, as they have been called ever since the time of Jacobi.  A theta function is a rapidly convergent Fourier series with quasi-periodic properties, and  the quotient of two such functions can represent an elliptic function.

Let us give  an example of two such functions whose quotient is an elliptic function.  Our notation differs from that used by Jacobi, but it is the same thing mathematically. Let 
\[
\theta(z;\tau):= \sum_{n=-\infty}^{\infty} e^{\pi in^2\tau + 2\pi i z}
\]
be defined for \(z\in \BC, \tau\in\BC\) with \(\Im(\tau)>0.\).  We consider \(\theta(z;\tau)\) as a function of \(z\), with \(\tau\) as a parameter. Since \(\Im \tau>0\), it follows that
\[
|e^{\pi in^2\tau}|\le e^{-\pi n^2 \Im \tau},
\]
which shows that, for fixed \(\tau\), the coefficients of the Fourier series converge to zero very rapidly, and hence \(\theta(z;\tau\)) is a holomorphic function of \(z\) for fixed \(\tau\). Moreover, it is clear from the definition that
\[
\theta(z+m;\tau) = \theta(z;\tau),
\]
so \(\theta(z;\tau)\) is periodic with period 1. Now we compute the behavior with respect to multiples of \(\tau\), i.e., we want to calculate \(\theta(z+k\tau;\tau)\). We find
\bean
\theta(z +k\tau;\tau) &=&\sum_{n=-\infty}^{\infty} e^{i\pi n^2 \tau +2\pi i n(z+k\tau)},\\
&=& e^{-i \pi \tau k^2}\sum_{n=-\infty}^{\infty} e^{i\pi\tau(n^2+2nk+k^2)+2\pi i n z},
\eean
and letting \(l=n+k\), we have
\bea
\theta(z+k\tau;\tau)&=& e^{-i\pi\tau k^2-2\pi ik z}\sum_{l=-\infty}^{\infty} e^{i\pi\tau l^2 + 2\pi ilz},\nonumber\\
\label{eqn:mult-factor}&=& e^{-i\pi\tau k^2-2\pi ik z}\theta(z;\tau).
\eea
This is the quasiperiodicity alluded to above. Except for the factor \(e^{-i\pi\tau k^2-2\pi k z}\), \(\theta(z;\tau)\) seems to be periodic in the direction \(\tau\). How can we exploit this? Let's consider a second such function
\[
\theta_1(z;\tau):= \theta(z+\frac{1}{2},m\tau).
\]
This is also holomorphic and periodic with period 1 in \(z\).  What is the periodicity in the direction \(\tau\)? Again we compute and find
\bean
\theta_1(z+k\tau;\tau)&=&\sum_{n=-\infty}^{\infty} e^{\pi i n^2+2\pi i n(z+\frac{1}{2}+k\tau}),\\
&=& e^{-\pi i k^2}\sum_{n=-\infty}^{n\infty}e^{\pi i \tau(n^2+2nk+k^2)+2\pi i n(z+\frac{1}{2})},
\eean
which gives, letting \(n= l-k\), as before
\bean
\theta_1(z+k\tau;\tau)&=& e^{-\pi i \tau k^2-2\pi i k z}\dot e^{-\pi ik}\theta_1(z;\tau),\\
&=& (-1)^k e^{-\pi i \tau k^2-2\pi ikz}\theta_1(z;\tau).
\eean

Thus the multiplicative factor here is the same as in (\ref{eqn:mult-factor}), except for the factor of \((-1)^k\). Therefore, if we form the quotient,
\[
e(z;\tau):= \frac{\theta(z,\tau)}{\theta_1(z;\tau)},
\]
we see that 
\[
e(z+m+k2\tau;\tau)=e(z;\tau).
\]
Thus we see that \(e(z;\tau)\) is a doubly-periodic function with periods \((1,2\tau)\), where \(\Im \tau>0\). By modifying suitably the choice of such theta functions, one can construct all of the Jacobi elliptic functions (again see either \cite{whittaker-watson1962}, or \cite{hurwitz-courant1964} or any other standard reference on elliptic functions).

\section{Holomorphic functions and mappings}
\label{sec:holomorphic}
So far in this paper we  have seen the development of the complex plane, Abel's theorems and the  creation of the theory of elliptic functions of a complex variable. Now we turn to a set of ideas which also started in the early decades of the 19th century, and which had very important developments throughout the course of this century, and this was the creation of the fundamental concepts of what we call today complex analysis, or function theory, as it was often called in the 19th century.  The fundamental concept is the study of special classes of complex-valued functions of a complex variable which are known today as {\it holomorphic} and {\it meromorphic functions}. We will see how these concepts arose out of the work of various mathematicians over a long period of time.

The fundamental innovators in the creation of function theory were Cauchy, Riemann and Karl Weierstrass (1815--1897), and we will discuss their respective contributions in some detail below.  Today a course in complex analysis is considered an essential part of undergraduate education, and over the course of the 20th century (and indeed towards the end of the 19th century) a number of texts evolved to explain this important subject, for instance, Hurwitz and Courant from 1922 (\cite{hurwitz-courant1964}), but often updated (for instance in 1964), and still in print, the classic text by Ahlfors \cite{ahlfors1953}, which many American undergraduates learned from, and there are many other fine more recent texts on the subject.

We start with the fundamental contributions of Cauchy, who contributed to the development of complex analysis throughout most of his very productive career.  His collected works consist of two series, each with about 12 volumes and approximately 500 pages per volume; this includes his published papers as well as a number of monographs and textbooks.  He worked on numerous fields of mathematics, including differential geometry, number theory, mathematical physics, and a variety of other areas. The first paper in complex analysis \cite{cauchy1814} was presented to the Academie des Sciences in 1814 and finally published in 1827.  Several footnotes added to the published version indicate some conceptual progress he had made in going from the real to the complex setting. 

In this paper \cite{cauchy1814} Cauchy considers a function \(f\) of a {\it real} variable \(z\), and shows that if \(z\) is considered to be a function of two other real variables \(x\) and \(y\),%
\footnote{We've used the now standard notation \(z, x\), and \(y\) for these variables, where \(z=x+iy\); Cauchy used \(y\), \(z\) and \(x\) with \(y=z+\sqrt{-1}x\) in this paper.  In his later papers he used the now standard notation, and in his earlier papers he used \(\sqrt{-1}\) instead of the symbol \(i\) for the imaginary unit, which he also used later.}
then 
\be
\label{eqn:c-1}
\pd{}{x}\left(f(z)\pd{z}{y}\right)= \pd{}{y}\left(f(z)\pd{y}{x}\right),
\ee
which is easy to verify. Namely
\bea
\label{eqn:c-2}
\pd{}{x}\left(f(z)\pd{z}{y}\right)&=& f'(z)\pd{z}{x}\pd{z}{y} +f(z)\frac{\partial^2z}{\partial x\partial y},\\
\label{eqn:c-3}
\pd{}{y}\left(f(z)\pd{y}{x}\right) &=& f'(z)\pd{z}{y}\pd{z}{x} + f(z)\frac{\partial^2z}{\partial y\partial x},
\eea
and since
\[
\frac{\partial^2z}{\partial x\partial y} = \frac{\partial^2z}{\partial y\partial x},
\]
we see that (\ref{eqn:c-1}) is satisfied.

Now we let \(z\) be a {\it particular} function of the two real variables \(x\) and \(y\) using complex numbers, namely let
\[
z=x+iy,
\]
where \(i=\sqrt{-1},\) and let \(f(z\)) take on complex values, where we let
\[
f= u+iv,
\]
where \(u\) and \(v\) are real-valued functions.
Then (\ref{eqn:c-1}) becomes, noting that \(\pd{z}{x}=1,\) and \(\pd{z}{y}=i\),
\be
\label{eqn:cr-1a}
i\pd{f}{x} (z)=1\pd{f}{y}(z),
\ee
that is,
\[
i\pd{}{x} (u+iv)=1\pd{}{y}(u+iv),
\]
which becomes, upon setting real and imaginary parts equal to each other,
\bea
\label{eqn:cr-1}
\pd{u}{v}&=&\pd{v}{y},\\
\label{eqn:cr-2}
\pd{v}{x}&=& -\pd{u}{y},
\eea
and this is the first appearance in Cauchy's work of the well known {\it Cauchy-Riemann equations}.%
\footnote{These equations had appeared earlier in the work of d'Alembert in the context of fluid dynamics and in the work of Euler  and Laplace for the evaluation of certain definite integrals.  See Kline \cite{kline1972}, pp. 626-628 for a discussion of this point. This is the beginning of his very well written historical chapter on the history of function theory.}
Cauchy remarks at this point in his paper (p. 338):
\begin{quote}
Ces deux \'{e}quations renferment toute la th\'{e}orie du passage du r\'{e}el \`{a} l'imaginaire, et il ne nous reste plus qu'\`{a} indiquer la mani\`{e}re de s'en servir.%
\footnote{These two equations contain all the theory of passing from the real to the imaginary, and it only remains for us to indicate how this can be used.}
\end{quote}
Thus Cauchy indicates that he understood the significance of these equations, and his work over the next 30 years certainly bears this out.  The implicit assumption that Cauchy make here is that the derivative \(f'(z)\) in (\ref{eqn:c-2}) and (\ref{eqn:c-3}) {\it exists}, as a generalization of \(f'(z)\), when \(z\) was a real variable.  This means that the limit
\[
f'(z)=\lim_{\e\rightarrow 0} \frac{f(z+\e)-f(z)}{z-\e},
\]
exists and is well defined, for small complex-valued \(\e.\)  We shall return to this point later when we look at Riemann's work.

\begin{figure}
\vspace{6pt}
\centerline{
	\includegraphics[width=10cm]{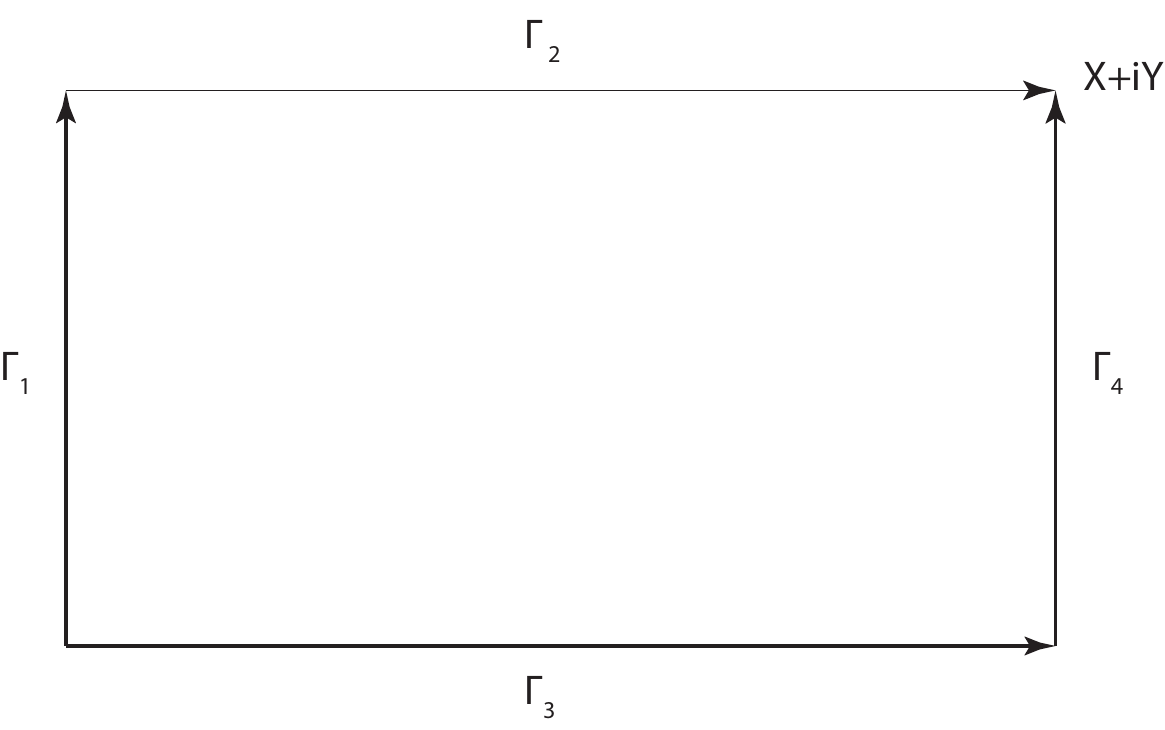}}
 \caption{A rectangle \(R\) in the complex plane}
 \label{fig:rectangle}
\end{figure}

The next step Cauchy takes is to integrate both sides of (\ref{eqn:cr-1a}) over a rectangle in \(\BR^2\), which we take to be the rectangle \(R\) defined as the product of the two intervals \([0,X]\) on the \(x\)-axis and \([0,Y]\) on the \(y\)-axis, as pictured in Figure \ref{fig:rectangle}.  Thus we have
\[
i\int_R \pd{f}{x} dxdy = \int_R\pd{f}{y} dxdy,
\]
and, assuming that these partial derivatives are continuous on the rectangle \(R\), we can evaluate these area integrals in terms of interated integrals, obtaining,
\[
i\ \int_0^Y\left(\int_0^X \pd{f}{x}dx\right)dy = \int_0^X\left(\int^Y_0 \pd{f}{y} dy\right)dx,
\]
which gives, using the fundamental theorem of calculus,
\be
\label{eqn:ci-1}
 i\int_0^Y[f(X,y)-f(0,y)]dy = \int_0^X[f(x,Y)-f(x,0]dy.
\ee
If we denote by \(\G_1+\G_2\) and \(\G_3+\G_4\) the two path along the edges of the rectangle from \(0\) to \(X+iY\) as indicated in Figure \ref{fig:rectangle},
then we see that (\ref{eqn:ci-1}) becomes
\[
\int_0^Xf(x,0)dx + \int^Y_0fX,y)d(iy) = \int_0^Yf(0,y)d(iy)+\int_0^Xf(x,Y)dx,
\]
which is the same as
\be
\label{eqn:ci-2}
\int_{\G_3+\G_4}f(z)dz = \int_{\G_1+\G_2} f(z)dz.
\ee
This equation says that the path integrals of \(f(z)\) in the complex plane along the two paths \(\G_1+\G_2\) and \(\G_3+\G_4\) have the same value.  If we let \(\G\) be the closed path \(\G_1+\G_2-\G_3-\G_4\) around the boundary of the rectangle \(R\), then we have
\be
\label{eqn:ci-3}
\int_\G f(z) dz.
\ee
This is the famous {\it Cauchy integral theorem} in this important special case.  

\begin{figure}
\vspace{6pt}
\centerline{
	\includegraphics[width=10cm]{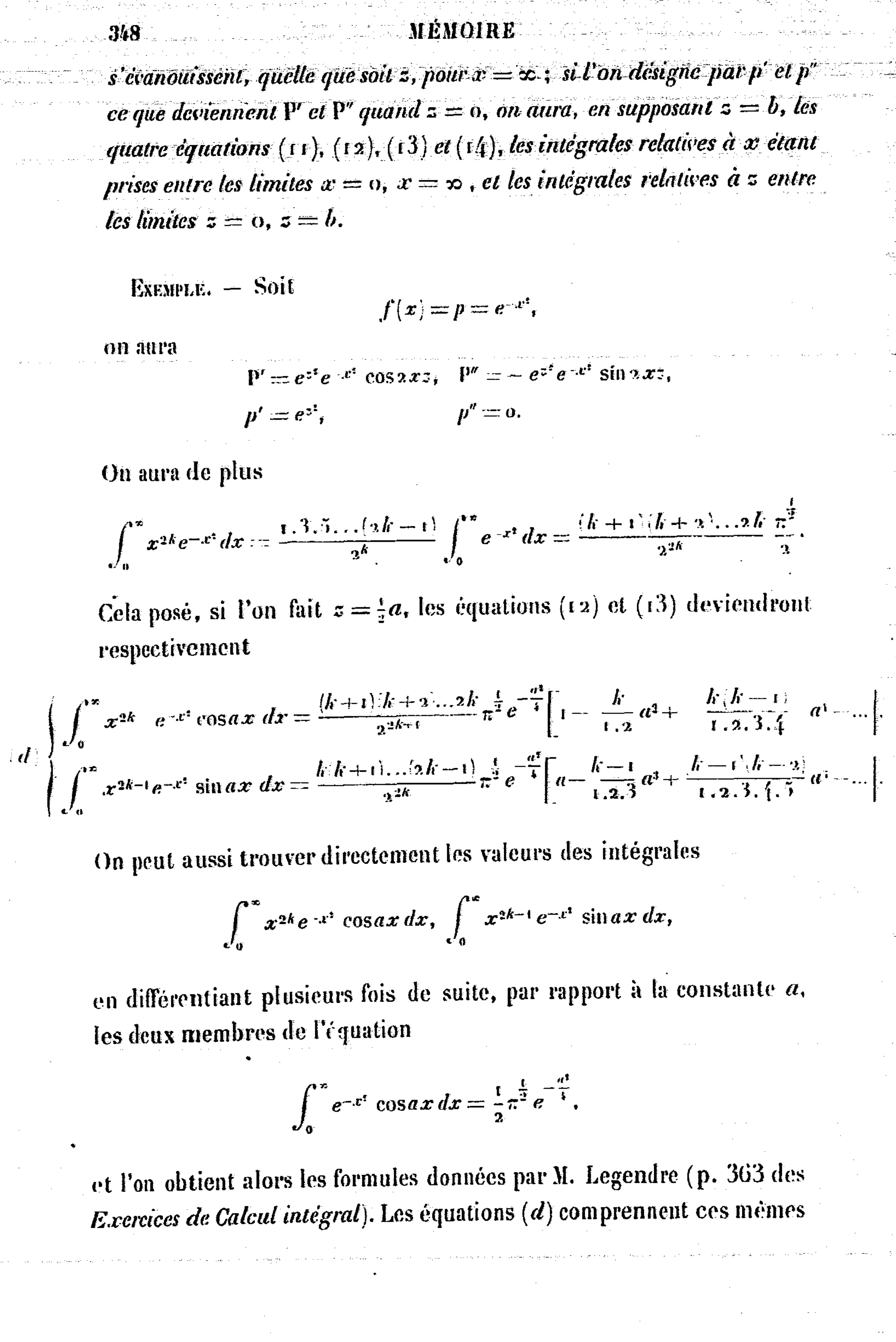}}
 \caption{Page 348 of Cauchy's paper \cite{cauchy1814} showing the evaluation of definite integrals using the first version of Cauchy's Integral Theorem}
 \label{fig:p348-Cauchy1814}
\end{figure}

Cauchy expressed this theorem in terms of the real-valued functions \(u\) and \(v\), and only later, when this 1814 paper was printed in 1827 he added footnotes indicating how using the complex variables notation the work could be simplified, as we have done here. He used these results to compute various examples of definite integrals, usually of the form \(\int_{-\infty}^\infty f(z)dz\), where, for instance, the vertical integrals vanished asymptotically, and one was left with something like
\[
\int_{-\infty}^\infty f(x+ib)dx = \int_{-\infty}^\infty f(x)dx,
\]
the integration being shifted from the \(x\)-axis to a translate of the \(x\)-axis in the complex plane, which could often be simpler to compute.  He also concerned himself with a variety of singular integrals, and proper values of integrals. In Figure \ref{fig:p348-Cauchy1814} we see a sampling of the evaluation of such integrals.  

Cauchy wrote a number of papers and books on this topic over the decades following his seminal 1814 paper, and his writings became the basis for a significant amount of what became known as function theory.  We will discuss briefly some of these in the next several paragraphs.  Perhaps his most important paper is a fundamental paper written in 1825 entitled {\it M\'{e}moire sur les int\'{e}grales, d\'{e}finies entre les limites imaginaires} \cite{cauchy1825}, which, due to its significance, was reprinted in 1874.  Curiously, it doesn't appear in his collected works.  We recall that at this time the mathematicians still used the term ``imaginary number" for what we now call complex numbers. 

Cauchy considered in \cite{cauchy1825} complex-valued functions of a complex variable to have a well-defined derivative at each point where it had a finite value, and at any point where the function became infinite, he considered the function to be locally the reciprocal of a function with a zero of finite order. Today we call such a function a {\it meromorphic function}. For such a function \(f(z)\), he defined the path integral 
\[
\int_{z_0}^{Z} f(z)dz,
\]
where 
\(\f(t)\) is a smooth curve in the complex plane parametrized by a parameter \(t\) varying between \(t_0\) and \(T\), and where \(z_0=\f(t_0),\) and \(Z=\f(T)\), to be 
\[
\int_{t_0}^T f(\f(t))\f'(t)dt,
\]
provided the function is finite at all points of the path.  

Then Cauchy shows by a calculus of variation technique that if one perturbs the path suitably, then the first variation of the perturbation vanishes, indicating to Cauchy that the path integral for the perturbed paths are the same as for the original unperturbed path, an infinitesimal version, so to speak, of the Cauchy Integral Theorem in this case. Then Cauchy considers the case where the path encounters a point where \(f(z)\) becomes infinite, and he introduces the notion of {\it proper value} for such a singular integral, being a limit of a specific perturbation of the integral.  

Cauchy defines, for a function which is infinite at a point \(z_0\) of finite order, which means that \(f(z)\) can be represented near \(z_0\) in the form
\be
\label{eqn:residue}
f(z) = \frac{a_{-m}}{(z-z_0)^m} + ... + \frac{a_{-1}}{(z-z_0)} + g(z),
\ee
where \(g(z)\) is finite at \(z_0\), the {\it residue of \(f(z)\) at \(z_0\)} to be the coefficient \(a_{-1}\) in (\ref{eqn:residue}). Let us denote  this residue by \(\textrm{Res}f(z_0)\) (a notation similar to that which Cauchy uses in his later papers).

He formulates and proves, using various perturbations of path integrals,  the  {\it Cauchy Residue Theorem} for a large rectangle \(R\) containing a finite number of singular points of \(f(z)\), \(z_0, z_1, ... , z_N,\) to be
\[
\int_{\partial R} f(z)dz = 2\pi i\sum_{k=0}^N\textrm{Res}f(z_k).
\]
He uses this for numerous examples of calculations of specific definite integrals, similar to what he had done in his 1814 paper.

This was a very important breakthrough for complex analysis, and he further developed this theory with numerous example in his four volume work {\it Exercices de Math\'{e}matiques}  \cite{cauchy1826}, published between 1826 and 1829.  In a slightly later work {\it Exercices d'Analyse et de Physique Math\'{e}matique} \cite{cauchy1841} he developed the {\it Cauchy Integral Formula} as 
\[
f(z_0) = \frac{1}{2\pi i}\int_{\partial \D}\frac{f(z)dz}{z-z_0},
\]
where \(\D\) is a small circular disc centered at \(z_0\), and where \(f(z)\) has finite values on the closure of the disc.  In addition he created similar integral formulas for all of the coefficients of the Laurent series of a function \(f(z)\) with an isolated singularity at \(z_0\). Such a series was first formulated for an isolated singularity by Pierre Alphonse Laurent  (1813--1854) in 1843  \cite{laurent1843}, and  was developed in full in a paper published   posthumously in 1863 \cite{laurent1863}. The series has the form
\[
\sum_{n=-\infty}^\infty a_n(z-z_0)^n,
\]
and the coefficients can be computed by
\[
a_n= \frac{1}{2\pi i} \int_{\partial \D}\frac{f(z)dz}{(z-z_0)^{n+1}},
\]
which became an alternative to the usual formula for the Taylor series coefficients \(\frac{f^{(n)}(z_0)}{n!}\) in terms of derivatives, which don't make any sense when \(f(z)\) is singular at \(z_0\).

In 1846 \cite{cauchy1846a} Cauchy formulated but did not prove what has become known as {\it Green's Theorem}%
\footnote{George Green (1793-1841) in \cite{green1828} formulated and proved a version of Green's theorem in three dimensions, now known as the "divergence theorem" or often Gauss's Theorem, all special cases of the general Stokes' Theorem in \(n\)-dimensions, see e.g., \cite{spivak1999-1}},
namely, for two  continuously differentiable real-valued functions \(P(x,y),\) and \(Q(x,y)\) defined on the closure of a bounded domain \(D\)  in \(\BR^2\),
\be
\label{eqn:Green}
\int_D \left(\pd{Q}{x} - \pd{P}{y}\right)dxdy = \int_{\partial D} Pdx +Qdy.
\ee
He showed that this formula will imply a proof of the Cauchy Integral Theorem, and, as it is so simple and instructive, we indicate its proof here.  Let \(f(z)\) be a function which satisfies the Cauchy-Riemann equations on the closure of a domain \(D\), then we see that
\bean
\int_{\partial D} f(z)dz&=& \int_{\partial D} (u+iv)(dx+idy),\\
&=& \int_{\bd{D}} (udx-vdy) + i\int_{\bd{D}} (vdx+udy),\\
&=&\int_D\left(\pd{v}{x}+\pd{u}{y}\right)dxdy + i \int_D \left(\pd{u}{x} -\pd{v}{y}\right)dxdy,\\
&=& 0+i0,
\eean
by Green's Theorem (\ref{eqn:Green}) above. This shows the direct relationship between a function satisfying the Cauchy-Riemann equations and the Cauchy Integral Theorem.  Riemann proved Green's Theorem, and hence Cauchy's Integral Theorem,  in his dissertation \cite{riemann1851}.

Finally we mention that Cauchy wrote several papers dealing with the generalizations of his ideas to the case of multivalued functions, including applications to the study of elliptic integrals and functions (see, e.g.,  \cite{cauchy1846c}) and in \cite{cauchy1851} he extended his theory of multivalued function to be single-valued functions spread over a complex plane, a concept that Riemann developed much more fully, as we will see below. Cauchy introduced the notion of branch points and branch cuts, which again Riemann would develop more fully later. This paper of Cauchy is one of several of his dealing with path integrals of multivalued functions.

We now turn to Riemann's Inaugural Dissertation from 1851 \cite{riemann1851}, and we want to discuss it   in the context of the development of holomorphic functions and function theory in general.   Riemann starts his paper with the introduction of surfaces spread over domains in the complex plane (branched coverings), and elementary notions of connectedness for open sets of such surfaces.  This was followed up in his 1857 paper on Abelian functions \cite{riemann1857}, and this became the theory of Riemann surfaces and eventually developed into the theory of complex manifolds in general. 

He starts, as did Cauchy, by considering the class of complex-valued functions of a complex variable on an open set that have a well-defined derivative at each point.  He computes the derivative of a function \(f(z)=u(z)+iv(z)\) as
\be
\label{eqn:complex-derivative}
\frac{d(u+iv)}{dx+iy}=\frac{\left(\pd{u}{x}+i\pd{v}{x}\right)dx+i\left(\pd{v}{y}-i\pd{u}{y}\right)dy}{dx+idy},
\ee
and he argues that this is a well-defined complex number if and only if
\be
\label{eqn:CR2}
\pd{u}{x}=\pd{v}{y}, \;\textrm{and}\; \pd{v}{x}= -\pd{u}{y}.
\ee
This gives
\[
f'(z) =\frac{df}{dz}=\frac{d(u+iv)}{dx+iy}=\frac{1}{2}\left(\pd{u}{x} +\pd{v}{y}\right)+ i \frac{1}{2}\left(\pd{v}{x}-\pd{u}{y}\right).
\]
as the value of the derivative at the point \(z\) in terms of the real-valued partial derivatives, when this derivative exists.  

Riemann then shows that a function \(f(z)\) satisfying (\ref{eqn:CR2}) is conformal and orientation-preserving at each point of the domain where \(f'(z)\ne0\). Conformality means that if two smooth curves meet at a point \(z_0\), then the angle between their tangents at that point  is preserved under the mapping of the plane to the plane by the functions \(f(z)\) near the point \(z_0\).%
\footnote{Conformality was called by Gauss \cite{gauss1825} and by Riemann \cite{riemann1851} "im kleinsten Theilen \"{a}hnlich", or in English we might say "infinitesimally similar". We will use the now common term {\it conformal}.}
The mapping preserves orientation if the direction from one tangent line to another is also preserved.  Gauss had found necessary and sufficient conditions for a local mapping of \(\BR^2\) to be conformal in 1825 \cite{gauss1825}, which turned out, in the language of holomorphic mappings, a la Riemann in the paper we are discussing here, to mean a holomorphic mapping.%
\footnote{A conformal mapping in this same context, could also be an {\it anti-holomorphic mapping}, i.e., a mapping \(f(z\), such that \(\pd{f}{z}\) is zero, instead of the holomorphic mapping where \(\pd{f}{\overline z}\) vanishes, using the contemporary notation that \(\pd{}{z}= \frac{1}{2}\left(\pd{}{x}-i\pd{}{y}\right)\), and  \(\pd{}{\overline z}= \frac{1}{2}\left(\pd{}{x}+\pd{}{y}\right)\).}
That a holomorphic mapping at a point where \(f(z_0)\) is not zero preserves orientation which is, as Riemann points out, easy to prove, since the Jacobian determinant of the mapping at that point is given by
\[
\pd{(u,v)}{(x,y)} = \left|\ba{cc}
\pd{u}{x}&\pd{u}{y}\\
\pd{v}{x} & \pd{v}{y}
\ea\right|,
\]
which, by the Cauchy-Riemann equations (\ref{eqn:CR2}) has the value
\[
\left(\pd{u}{x}\right)^2 + \left(\pd{u}{y}\right)^2 = \left(\pd{v}{y}\right)^2 +\left(\pd{v}{x}\right)^2>0,
\]
unless \(f'(z_0)=0\). Riemann proves later in this paper that a holomorphic mapping from an open set in \(\BC\) to \(\BC\) has an image which is an open set.  Any mapping that maps open sets to open sets is called an {\it open mapping}, and this theorem of Riemann is often referred to in function theory as the {\it open mapping theorem}. In this context, one speaks of a holomorphic function defined on a domain where \(f'(z)\ne0\) as providing a {\it conformal mapping} from one domain to another. 

Riemann stated and gave a proof of Green's Theorem (\ref{eqn:Green}) for a domain in \(\BR^2\) and showed, as Cauchy had done that Cauchy's Integral Theorem was valid and that an integral of a holomorphic function in a simply-connected domain was path independent. He also used a version of Green's Theorem to develop some fundamental properties of harmonic functions, which play an important role in his later development in this same paper of what we now call the Riemann mapping theorem, which we will discuss shortly.  {\it Harmonic functions} are solutions of the equation 
\[
\D u= 0,
\]
where 
\[
\D:= \frac{\partial^2}{\partial x^2} + \frac{\partial^2}{\partial y^2}
\]
is the {\it Laplacian differential operator} in two dimensions. This terminology, which is now standard, was not used by Riemann and was introduced by William Thomson (Lord Kelvin, 1824-1907) in the mid 19th century (see Kline \cite{kline1972}, p. 685).

The main tool Riemann uses is a variation of what are now known as Green's Formulas which he proved as a consequence of Green's Theorem.  Namely, let \(X\) and \(Y\) be two continuously differentiable functions on a domain \(T\subset \BR^2\), which satisfies on  \(T\) the equation
\[
\pd{X}{x} + \pd{Y}{y} = 0,
\]
then. by Green's Theorem (\ref{eqn:Green}), the boundary integral
\[
\int_{\bd T}Ydx-Xdy=0.
\]
Riemann introduces normal and tangential coordinates along a neighborhood of the boundary curve \(\bd T\) in the following manner. He lets \(s\) be the arc length from a fixed point on the boundary to a variable point \(P\) on the boundary and let \(p\) be the distance from that point \(P\) along an inner directed normal to a point \(z=x+iy\) on the interior of \(T\).  Then, by letting \(\xi\) be the angle the normal at \(P\) makes with the \(x\)-axis, and \(\eta\) the angle the normal makes with the \(y\)-axis, the version of Green's Theorem that Riemann uses becomes
\[
\int_T \left(\pd{X}{x} +\pd{Y}{y}\right)dT = -\int_{\bd T} (X\cos \xi + Y \cos \eta)ds.
\]
Riemann computes the change of variables formulas
\[
\ba{ll}
\pd{x}{p} = \cos \xi,& \pd{y}{p} = \cos \eta,\\
\pd{x}{s} = \cos \eta,& \pd{y}{s} = -\cos \xi,
\ea
\]
using a positive orientation of the coordinates \((s,p)\) with respect to the standard orientation of the \((x,y)\) plane.  Green's Theorem then becomes
\bea
\int_T(\pd{X}{x}+ \pd{Y}{y})dt &=& -\int_{\bd T} \left(X\pd{x}{p} + Y\pd{y}{p}\right)ds,\\
\label{eqn:Green2}
	&=& \int_{\bd T}\left(X\pd{y}{s} - Y\pd{x}{s}\right)ds.
\eea

 If \(\pd{X}{x}+ \pd{Y}{y}=0\) in \(T\), then these boundary integrals are zero.  If \(\pd{X}{x}\) or \(\pd{Y}{y}\) have singular points at some finite set of points in \(T\), say \(z_1,...,z_N\), and 
\[
\pd{X}{x}+\pd{Y}{y} = 0 \; \textrm{on}\; T-\{z_1,...,z_N\},
\]
then
\be
\label{eqn:green-sing}
\int_{\bd T} (X\pd{x}{p}+Y\pd{y}{p})ds = -\sum_{j=1}^N \int_{\bd \D_j}(X\pd{x}{p} + Y\pd{y}{p})ds,
\ee
where \(\D_j\) are small nonintersecting discs centered at \(z_j\), such that  the closure of each disc is contained in the open set \(T\). This is somewhat parallel to Cauchy's Residue Theorem in this context, where Cauchy gave meaning to the localized integrals in terms of residues of the holomorphic functions at such singular points. 

Riemann then considers two functions \(u\) and \(\tilde u\) which are \(C^2\) on \(T\) and which have  a continuous extension along with their first derivatives to \(\bd T\).%
\footnote{Here we use the notation \(C^k\) to denote functions which have continuous derivatives of order \(k\) and \(C^\infty\) will mean continuously differentiable of all orders.} 
Now suppose that both \(u\) and \(\tilde u\) are harmonic on \(T\), then setting  
\bean
X&=& u\pd{\tilde u}{x}-\tilde u\pd{u}{x},\\
Y&=& u\pd{\tilde u}{y}-\tilde u\pd{u}{y},
\eean
it follows that
\[
\pd{X}{x}+\pd{Y}{y} = u\D \tilde u-\tilde u\D u,
\]
and using Green's Theorem (\ref{eqn:Green2}) above gives
\be
\label{eqn:green3}
\int_{\bd T}  \left(u\pd{u'}{p} - u'\pd{u}{p}\right)ds =0.,
\ee
 Riemann considers two particular cases for \(\tilde u\), which lead to interesting and useful results.  The first case is to simply take \(\tilde u \equiv 1,\) and it follows that 
\be
\label{eqn:green4}
\int_{\bd T} \pd{u}{p}ds = 0.
\ee
In the second case, for any particular point \(z_0\in T,\) Riemann chooses polar coordinates \(z-z_0=r^{i\f}\) and sets
\[
\tilde u(z):= \log r = \log |z-z_0|.
\]
Then using the extension of Green's theorem to the case of singularities  (\ref{eqn:green-sing}) it follows that
\[
\int_{\bd T} \left(u\pd{\log r}{p} -\log r \pd{u}{p}\right)ds = \int_{\D_{\e}}\left(u\pd{\log r}{p} -\log r \pd{u}{p}\right)ds,
\]
where \(\D_{\e}\) is a small disc of radius \(\e\), centered at \(z_0\).  Note that on the boundary of \(\D_{\e}\)
\[
\pd{\log r}{p} = -\pd{\log r}{r} = -\frac{1}{r},
\]
and thus
\bea
\int_{\bd T} \left(\log r\pd{u}{p} -u\pd{\log r}{p}\right)ds &=& \int_0^{2]pi}u(\e e^{i\f}df +\log r\int_{\D_{\e}} \pd{u}{p}ds,\nonumber\\
\label{eqn:green5}
&=&\int_0^{2\pi} u(\e e^{i\f}d\f,
\eea
since the second term on the right hand side vanishes by (\ref{eqn:green4}) (where we let \(T\) be \(\D_{\e}\)). Now letting \(\e \rightarrow 0\) in (\ref{eqn:green5}) Riemann obtains the following formula
\be
\label{eqn:harmonic-integral}
u(z_0) = \frac{1}{2\pi}\int_{\bd T} \left(\log|z-z_0|\pd{u}{p} - u \pd{\log |z-z_0|}{p}\right)ds,
\ee
which represents the value of the harmonic function \(u(z)\) at an interior point \(z_0\) of \(T\) in terms of the boundary integral on \(\bd T\).%
\footnote{In more contemporary literature, the normal derivative of data along \(\bd T\) is usually denoted by \(\pd{}{n}\).}
This result and the use of the potential function \(\log r\) is similar to the work of Green from 1828 \cite{green1828}, in which the three-dimensional potential function \(1/r\) is used in \(\BR^3\).

Suppose we restrict the harmonic function \(u(z)\) above to a disc \(\D_{\e}\) of radius \(\e\) centered at \(z_0\) whose closure is contained in the domain \(T\).  Then the formula (\ref{eqn:harmonic-integral}) becomes
\bea
\nonumber
u(z_0) &=& \frac{1}{2\pi} \log r \int_{\bd \D_{\e}} \pd{u}{p}ds + \frac{1}{2\pi}\int_0^{2\pi} u(\e e^{i\f}) d\f\\
\label{eqn:MVT}
&=&\frac{1}{2\pi} \int_0^{2\pi} u(\e e^{i\f})d\f,
\eea
since the first term on the right hand side vanishes, by (\ref{eqn:green4}), that is to say \(u(z_0)\) is the {\it mean value} of its values on the boundary of the disc \(\D_{\e}\).  Equation (\ref{eqn:MVT}) is known as the {\it Mean Value Theorem} for harmonic functions (and this is true in all dimensions).  It has numerous consequences, as Riemann shows in his paper, and we list some of them here.  We refer the reader to, e.g. Ahlfors \cite{ahlfors1953} for proofs of these results, as well as to this paper of Riemann.

Let \(u\) be a harmonic function in a domain \(T\).
\bi
\item 
{\it Removable Singularity Theorem}: If \(u\) is potentially singular or undefined at a point \(z_0\), and if \(\r\) is the distance from a neighboring point \(z\) to \(z_0\), and if
\[
\r \pd{u}{x}\; \textrm{and} \r \pd{u}{y} \rightarrow 0,\; \textrm{as}\; \r \rightarrow 0,
\]
then \(u\) can be continued as a continuous function to \(z_0\) and \(u\) is harmonic in a neighborhood of \(z_0\).
\item
{\it Smoothness}:
The harmonic function \(u\) is \(C^\infty\) in all of \(T\) (this follows from (\ref{eqn:harmonic-integral}) by differentiation under the integral sign).
\item
{\it Maximum Principle}: The harmonic function cannot have a local maximum or minimum at any interior point of \(T\) unless \(u\) is a constant function.
\item
{\it Identity Theorem}: The harmonic function \(u\) in \(T\) is determined by the values of \(u\) and \(\pd{u}{p}\) on any arc segment in \(T\), and moreover, if  on a segment of an arc in \(T\)  \(u\equiv 0\) and \(\pd{u}{p}\equiv 0,\), then \(u\equiv 0 \) in \(T\).
\ei

Riemann remarks that many of these properties of harmonic functions carry over to holomorphic functions in a natural manner. For instance, 
\bi
\item
{\it Riemann removable singularity theorem}: If \(f\) is holomorphic on a punctured disc centered at \(z_0\),  \(\D-\{z_0\}\), and if \((z-z_0)f(z)\rightarrow 0\),
as \(z\rightarrow z_0,\) then \(f\) extends as a holomorphic function to \(\D\).
\item
{\it Smoothness}:  A holomorphic function in a domain is infinitely differentiable.
\item
{\it Maximum principle}:  The modulus of  holomorphic function \(f\) in a domain can take on a local maximum in the interior of the domain only if \(f\) is constant in the domain.
\ei.

In Section 15 of his paper (p. 28), Riemann formulates and proves the {\it Open Mapping Theorem} for holomorphic functions, mentioned earlier.  Namely, let \(T\) be a domain in \(\BC\), and \(U\) be any open subset of \(T\), and let \(f:U\rightarrow \BC\) be a nonconstant holomorphic function defined on \(U\), then \(f(U)\), the image of \(f\) under the mapping \(f\), is an open set in \(\BC\).%
\footnote{At the time of Riemann, the notion of open set was not yet a mathematical concept.  He formulated his theorem in terms of neighborhoods of points.  We are giving  the modern formulation of this important result.}
This is a strong property of holomorphic functions, and it is also proved in any standard complex analysis text (again, e.g., Ahlfors \cite{ahlfors1953}). We note that, in contrast, this is not true for real-valued smooth (or real-analytic) functions.  As a simple example, the mapping \(f(x): \BR \rightarrow \BR\), given by \(f(x) = \sin x\), has as an image of the open set \(\BR\) the closed set \([-1,1]\). 

In the last section of his paper Riemann comes to what is undoubtedly the deepest result in this very innovative paper, the {\it Riemann Mapping Theorem}. We note first that in the beginning of the paper Riemann formulated the concept of what is now known as a Riemann surface spread over a region of the complex plane, a branched covering of an open set in \(\BC\), and most of his results in this paper are formulated in this more general context. We have chosen to formulate his results for domains in the complex plane, as that is simple.  In Sections 5 and 6 of Riemann's Dissertation he formulates and proves a number of elementary results concerning the Riemann surfaces he introduces, and these include the important notions of connectedness of domains in the plane, for instance, simply-connectedness, which we have already had occasion to use, and more general connectedness of order \(n\), and we will come back to these concepts  in the next section.  Simply-connected is defined by Riemann for a domain in the plane as having the property that any curve in the plane joining any two boundary points will split the domain into two parts that are not connected to each other (he always means path-wise connected).  

Now we can give, in his own words, Riemann's formulation of the Riemann Mapping Theorem:
\begin{quote}
Zwei gegebene einfach zusammenh\"{a}ngende ebene Fl\"{a}chen k\"{o}nnen
stets so auf einander bezogen werden, dass jedem Punkte der einen
Ein mit ihm stetig fortr\"{u}ckender Punkt der andern entspricht und
ihre entsprechenden kleinsten Theile ahnlich sind; und zwar kann zu
Einem innern Punkte und zu Einem Begrenzuugspunkte der entsprechende
beliebig gegeben werden; dadurch aber ist f\"{u}r alle Punkte die
Beziehung bestimmt.%
\footnote{Two given simply-connected plane surfaces can always be related to one another, so that each point of one corresponds in a continuous manner to each point of the other and such that the corresponding smallest parts are similar [infinitesimally similar; conformal]; and indeed such that a given inner point and a given boundary point correspond to a specified interior and boundary point; with this last condition, the relationship is determined for all points.}
\end{quote}
Riemann immediately reduces this formulation to the simpler statement that any simply-connected domain can be conformally mapped onto the unit disc, with one interior point mapping to the center of the disc, and a boundary point mapping to a specified boundary point of the unit disc (e.g., \(z=1\)).

The proof that Riemann gives is based on what he terms the {\it Dirichlet Principle}, a name he gave to this principle in his follow-up paper on Abelian functions  (which used the same principle for additional results) in 1857 \cite{riemann1857}. We will discuss this principle here, and then indicate how this became a problem for Riemann and his proof of the Riemann mapping theorem, and how this was finally resolved some 50 years later by David Hilbert (1862--1943). We will return to the Riemann mapping theorem in our Conclusion, Section \ref{sec:conclusion}.

We now formulate an important special case of the Dirichlet Principle.  Let \(T\) be a bounded domain in \(\BR^2\) with a smooth boundary, and let \(f\) be a continuous function on \(\bd T\).  Consider the family \(\SF\) of real-valued functions that are continuously differentiable in \(T\) and continuous on \(\overline T\), such that \(u|_{\bd T} = f.\)  This is an infinite-dimensional family of functions.  Consider the Dirichlet integral
\be
\label{eqn:dirichlet-integral}
D(u):= \int_T \left[\left(\pd{u}{x}\right)^2+\left(\pd{u}{y}\right)^2\right]dxdy, \;\textrm{for}\; u\in \SF.
\ee
From the definition it is clear that
\[
0\le D(u) < \infty.
\]
It follows that
\[
m:=\inf_{u\in\SF} D(u) 
\]
is well defined and \(m\ge 0\). The Dirichlet principle asserts that there exists a unique \(u\in \SF\) such that
\[
D(u) = m,
\]
and moreover, that \(u\) is harmonic (hence \(C^\infty\) on \(T\)), and that \(u\) is continuous on \(T\) and that \(u|_{\bd T} =f\). This last statement says that \(u\) is a solution to the Dirichlet problems (find a harmonic function with given boundary values on \(\bd T\)). 

Today, there are many different proofs of the solution to the Dirichlet problem, but in the mid-nineteenth century, these did not  exist, and mathematicians and physicists used this principle to solve many difficult problems.  Specifically, Green and Thomson formulated and used this principle, Green in 1835
 \cite{green1835}
 in the context of gravitational attraction, and Thomson in 1848 
\cite{thomson1848}
 as a general mathematical principle (called Thomson's Principle in England for some time, see Kline
 \cite{kline1972}, p. 685). 
 In lectures in G\"{o}ttingen in 1856 concerning inverse square forces, which very likely Riemann attended, Dirichlet used this minimization principle for the existence and uniqueness of specific harmonic functions, i.e., to solve the Dirichlet problem.  Riemann used this principle in both his 1851 dissertation \cite{riemann1851} and his seminal paper on Abelian functions \cite{riemann1857}, wherein he denoted this principle as the {\it Dirichlet Principle}, and it has been called that ever since (in spite of the earlier work of Green and Thomson).  

Then in 1870 Weierstrass gave a lecture at the Berlin Academy of Sciences (which was  published in his collected works in 1895 
\cite{weierstrass1870})
 entitled {\it \"{U}ber das sogennante Dirichlet'sche Princip}.%
\footnote{{\it On the so-called Dirichlet principle}.}
As Weierstrass notes in this paper, he had a handwritten copy of lecture notes from Dirichlet's lectures which he had received from Dedekind.  Weierstrass quotes several pages from these notes, and then points out that Dirichlet's proof of the Dirichlet principle (namely the proof of the existence of a minimizing function) was incomplete and not rigorous, and then produces an example of a similar type of calculus of variations problem in one dimension in which he showed that the minimum value of a specific energy integral was {\it not} assumed by any function in the class of functions being considered.

In the same year as Weierstrass's lecture, Hermann Amandus Schwarz (1843--1921), a student of Weierstrass, published a paper \cite{schwarz1870} which included a rigorous proof of the Riemann mapping theorem in the special case of simply-connected domains in the complex plane (as stated above).  This proof did not apply to the more general case of a simply-connected domain in a Riemann surface, which Riemann had formulated but also did not prove (this was proved later in the uniformization theorem of Paul Koebe (1882--1945), a student of Schwarz, and we will  discuss this in the Conclusion of this paper). Schwarz followed the outline of the proof of Riemann, replacing the Dirichlet principle argument with a convergent iterative argument involving a sequence of harmonic functions defined on specified open subsets of the given domain. Finally, in 1904 Hilbert \cite{hilbert1904} gave a proof of this disputed Dirichlet principle in the special context of one of Riemann's existence theorems in his paper on Abelian functions \cite{riemann1857}, thus justifying Riemann's original argument. 

Now we turn more specifically to the third of our major contributors to function theory, whom we have mentioned several  times already, namely, Karl Weierstrass.  His work stretched over a number of decades in the latter half of the 19th century and set standards of rigor and methodology that became a major force in how function theory (and more generally limiting processes and analysis in general) was perceived and used in the 20th century. The first papers of Weierstrass in the 1840s concerned themselves with specific problems in the theory of elliptic functions following up on the pioneering work of Abel and Jacobi.  During these early years Weierstrass wrote several fundamental papers which were only published later.

The main influence of Weierstrass in function theory came via his lecture courses in Berlin in the 1860s, which were published at the time, and which are all included in his collected works.  The first three volumes of his collected works \cite{weierstrass1894}, \cite{weierstrass1895}, and \cite{weierstrass1903} contain primarily his original papers over his professional lifetime (including those papers mentioned earlier that weren't published when they were written), and the following three volumes contain reprints of his lecture notes from his lectures on Abelian functions \cite{weierstrass1902}, his lectures on elliptic functions \cite{weierstrass1915}, and his lectures on applications of elliptic functions \cite{weierstrass1915a}.  

Let us mention here some of his principle results which have become part of the standard repertoire in function theory. Weierstrass defined a holomorphic function%
\footnote{Weierstrass used the term analytic function instead of holomorphic function, which was used regularly in the 20th century as well.  Later this became known as complex-analytic functions to contrast with the similarly defined real-analytic function defined as a locally convergent real power series of real variables.  Today holomorphic refers to the complex-analytic case, and one still uses the term real-analytic for the case of analytic functions of a real variable or variables.}
 to be a locally defined function of a complex variable defined near a point \(z_0\in \BC\) of the form
\[
f(z)= \sum_{n=0}^\infty a_n(z-z_0)^n,
\]
which converges in some disc of radius \(R\) centered at \(z_0\). If two such functions \(f_1\) and \(f_2\) are defined in discs \(\D_1\) and \(\D_2\) centered at two points \(z_1\) and \(z_2\), and if these two discs intersect, and if
\[
f_1|_{\D_1\cap \D_2} = f_2|_{\D_1\cap \D_2},\]
then \(f_2\) is said to be the {\em analytic continuation} of \(f_1\), and vice versa.  Moreover,
\[
f:=\left\{\ba{c}
f_1, z\in \D_1,\\
f_2, z\in \D_2,
\ea
\right.
\]
is a holomorphic function in \(\D_1 \cup \D_2\). More generally, a holomorphic function in a domain \(D\) is a function which admits such a power-series expansion near each point of \(D\).%
\footnote{This definition extends naturally to Riemann surfaces spread over domains in \(\BC\) as well.}
This definition of holomorphic as formulated by Weierstrass is equivalent to the definition used by Cauchy and Riemann as solutions of the Cauchy-Riemann equations.

Weierstrass brought much needed rigor to mathematical analysis, not only in function theory.  For instance he showed that a sequence of continuous functions on a domain \(D \subset \BR^n\)
\[
f_1(x), f_2(x), ..., f_k(x), ...
\]
which converges uniformly on compact subsets of the domain \(D\) has a limit that is continuous on \(D\). In the holomorphic setting he showed that if 
\[
f_1(z), f_2(z), ... , f_k(z),...
\]
is a sequence of holomorphic functions in a domain \(D \subset \C^n\) which converges uniformly on compact subsets of \(D\), then the limiting function is holomorphic in \(D\). 

Here we note that, although our emphasis in this section has been on holomorphic functions of one variable, Weierstrass (and others) considered holomorphic functions of several variables as well.  For instance, the definition of a holomorphic function \(f(z_1,...,z_n)\) of several complex variables \((z_1,...,z_n)\) can either be that the function has a convergent power series expansion of the form
\[
f(z_1,...,z_n) = \sum_{i_1,...,i_n}a_{i_1...i_n}z_1^{i_1}...z_n^{i_n},
\]
near each point of a domain, where \((z_1,...,z_n)\) are coordinates in \(\BC^n\), or, alternatively, one can require that the Cauchy-Riemann equations in \(\BC^n\) are satisfied, i.e.,
\[
\pd{f}{\overline z_j} (z_1,...,z_n) = 0,
\]
where
\[
\pd{}{\overline z_j} = \frac{1}{2}\left(\pd{}{x_j}+ i \pd{}{y_j}\right), j= 1,..., n.
\]
We note here that Riemann, Weierstrass and others were very interested in {\em Abelian functions}, which were functions of several complex variables in \(\BC^n\) which generalized of elliptic functions of one complex variable that we will encounter in the next section.

A given holomorphic function with a power series expansion at a given point has a {\it radius of convergence} for the series, and there are various criteria and descriptions of how to compute this due to Cauchy and others.  In particular, if the radius of convergence of a function \(f\) at a given point \(z_0\) is a finite number \(R<\infty\), then there is a least one boundary point \(z_1\) on the disc of radius \(R\) centered at \(z_0\) where \(f\) is singular and doesn't admit any analytic continuation to a larger open set containing that point.  For instance the function 
\[
f(z)=\frac{1}{z-1}
\]
has an expansion in the unit disc (the geometric series), and this does not converge at the point \(z=1\), and the function cannot be analytically continued beyond (or through) that point.  However, this function does have an analytic continuation through all other points on the boundary of the unit disc, as is very easy to see, since \(f(z)\) is holomorphic on \(\BC -\{1\}.\)  Weierstrass was the first to describe a function holomorphic on the unit disc which is singular at every boundary point of the unit disc.  

Here's a simple example of such a function given by a lacunary series,
\[
f(z)=\sum_{n=0}^\infty z^{2^n},
\]
and it easy to see that the series is divergent on the dense set of roots of unity of all orders on the boundary of the unit disc. Namely, for \(z=1\)
\[
f(1) = 1+1+1+\cdots +1\cdots,
\]
which diverges,
and for the two roots \(\e_1, \e_2\) of \(z^2=1\),
\bean
f(\e_j) &=& \e_j^1 + \e_j^2+\e_j^4 +\cdots +\e_j^{2^n}+\cdots,\\
&=& \e_j +1+1+\cdots +1+\cdots,
\eean
and for the four roots \(\e_1, \e_2,\e_3, \e_4\) of \(z^{2^2}=1\) we see that
\bean
f(\e_j) & =& \e_j^1+\e_j^2 +\e_j^4 + \e_j^8 +1+\dots +\e_j^{2^n}+\dots,\\
&=& \e_j +e_j^2 +1+1+\cdots +1+\cdots.
\eean
Using an induction argument, we see the series is divergent on this dense set, and hence at each of the boundary points (since no point where a function is holomorphic can be a limit point of singular points). 

By the Riemann mapping theorem, any simply connected domain has a function singular at every boundary point, and by a result of Weierstrass and G\"{o}sta Mittag-Leffler (1846-11927), which we mention below, this is true for all domains in \(\BC\).  However, a striking result for holomorphic functions of several complex variables due to Friedrich Moritz Hartogs (1874--1943) \cite{hartogs1906} at the beginning of the 20th century shows that this is {\em not} true for functions of two or more variables, and this led to a major new direction of research for complex analysis for functions of several complex variables in the 20th century.

As we mentioned earlier, Weierstrass started his mathematical career by studying elliptic functions as was formulated by Jacobi, and we met these functions earlier as \(\sn z, \cn z\), and \(\dn z\).  In his later work on elliptic functions Weierstrass (see his lectures on elliptic functions \cite{weierstrass1915}), introduced a new way to describe elliptic functions, which has now become one of the two standard ways of approaching elliptic functions (the other being that of Jacobi). We will describe this briefly here, and return to this when we discuss algebraic curves in the next section (see \cite{hurwitz-courant1964} or \cite{whittaker-watson1962} for a complete discussion of classical elliptic function theory). If \(\om_1,\om_2\) are two complex numbers with \(\Im \om_1/\om_2\ne 0\), then Weierstrass defines the Weierstrass \(\wp\)-function as
\be
\label{eqn:weierstrass-pfunction}
\wp(z) := \frac{1}{z^2} + \sum_{m^2+n^2 >0} \frac{1}{(z-m\om_1-n\om_2)^2} -\frac{1}{(m\om_1+n\om_2)^2}.
\ee
The series converges due to the extra term
\[
\frac{1}{(m\om_1+n\om_2)^2},
\]
which is added as a Mittag-Leffler correction term (we will discuss this somewhat later in this section).  It is easy to see that this function has the following properties (assuming the convergence, which requires some work):
\bi
\item \(\wp (z+m\om_1+n\om_2) = \wp(z)\), i.e., \(\wp(z)\) is doubly-periodic in \(\BC\) with periods \(\om_1\) and \(\om_2\).
\item \(\wp(z)\)  is a meromorphic function on \(\BC\) with a single double pole in each period parallelogram.%
\footnote{The period parallelograms defined by the periods \(\om_1, \om_2\) are the translates in the complex plane by integers of the form \(m+in\) of the fundamental period parallelogram with the four vertices \(0, \om_1, \om_2,\) and \(\om_1+\om_2\) (see, e.g., \cite{hurwitz-courant1964}).}
\item The derivative \(\wp'(z)\) of the Weierstrass \(\wp\)-function is also a doubly periodic function with periods \(\om_1\) and \(\om_2\).
\ei
The two functions \(\wp(z)\) and \(\wp'(z)\) play an analogous role in elliptic function theory to the original doubly-periodic functions \(\sn(z)\) and \(\cn(z)\) of Jacobi, and we will return to them later.  A classical reference for this infinite series approach to elliptic functions is the text by Karl Boehm \cite{boehm1908}.

The final result of Weierstrass we want to mention in this section is often referred to as the {\it Weierstrass factorization theorem} or the {\em Weierstrass product theorem} which was published in 1876 \cite{weierstrass1876}.  This result is an important  generalization of the fundamental theorem of algebra, which asserts that any polynomial can be expressed in terms of factors
\be
\label{eqn:fundamental-algebra}
p(z) = c(z-a_1)^{m_1}\cdots(z-a_k)^{m_l},
\ee
where \(a_1,\dots,a_k\) are the roots of the polynomial with multiplicities \(m_1,\dots,m_k\) and \(c\) is a constant. Let \(f(z)\) be a holomorphic function in the complex plane \(\BC\) (such an \(f\) is called an {\em entire function}) with zeros at a possibly infinite set of points \(a_1,\dots,a_k,\dots,\), then the Weierstrass factorization theorem asserts that \(f\) can be represented as an infinite product similar to the finite product in (\ref{eqn:fundamental-algebra}),
\be
\label{eqn:weierstrass-product}
f(z) = z^m e^{g(z)}\prod_{n=1}^\infty \left(1-\frac{z}{a_n}\right)e^{\left(\frac{z}{a_n}\right)
+\frac{1}{2}\left(\frac{z}{a_n}\right) + \cdots+\frac{1}{m_n}\left(\frac{z}{a_n}\right)^{m_n}},
\ee
where \(m, m_n\) are integers, and \(g\) is an entire function (see Ahlfors \cite{ahlfors1953} or any standard complex analysis text for a discussion and proof of this theorem).  The result is often formulated in the following manner.  Let \(a_1,\dots,a_k,\dots\) be any sequence of points in the plane such that \(\lim a_k = \infty\), then there exists an entire function with zeros at precisely these points (namely use the formula (\ref{eqn:weierstrass-product})). Weierstrass introduced the exponential factors in the infinite product to insure its convergence.

A consequence of this theorem, as observed by Weierstrass, is that any meromorphic function in \(\BC\) can be expressed as the quotient of two entire functions.  These results of Weierstrass were generalized by Mittag-Leffler in 1884 \cite{mittag-leffler1884} from the case of the complex plane to an arbitrary domain in the following sense.  Let \(D\) be a domain in the complex plane and let \(a_1,\dots,a_k,\dots\) be an infinite sequence of points in \(D\) with no accumulation points in \(D\), i.e., each point \(a_k\) is isolated in \(D\), then there exists a function \(f(z)\) holomorphic in \(D\) which has zeros precisely at the points \(a_k\).  If we consider a domain \(D\) with a boundary \(\bd D\), and let \(a_k\) be a set of points in \(D\) again with no accumulation points in \(D\) and which is dense on the boundary of \(D\), then the Weierstrass function with these zeros has no analytic continuation beyond any boundary point (a result we alluded to above in the context of the Riemann mapping theorem and lacunary series for simply connected domains).  

We now turn to the final topic of this section on holomorphic functions, the Mittag-Leffler Theorem, mentioned briefly above.  As we just saw, the Weierstrass factorization theorem showed that for given prescribed zeros, one can find a holomorphic function with those zeros.  By taking reciprocals, one could find a meromorphic function with poles of a certain order at those same points by prescribing the multiplicity of the zeros or the order of the poles.  A variation on this question was raised and solved in two earlier papers of Mittag-Leffler's in 1877 \cite{mittag-leffler1877a,mittag-leffler1877b} which we formulate here.  A meromorphic function \(f(z)\) near a pole at a point \(z_0\) of order \(m\) has a Laurent expansion at \(z_0\) of the form
\[
\frac{a_{-m}}{(z-z_0)^m} +\cdots+\frac{a_{-1}}{(z-z_0)} + \sum_{n=0}^\infty a_n(z-z_0)^n,
\]
where the infinite series converges in the neighborhood of \(z_0\), and the finite number of terms of powers of \(\frac{1}{(z-z_0)}\) all converge to \(\infty\) at \(z_0\) and represent what is called the {\em principal part} of the meromorphic function at \(z_0\).  We saw this earlier, when we identified \(a_{-1}\) as the residue of \(f(z)\) at \(z_0\).  More generally, if the meromorphic function \(f(z)\) has poles at \(z_k\) in a domain \(D\subset\BD\), then there are functions \(p_k(z)\) and  \(g_k(z)\) defined near each point \(z_k\) such where \(p_k(z)\) has the form
\be
\label{eqn:principal-part}
p_k(z) = \frac{a^k_{-m_k}}{(z-z_k)^{m_k}} +\cdots+\frac{a^k_{-1}}{(z-z_k)},
\ee
where  \(p_k(z)\) is the principal part of the meromorphic function \(f(z)\) at \(z_k\) and \(f(z)-p_k(z) := g_k(z)\) is holomorphic near \(z_k\).  

The question Mittag-Leffler raised  was the following.  Given a discrete set of points \(z_k\) in a domain \(D\) and for each point \(z_k\) given a polynomial of the form given in (\ref{eqn:principal-part}), does there exist a meromorphic function \(f(z)\) in \(D\) and locally defined holomorphic functions \(g_k(z)\) defined near each point \(z_k\) such that near \(z_k\)
one has
\[
f(z)-p_k(z) = g_k(z).
\]
Mittag-Leffler showed that this was true, and the result is known as the {\it Mittag-Leffler theorem}.

Let us sketch Mittag-Leffler's proof of this in the simple case where \(z_k\) is a discrete sequence in the complex plane and 
the principal parts which are given are simple poles of the form
\[
p_k(z) = \frac{a_k}{(z-z_k)}.
\]
Assume the points \(z_k\) are ordered such that
\[
|z_1| \le |z_2| \le \cdots \le |z_k|\cdots,
\]
and let \(\D_k\) be concentric circles centered at the origin whose radii increase in such a fashion that
\[
z_j \notin \overline{\D_k} ,\;\textrm{for} j\ge k,
\]
where we use the standard notation \(\overline{K}\) to mean the closure of a set \(K\) in \(\BR^n\) (the set of all accumulation points of \(K\)). Thus each \(p_k(z)\) is holomorphic on a neighborhood of \(\D_k\).  This implies that \(p_k(z)\) can be expanded in a power series centered at the origin which converges on a neighborhood of \(\overline \D_k\), and hence \(p_k(z)\) can be approximated on \(\overline{\D_k}\) by a polynomial \(h_k(z)\) such that
\[
|p_k(z)-h_k(z)|<\frac{1}{2^k},\;\textrm{for}\; z\in\overline{\D_k}.\]
It follows that the series
\bean
f(z) &=& \sum_{k=1}^\infty (p_k(z)-h_k(z)),\\
&=& \sum_{k=1}^\N (p_k(z)-h_k(z))+\sum_{k=N+1}^\infty (p_k(z)-h_k(z)).
\eean
converges uniformly on \(\D_n, N=1,2,\dots,\)and \(f(z)\) is a well-defined meromorphic function on \(\BC\) with principal parts \(p_k(z)\) near each pole \(z_k\).

Mittag-Leffler first proved this result in his two papers (written in Swedish) \cite{mittag-leffler1877a,mittag-leffler1877b} for the case of \(D=\BC\), and in his longer {\em Acta Mathematica} paper%
\footnote{{\em Acta Mathematica} was founded by Mittag-Leffler in 1882, and initially, for a number of years, all the papers were in French, the international language of its time.}
\cite{mittag-leffler1884} for arbitrary domains. Carl Runge (1856-1927) gave a new and simpler proof of Mittag-Leffler's theorem in 1885 \cite{runge1885} which involved a new approximation theorem, now called Runge's theorem or Runge's approximation theorem, which showed how one can approximate holomorphic functions on a multiply-connected domain \(D\) uniformly on compact subsets of the domain by rational functions with poles in the bounded components of the complement of \(D\).  If \(D\) is simply connected, then the function can be approximated by a polynomial in the same manner. These results of Weierstrass, Mittag-Leffler and Runge all have generalizations to holomorphic functions of more than one complex variable, and they play an important role in the further development of this field of research as it developed in the 20th century.

The work of Cauchy, Riemann, and Weierstrass  generated immense interest in the mathematicians of the second half of the 19th century.  Initially there were more or less three schools of thought following these three innovators (the integral theorems of Cauchy, the differential equations of Riemann and the power series of Weierstrass), but at the end of the 19th century the concept of {\em function theory}, the theory of holomorphic functions of a complex variable, reached a significant stage of maturity and used all of the tools available to study new levels of problems which arose. At the turn of the 19th to the 20th century various mathematicians began the study of holomorphic functions of several complex variables, and new phenomena (Hartogs' theorem) changed the direction of research in this higher dimensional setting.

\section{Riemann surfaces}
\label{sec:riemann-surfaces}
\label{sec:riemann-surfaces}
We now turn to our final topic of this paper, the creation of the theory of Riemann surfaces. This singular creation by Riemann in his dissertation of 1854 \cite{riemann1854} and his papers on Abelian functions in 1857 \cite{riemann1857a,riemann1857b,riemann1857c,riemann1857d}%
\footnote{In Riemann's collected works \cite{riemann1876} these four papers are published together under the heading of a single paper entitled {\em Theorie der Abel'sche Functionen}.  The first three papers summarize and clarify concepts developed in his dissertation from 1854 as tools for his detailed study of Abelian integrals and Abelian functions in the fourth paper. We will often refer simply to his Abelian functions paper of 1857.} 
developed over the next century into the very rich subject of {\em complex manifolds} of arbitrary dimension (Riemann surfaces being the case of a one-dimensional complex manifold), with strong overlaps with algebraic geometry, as we will indicate later.

Riemann's motivation for his creation of Riemann surfaces arose from the study of multivalued functions, and in particular in the multivalued functions that arose in Abel's work on generalizations of elliptic integrals that we discussed in Section \ref{sec:abel-theorem}.  Multivalued functions had been a topic that had occupied mathematicians a great deal during the several centuries preceding Riemann's work, and the ambiguities that arose was a major concern.  A fundamental example that arose in Abel's work was the study of a  \(y(x)\) which arose as a solution of the algebraic equation
\[
y^n+a_{n-1}y^{n-1} + \cdots+a_0(x)=0,
\]
where \(a_k(x)\) are polynomials in \(x\). A different and familiar set of examples are given by the inverses of elementary transcendental functions such as \(e^z\) and \(\sin z\), which we denote by \(\log z\) and \(\textrm{arcsin} z\), and which were intensely studied in the 18th century. These particular multivalued functions have an infinite number of different values at a given point.

A key ingredient in Riemann's creation of Riemann surfaces  was the notion of analytic continuation of a holomorphic function, which we discussed briefly earlier.  As Riemann observes on p. 102 of his Abelian function paper \cite{riemann1857a}, the function \(\log (z-a)\), a well-known multivalued function, when continued analytically on a simple closed path around the point \(a\), increases or decreases its value by \(2\pi i\), depending on the direction of the path. If we let \(z-a=re^{i\th}\) be polar coordinates at the point \(a\), then
\[
\log(z-a) = \log r e^{i\th}= \log r + i\th,
\]
and as \(\th\) varies from \(0\) to \(2\pi\), \(\log (z-a) \) varies from \(\log r\) to \(\log r +2\pi i\). This well known phenomenon played a major role role in Riemann's work. 

Riemann considered the possible different analytic continuations of a given holomorphic function to be {\em branches} (Zweige) of the function, and he defines a {\em branch point} (Verzweigungsstelle) as a point around which one branch moves into another (in our example \(a\) is a branch point for the multivalued function \(\log(z-a)\)). He then describes (on pages 103-4) of \cite{riemann1857a} the surfaces he wants to consider:
\begin{quote}
F\"{u}r manche Untersuchungen, namentlich f\"{u}r die Untersuchung algebraischer
und Abel'scher Functionen ist es vortheilhaft, die Verzweigungsart
einer mehrwerthigen Function in folgender Weise geometrisch darzustellen.
Man denke sich in der \((x,y)\)-Ebene eine andere mit ihr zusammenfallende
Fl\"{a}che (oder auf der Ebene einen unendlich d\"{u}nnen K\"{o}rper) ausgebreitet,
welche sich so weit und nur so weit erstreckt, als die Function gegeben ist.
Bei Fortsetzung dieser Function wird also diese Fl\"{a}che ebenfalls weiter ausgedehnt
werden. In einem Theile der Ebene, f\"{u}r welchen zwei oder mehrere
Fortsetzungen der Function vorhanden sind, wird die Fl\"{a}che doppelt oder
mehrfach sein; sie wird dort aus zwei oder mehreren Bl\"{a}ttern bestehen, deren
jedes einen Zweig der Function vertritt. Um einen Verzweigungspunkt der
Function herum wird sich ein Blatt der Fl\"{a}che in ein anderes fortsetzen, so
dass in der Umgebung eines solchen Punktes die Fl\"{a}che als eine Schraubenfl\"{a}che
mit einer in diesem Punkte auf der \((x,y)\)-Ebene senkrechten Axe
und unendlich kleiner H\"{o}he des Schraubenganges betrachtet werden kann.
Wenn die Function nach mehren Uml\"{a}ufen des um den Verzweigungswerth
ihren vorigen Werth wieder erh\"{a}lt (wie z.B. \((z-a)^{\frac{m}{n}}\), wenn \(m,n\) relative
Primzahlen sind, nach \(n\) Uml\"{a}ufen von \(z\) um \(a\)), muss man dann freilich annehmen,
dass sich das oberste Blatt der Fl\"{a}che durch die \"{u}brigen hindurch
in das unterste fortsetzt.

Die mehrwerthige Function hat f\"{u}r jeden Punkt einer solchen ihre
Verzweigungsart darstellenden Fl\"{a}che nur einen bestimmten Werth und kann
daher als eine v\"{o}llig bestimmte Function des Orts in dieser Fl\"{a}che angesehen
werden.%
\footnote{
For many investigations, namely for the investigation of algebraic and Abelian functions it is advantageous to geometrically represent the branching nature of a multivalued function in the following manner.  One imagines a surface (or an infinitesimally thin body) coinciding with and spread over the \((x,y)\)-plane, which is extended as far as, and only as far as, the function is given. As the function is analytically continued, the surface will be further extended. In a region of the plane, for which two or more continuations are present, the surface will be covering twice or more times the region; it will consist of two or more sheets, each of which will represent a branch of the function.  Around a branch point the function will continue from one sheet of the surface so that in the neighborhood of such a point the surface can be considered as a helicoid with a vertical axis through this point, and infinitesimally small heights of the screw thread from one revolution to another.  If the function comes back to its same value after several such revolutions (as happens, for instance with \((z-a)^{\frac{m}{n}}\), if  \(m,n\) are relatively prime numbers, after \(n\) cycles of \(z\) around \(a\)), then one has to assume that the upper sheet moves through the other sheets to the bottom sheet.

The multivalued function has, for every point of such a surface representing its branching, only one definite value, and can thereby be a completely well-determined function of position on this surface.}
\end{quote}

This description of a Riemann surface and Riemann's further use of it in these four papers became the standard way to describe Riemann surfaces for the next half century until Hermann Weyl (1885--1955) introduced the first abstract version of a Riemann surface as a topological manifold with a complex structure in 1913 \cite{weyl1913}.  Fundamentally, the Riemann surface as a covering of the extended complex plane gave local coordinates at each point of the surface except at the branch points, and at a branch point \(a\) of the type described in the quote above from Riemann one can use \(\z=(z-a)^\frac{1}{n}\) as a local coordinate chart at this point. Riemann also added the point at infinity, \(\infty\), to each sheet of the Riemann surface, thus giving rise to a closed or compact Riemann surface, with the local coordinate system \(\z=1/z\) at the point at infinity. The system of local coordinate charts for points of a Riemann surface was formalized by Hermann Weyl in his book mentioned above, but the 19th century mathematicians worked quite well with the structure Riemann set up and which we have summarized here.

Riemann's second paper in this series of four \cite{riemann1857b} has the title {\em Lehrs\"{a}tze aus der analysis situs f\"{u}r die Theorie der Integrale von zweigliedrigen vollst\"{a}ndigen Differentialen}%
\footnote{{\em Theorems from analysis situs for the theory of integrals of two-fold complete differentials}}.
{\em Analysis situs} was a somewhat common name in the 19th century for what became {\em topology} (or what became {\em algebraic topology}, more specifically) in the 20th century (we will discuss the origin of the word "topology" shortly).  The term {\em analysis situs} originated in work of Gottfried Wilhelm Leibniz (1646-1716) which was contained in correspondence between Leibniz and Christiaan Huygens (1629 – 1695) with the first and most fundamental letter being from Leibniz to Huygens on 8 September 1679 \cite{leibniz1679} in which he compared geometry of magnitude with geometry of position ({\em situm}), and felt that he could contribute to this new way of thinking by expressing positions of geometric objects and their relationships with symbols, just as algebra used symbols to represent relationships between numbers.  Leibniz felt this was very important, but Huygens remained skeptical of this optimistic young mathematician's ideas in this direction, while recognizing the significance of his work on infinitesimal analysis. The recent book by Vincenzo Risi {\em Geometry and Monadology} \cite{risi2007} has a very interesting analysis of Leibniz's work on analysis situs.

\begin{figure}
\vspace{6pt}
\centerline{
	\includegraphics[width=10cm]{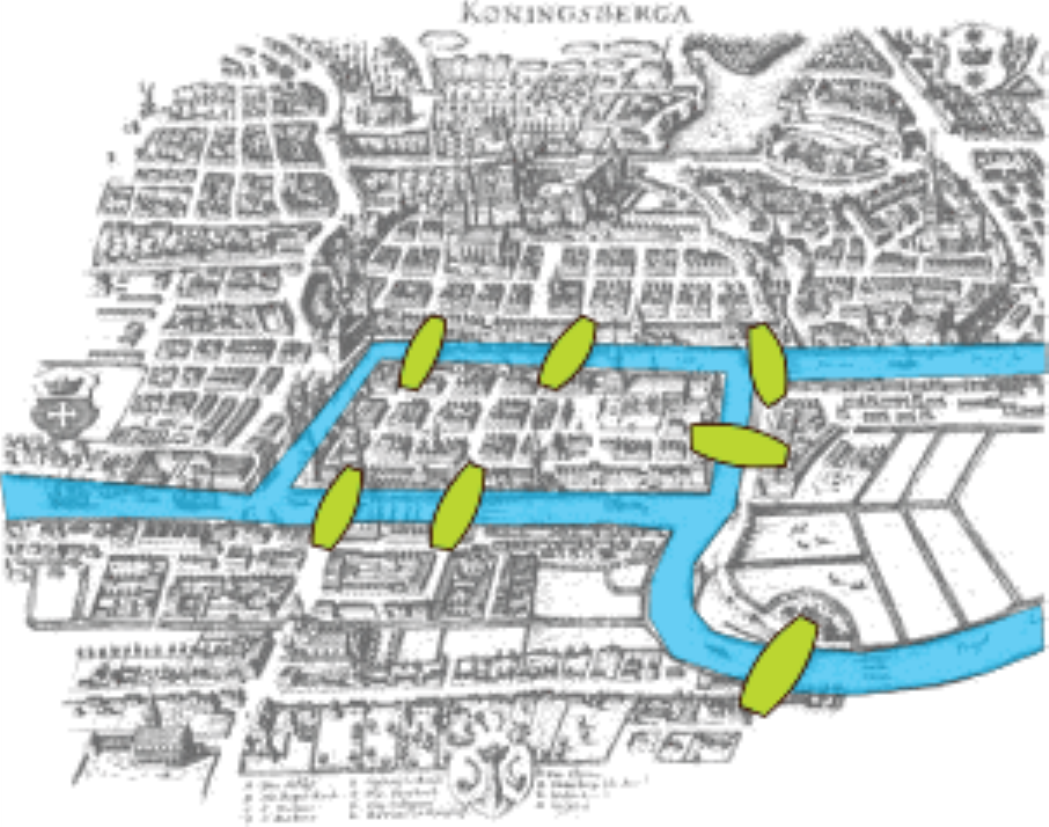}}
 \caption{The K\"{o}nigsberg Bridges Problem} 
 \label{fig:koenigsberg-bridges}
\end{figure}

The work of Leibniz forms part part of the inspiration for the book {\em G\'{e}om\'{e}rie de Position} by Lazare Nicolas Marguerite Carnot (1753--1823) \cite{carnot1803} which had a major influence on projective geometry and where a major impulse was to investigate geometric phenomena that were not dependent on measurement of distances. The first definitive work on topology after the initial impetus by Leibniz came from Euler in 1735 in his famous solution of the {\em Seven Bridges of K\"{o}nigsberg} problem \cite{euler1735}.  Note the title of this paper, {\em Solutio problematis ad geometriam situs} contains the phrase {\em geometriam situs}, exactly the term used by Leibniz, which Euler cited and which Carnot used in the title of his book. Euler was able to abstract the problem (Can one find a path crossing all seven bridges precisely once? See Figure \ref{fig:koenigsberg-bridges}) 
to be a problem in what became graph theory and showed there was no such path.  Somewhat later Euler gave the first example of what became the {\em Euler characteristic} of a polyhedral surface \cite{euler1752}, which, in all of its generalizations in algebraic topology and other forms of geometry, has played  an important role in mathematics. In the 19th century Johann Benedict Listing (1808-1992) published a short book entitled {\em Vorstudien zur Topologie} \cite{listing1848}, which followed up on the work of Euler and developed a theory of knots,  a special topic in algebraic topology today concerning curves embedded in (usually) some Euclidean space, and in Listing's case, in \(\BR^3\).  Then a few years later came Riemann's work that we are discussing here on the algebraic topology of surfaces. 

\begin{figure}
\vspace{6pt}
\centerline{
	\includegraphics[width=12cm]{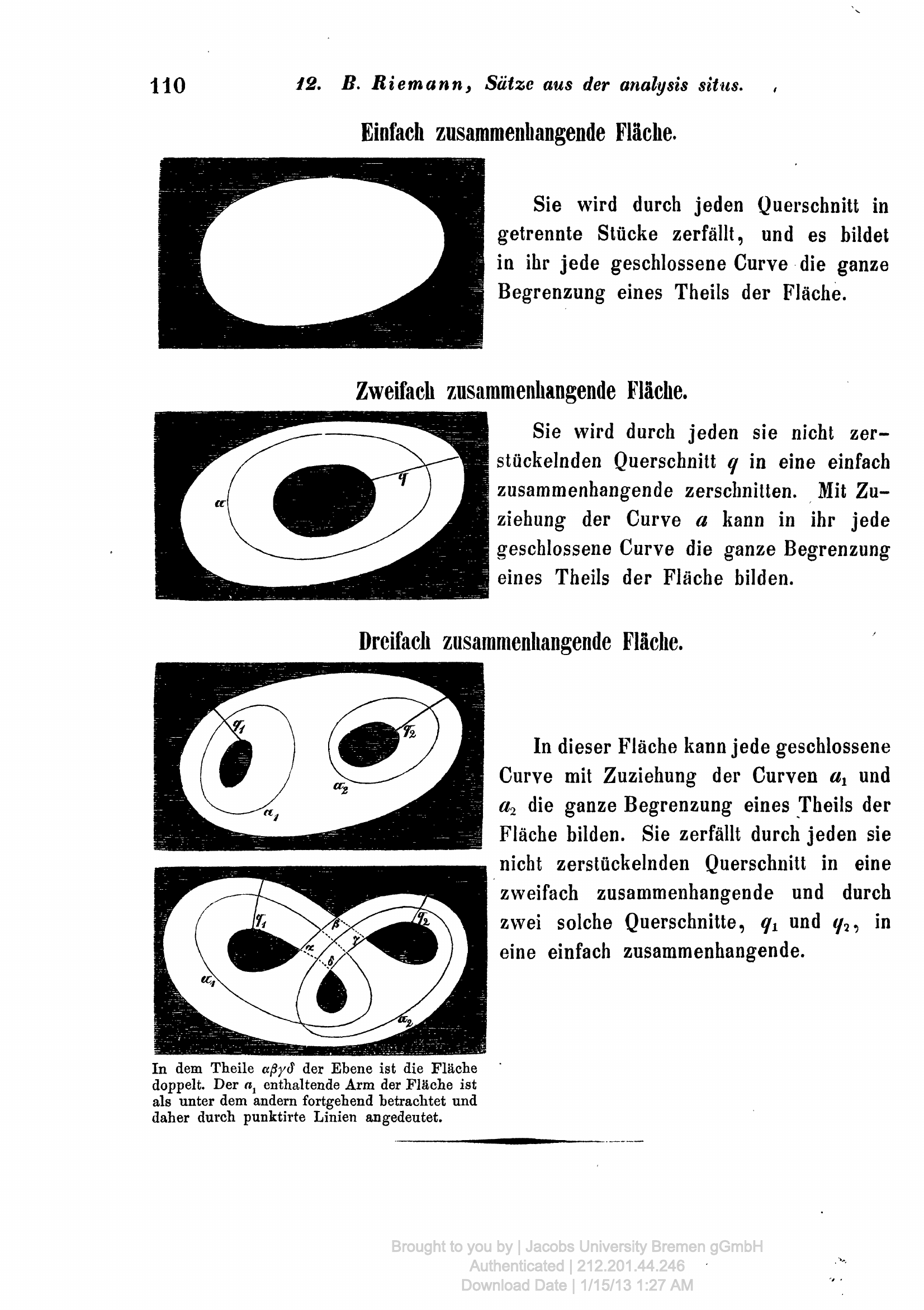}}
 \caption{Riemann's connectivity in his Analysis situs paper \cite{riemann1857b}} 
 \label{fig:riemann-connectivity}
\end{figure}

Riemann defined both in his dissertation \cite{riemann1851} and in the second paper in the Abelian functions series \cite{riemann1857b} the notion of {\em connectivity} of a surface.  First he defined a {\em simply-connected} surface \(S\) to be any surface such that the complement of any closed curve (or a curve from one boundary point to another for a surface with a boundary) in the surface consisted of two components.  He then defined an \(n+1\)-connected surface to be a surface \(S\) whereby \(n\) suitably chosen curves deleted from the surface would give a simply-connected subdomain \(S'\subset S\).  He showed that this concept was well-defined and used it extensively in the remainder of his Abelian functions papers.  In his dissertation he had used the notion of a simply-connected domain to formulate his Riemann mapping theorem, which we discussed earlier.  In his Abelian functions paper, the non-simply-connected compact surfaces play the most important role, as that is where the theory of Abelian functions can be developed.  In Figure \ref{fig:riemann-connectivity} we see diagrams from Riemann's paper \cite{riemann1857b} which illustrates the notion of connectivity.
Note that the illustration at the bottom of Figure \ref{fig:riemann-connectivity} shows the possibility of two sheets of a surface overlapping as described in his definition of a Riemann surface as a branched covering.

\begin{figure}
\vspace{6pt}
\centerline{
	\includegraphics[width=12cm]{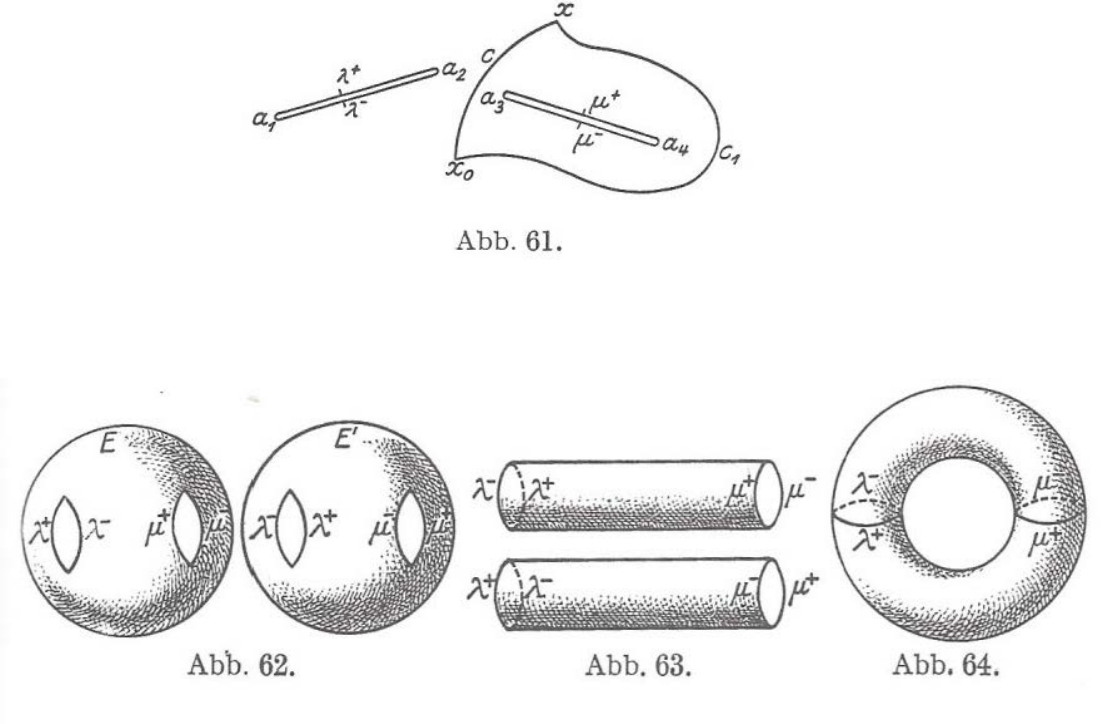}}
 \caption{Example of the Riemann surface of \(w^2=(z-a_1)(z-a_2)(z-a_3)(z-a_4)\) from pp. 236-237 of Hurwitz and Courant \cite{hurwitz-courant1964}} 
 \label{fig:hurwitz-courant}
\end{figure}

Let's consider a simple example illustrating the topology of a compact Riemann surface arising from an algebraic function.  Let the Riemann surface \(S\) be defined by the polynomial
\be
\label{eqn:elliptic-curve-example}
w^2=(z-a_1)(z-a_2)(z-a_3)(z-a_4),
\ee
where the points \(a_1, a_2, a_3, a_4\) are distinct, and let \(\l\) and \(\m\) be two nonoverlapping cuts joining \(a_1\) to \(a_2\) and \(a_3\) to \(a_4\) as illustrated in the top part of Figure \ref{fig:hurwitz-courant}.
Here \(\l^\pm, \m^\pm\) indicate the two sides of the cuts \(\l\) and \(m\).  By opening the cuts as indicated in the bottom portion of Figure \ref{fig:hurwitz-courant}, we see that we can glue the two copies of the Riemann sphere \(E\) and \(E'\) along cylinders joining them to create the torus on the lower right hand sided of the figure. Thus we conclude that the topology of the Riemann surface defined by (\ref{eqn:elliptic-curve-example}) is that of a torus.

In the third of these preliminary papers \cite{riemann1857c} Riemann revisits the Dirichlet principle and shows that there exists meromorphic functions with prescribed poles or logarithmic singularities.  Strictly speaking, a meromorphic function, by definition, does not have logarithmic singularities, but Riemann is interested also in the {\em integrals} of meromorphic functions from a specified point to an indefinite point, that is he is interested in Abelian integrals on a Riemann surface, and such integrals can have a logarithmic singularity, e.g., 
\[
\int_{z_0}^z \frac{dz}{z} =\log (z-z_0).
\]

In the fourth paper in this series {\em Theorie der Abel'schen Functionen} \cite{riemann1857d} Riemann develops his theory of {\em Abelian functions}, a vast generalization of elliptic functions which are defined as several variable inverses of Abelian integrals on Riemann surfaces of genus \(> 1\). Riemann's work followed up on some announced results of Weierstrass whose proofs hadn't yet appeared at the time  Riemann wrote his paper.  The basic theory which evolved became known as the theory of Abelian functions (this name was adopted by Riemann in the paper we are discussing), and was further developed by numerous mathematicians in the following century.  We will discuss this briefly later in this section.

\begin{figure}
\vspace{6pt}
\centerline{
	\includegraphics[width=12cm]{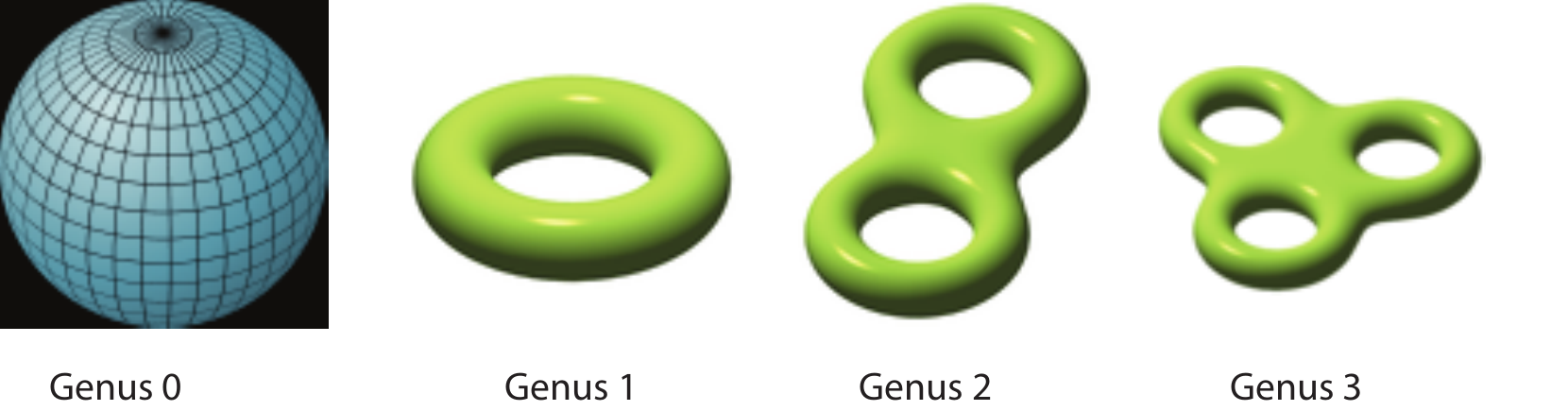}}
 \caption{Surfaces of genus 0, 1, 2, and 3}
 \label{fig:surfaces-genus}
\end{figure}

We want to discuss some aspects of this paper which directly relate to the theory of complex manifolds in the 20th century.  Let now \(S\) be a Riemann surface defined as a branched covering defined by the algebraic function
\be
\label{eqn:polynomial-equation}
F(z,w)= w^n + a_{n-1}w^{n-1} + \cdots+a_0=0,
\ee
where \(F(z,w)\) is an irreducible polynomial of degree \(n\) in \(w\) and of degree \(m\) in \(z\). Then the function \(w(z)\) defined by (\ref{eqn:polynomial-equation}) is a single-valued function on \(S\), and \(S\) is a compact Riemann surface, where we have added a point at infinity to each sheet of the Riemann surface, as did Riemann.  Let us suppose that \(S\) is \((2p+1)\)-connected, for \(p\ge 0\) (that the connectivity is odd for compact Riemann surfaces was shown by Riemann in his paper). In Figure \ref{fig:surfaces-genus}
we see an example of illustrations of surfaces of genus 0, 1, 2, and 3.  

Now consider the Abelian integral
\[
A(z)= \int_{z_0}^z R(z,w(z))dz,
\]
along any path in \(S\) joining \(z_0\) to \(z\), where \(R(z,w)\) is a rational function of \(z\) and \(w\) as we discussed in Section \ref{sec:abel-theorem}. Note that in Section \ref{sec:abel-theorem} we dealt with real integrals (with singularities in the integrand), and here we are dealing with path integrals in the complex plane or on a Riemann surface spread over the complex plane (also with possible singularities in the integrand, depending on the path). We define as before that  \(A(z)\) is an {\em Abelian integral of the first kind} if \(A(z)\) is finite for all points \(z\in S\).  We will see examples of this momentarily. We note also that \(A(z)\) is in principle a {\em multivalued} function if the genus \(p> 0\), since different paths along portions of a cycle which is not homologous to zero can have different values.  

If \(p=0\), then there are no Abelian integrals of the first kind.  To see this, suppose that \(A(z)\) is an Abelian integral on \(S\), where \(p=0\).  Then \(S\) is simply-connected and hence \(A(z)\) is single-valued and holomorphic at each point of \(S\).  But \(S\) is compact and therefore at some point of \(S\), by continuity, \(|f(z)|\) assumes a maximum value, which is necessarily an interior point, and hence \(A(z)\) must be constant.  However,  a constant cannot be an Abelian integral, since it's derivative would be zero, which cannot be the integrand of an Abelian integral as defined above.  

Abelian integrals of the second kind are characterized by  \(A(z)\) having  poles at some points of \(S\) and Abelian integrals of the third kind are characterized by  \(A(z)\) having logarithmic singularities at some points of \(S\).  We will not be discussing these further here, but we will see that Abelian integrals of the first kind are intimately connected to the topology of the surface \(S\).

Riemann shows in his paper \cite{riemann1857d} that on the Riemann surface \(S\) of genus \(p\) as defined above, there are precisely \(p\) linearly independent Abelian functions of the first kind \(A_1(z), A_2(z),\dots,A_p(z)\), and that if \(A(z)\) is any Abelian function of the first kind on \(S\), then there are constants \(\a_1,\dots,\a_p\) such that
\[
A(z) = \a_1A_1(z) +\cdots \a_pA_p(z) + \textrm{const}.
\]
Thus the number of linearly independent Abelian functions of the first kind is the genus \(p\)%
\footnote{In a later text on Abelian functions by Clebsch and Gordan written some 10 years after Riemann's work \cite{clebsch-gordan1866}, they compute the number \(p\) for a Riemann surface defined by a polynomial of the type \(F(z,w)\) in terms of degrees of the polynomial and numbers of double points and cusps and refer to this simply as the number of Abelian integrals of the first kind.}.  Let's give a couple of examples to illustrate this.

We consider  \(F(z,w)\) to be of the form, for some \(k>0, |k|\ne 1\), 
\be
\label{eqn:polynomial-elliptic}
F(z,w)= w^2-(z^2-1)(z^2-k^2),
\ee
 letting \(R(z,w)= 1/w\), and  setting
\[
A_1(z) = \int^z_0\frac{dz}{w(z)},
\]
(this is the elliptic integral considered by Abel (\ref{eqn:legendre-integral}) in our discussion of Abel's work on elliptic functions).  If \(z\ne \pm 1, \pm k \), then \(A_1(z)\) is finite.  Now if \(\g\) is a path from 0 to \(\infty\) not passing through these same  points, then 
\[
\lim_{z\rightarrow\infty}A(z)=v,\; \textrm{and}\;|v|< \infty,
\]
which is easy to verify.  Namely any path \(\g\) going from 0 to \(\infty\) and missing the points \(\pm 1,\pm (1/k)\) can be modified to be a path from 0 to a point \(R_0> \max \{1, 1/k\}\) and from \(R_0\) to \(\infty\) via a path on the positive real axis.  Then the resulting path integral part on the real axis would be
\[
\int_{R_0}^\infty \frac{dx}{\sqrt{(x^2-1)(x^2-k^2)}}.
\]
But, for \(x\in [R_0,\infty)\), there is a constant \(K\) such that
\[
\frac{1}{\sqrt{(x^2-1)(x^2-k^2)}}\le K\frac{1}{x^2},
\]
and thus \(\lim_{z\rightarrow \infty} A_1(z)\) exists and is finite.

Now we have to examine the behavior of \(A_1(z)\) at the singular points of the integrand.  Consider the point \(z=1\), and we want to show that \(\lim_{z\rightarrow 1} A_1(z)\) exists and is finite.  To do this we make a change of variable at this point of the form
\[
\z=(z-1)^{\frac{1}{2}},
\]
then
\[
d\z=\frac{1}{2}(z-1)^{-\frac{1}{2}}dz,
\]
which gives
\bean
dz&=&2\z d\z,\\
z&=&1+\z^2,
\eean
and hence
\bean
A_1(z) =A_1(1+\z^2)&=&\int_i^{1+\z^2} \frac{2\z d\z}{\sqrt{\z^2(\z^2+2)(\z^2+1+k)(\z^2+1-k)}},\\
	&=& \int_i^{1+\z^2} \frac{2 d\z}{\sqrt{(\z^2+2)(\z^2+1+k)(\z^2+1-k)}},
\eean
which has a nonsingular integrand near \(\z=0\), and hence \(A_1(z)\) is holomorphic in a neighborhood of \(z=1\).  A similar argument holds for the other singular points \(\{-1,k,-k\}\), and hence \(A_1(z)\) is finite at all points of \(S\), as Legendre and others knew, long before Riemann, but in the real variable context.

Now the Riemann surface for the polynomial (\ref{eqn:polynomial-elliptic}) has genus \(p=1\), as we discussed earlier, and, thus, according to Riemann, the only Abelian integrals of the first kind in this case are constant multiples of \(A_1(z)\) plus a possible constant.

A second example is the case of a hyperelliptic Riemann surface defined by an equation of the form
\[
F(z,w)=w^2-(z-a_1)(z-a_2)\cdots(z-a_{2p}),
\]
where \(a_k \ne 0\) are distinct complex numbers, then 
\be
\label{eqn:abelian-integrals}
A_k(z) := \int_0^z \frac{z^k dz}{w(z)}, k = 0,\dots,p-1,
\ee
are \(p\) distinct Abelian integrals of the first kind, as we mentioned in our discussion of Abel's work in Section \ref{sec:abel-theorem}.  The proof that these are Abelian integrals of the first kind is similar to that given above in the elliptic case, and also it is easy to verify that the genus for this Riemann surface is also \(p\). Namely, there are \(2p\) branch points for this two-sheeted Riemann surface. Hence the Abelian integrals in (\ref{eqn:abelian-integrals}) form a basis for the vector space of Abelian integrals of the first kind.

Let's look at the integrands of these Abelian integrals of the first kind. Notice that originally in the work of Euler, Legendre, Abel and others the integrals were integrals on the real axis with singularities (with the inherent multivaluedness).  For instance consider the integral 
\[
A_1(x)= \int_0^x \frac{dx}{\sqrt{(x^2-1)(x^2-k^2)}},
\]
in the elliptic case. The integrand is
\[
\frac{dx}{\sqrt{(x^2-1)(x^2-k^2)}},
\]
and this becomes, going to the complex plane
\[
\frac{dz}{\sqrt{(z^2-a)(z^2-k^2)}},
\]
which has the form
\be
\label{eqn:one-form-example}
f(z)dz,
\ee
where\(f(z)\) is a multivalued meromorphic function on \(\BC\).  We shall see below how (\ref{eqn:one-form-example}) can be interpreted as a single-valued and, in fact, holomorphic one form on the Riemann surface defined by \(w^2= (z^2-1)(z^2-k^2)\).

Let \(S\) be a Riemann surface spread over the extended complex plane of the type described by Riemann with local coordinates \(\z=z-a\) at nonbranching points, \(\z=(z-a)^\frac{1}{k},\) at a branching point of order \(k\), and \(\z = 1/z\) at infinity, assuming, for simplicity that the point at infinity is not a branch point.  A {\em meromorphic one-form on \(S\)} is defined to be a one form defined with respect to any local coordinate chart \(\z\) as above to be of the form
\[
\om=f(\z)d\z,
\]
where \(f(\z)\) is meromorphic in \(\z\).  And, if we transform from one coordinate system \(z\) to another by a change of coordinates \(\tilde \z (\z)\), then there is a meromorphic function \(\tilde f(\tilde \z)\) such that
\[
\om = \tilde f(\tilde z) d \tilde \z = f(\z)d\z,
\]
where \(d\tilde \z = \tilde \z'(\z)d\z,\) and \(\tilde \z'(\z)\) is a holomorphic function of \(\z\).

We will say that a meromorphic one-form \(\om\) on \(S\) is a {\em holomorphic one-form} on \(S\) if, for each local coordinate \(\z\), \(\om = f(\z)d\z\), where the coefficient function \(f(\z)\) is holomorphic.

If we look at the integrands of the examples of Abelian integrals that we discussed above, then it is easy to check that
\bean
\om_1&=&\frac{dz}{\sqrt{(z^2-1)(z^2-k^2)}},\\
\om_k &=& \frac{z^kdz}{\sqrt{(z-a_1)(z-a_2)\cdots(z-a_{2p})}}, k=0,\dots,p-1
\eean
are indeed holomorphic one forms on the respective Riemann surfaces of genus 1 and genus \(p\).  First of all, these are indeed meromorphic one-forms as they are defined, and it remains to show that they are holomorphic near any of the singular points (where these one-forms have potential poles). The calculations are essentially the same that we used to show that these integrals had well-defined values at the singular points. 

For instance, near the singular point \(z=a_1\) for the one form \(\om_k\) in the second case, we have
\bean
\z&=& (z-a_1)^{\frac{1}{2}},\\
d\z &=& \frac{1}{2}(z-a_1)^{-\frac{1}{2}}dz,
\eean
which gives
\bean
z&=& \z^2 +a_1,\\
dz&=& 2\z d\z,
\eean
and hence
\[
\om_k = \frac{2(\z^2+a_1)^kd\z}{\sqrt{(\z^2+a_1-a_2)(\z^2+a_1-a_3)\cdots(\z^2+a_1-a_{2p})}},
\]
which is holomorphic near \(z=0\).  It still has to be checked that \(\om_k\) is holomorphic at \(\infty\), and this is a similar calculation.

Using this terminology we see that the number of linearly independent Abelian integrals of the first kind on a Riemann surface \(S\) is the same as the number of linearly independent holomorphic one-forms on \(S\).  This point of view became standard in modern treatments of Riemann surfaces using differential forms to represent the connectivity (the genus) using de Rham's theorem, and using the holomorphic one-forms to be a way of representing this topology when the complex structure is assumed.  This is the essence of Hodge theory on general K\"{a}hler manifolds (see, e.g., \cite{wells2008}).

In the remainder of his main Abelian function paper \cite{riemann1857d}, Riemann formulated and solved several geometric problems which have come to have major significance in the subsequent development of Riemann surfaces, algebraic geometry and complex manifolds.  In addition he formulated and gave his version of a solution to the Jacobi inversion problem, which was concerned with the  generalizations of the inverses of elliptic integrals to inverses of Abelian integrals.  We will look at the geometric problems first and subsequently return to the Jacobi inverse problem.

The first geometric problem he formulated and resolved was to show that the Riemann surfaces of the kind he had formulated as a branched covering of \(\overline \BC\) could be realized as the Riemann surface defined by an algebraic function \(F(z,w)\).  In modern terms the question he raised could be reformulated to ask if an abstract Riemann surface could be realized as projective algebraic submanifold of complex projective space.  The answer to this question is that it is indeed possible, and this is a special case of and simple consequence of the Kodaira embedding theorem for compact complex manifolds (see, e.g., \cite{wells2008}).

The second problem he formulated was to consider the birational equivalence of solutions of equations of the form \(F(z,w)=0\) for different choices of the polynomial \(F(z,w)\).  Let's formulate this question somewhat more precisely.  Let
\[
C= \{(z,w)\in \BC^2 : F(z,w)=0\},
\]
where \(F(z,w)\) is a polynomial in the variables \(z\) and \(w\).  We call such a \(C\) an {\em algebraic curve}, and we include in \(C\) the extension of \(C\) to include points at infinity (which is easy to do using homogeneous coordinates, and became standard in algebraic geometry shortly after the time of Riemann%
\footnote{For instance, the text by Clebsch and Gordan on Abelian functions \cite{clebsch-gordan1866}, published in 1866, nine years after Riemann's fundamental papers of 1857 formulated this theory in terms of homogeneous coordinates.}%
).

Two algebraic curves
\bean
C &=& \{F(z,w)=0\},\\
C_1&=&\{F_1(z_1,w_1)=0\}\},
\eean
are {\em birationally equivalent} if there exist rational mappings
\[
z_1(z,w),z_2(z,,w),
\]
with inverse rational mappings 
\[z(z_1,w_1), w(z_1,w_1),
\]
such that 
\[
F_1(z_1,w_1)= F(z(z_1,w_1),w(z_1,w_1)),
\]
and 
\[
F(z,w) = F_1(z_1(z,w),w_1(z,w)).
\]
Riemann formulated the problem of classifying equivalence classes of birationally equivalent algebraic curves, and, by dimensional analysis of the parameters, he concluded that for a given algebraic curve of genus \(p\) (defined by an equation \(F(z,w)=0\)):
\bi
\item For \(p=0\), all curves are birationally equivalent.
\item For \(p=1\), there is a one (complex) parameter family of birationally inequivalent algebraic curves.
\item For \(p>1\), there is a \(3p-3\)-parameter family of birationally inequivalent algebraic curves.
\ei
Riemann termed these parameters the {\em moduli} of the algebraic curves, and a major problem in mathematical research over the next century became to understand the nature and representations of these moduli.  For Riemann the classification of algebraic curves was equivalent to the classification of Riemann surfaces associated to these algebraic curves.  

In algebraic geometry the classification of algebraic varieties of one or more dimensions has been a very important research topic ever since the time of Riemann.  In the theory of complex manifolds (most complex manifolds do not correspond to solutions of algebraic equations) the deformations of complex structures on manifolds of one or more dimensions has been an equally rich field of research, with the theory of Teichm\"{u}ller spaces playing an important role in the contemporary theory of moduli of Riemann surfaces.

Finally, we discuss briefly the main topic Riemann was addressing in this series of papers, the Jacobi inversion problem.  Namely, let
\[
F(z,w)=0,
\]
define a Riemann surface \(S\) of genus \(p\), and let \(A_1(z),\dots,A_p(z)\) be \(p\) linearly independent Abelian integrals of the first kind on \(S\).  Each of these are multivalued functions and the value of each of these functions at a given point \(z\) depends on the path of integration from some fixed initial point to the point on the Riemann surface whose coordinate in the extended complex plane is \(z\), which is implicit in the definition of each \(A_k\).  Define the functions 
\be
\label{eqn:abelian-integral-sum}
\ba{c}
v_1(z_1,\dots,z_p) = A_1(z_1)+\cdots + A_1(z_p),\\
v_2(z_1,\dots,a_p) = A_2(z_1)+\cdots+A_2(z_p),\\
\vdots\\
v_p(z_1,\dots,z_p)= A_p(z_1)+\cdots+A_p(z_p).
\ea
\ee
The {\em Jacobi inversion problem} is to find an inverse to this mapping and determine its properties.  An inverse mapping here would be \(p\) functions
\[
\ba{c}
z_1(v_1,\dots,v_p),\\
z_2(v_1,\dots,v_p),\\
\vdots\\
z_p(v_1,\dots,v_p),
\ea
\]
which provide an inverse to the mapping described in (\ref{eqn:abelian-integral-sum}). For the case of \(p=1\), we have only one Abelian integral to deal with, and this is an elliptic integral, and its inverse is an elliptic function as discovered by Abel and Jacobi, and which we discussed at length in Section \ref{sec:abel-theorem}.  Thus the inverse functions \(\{z_1(v_1,\dots,v_p),z_2(v_1,\dots,v_p),\dots\}\) would be functions of \(p\) complex variables which would be generalizations of elliptic functions, presumably with some periodicity properties of the same sort as elliptic functions have.  We shall see that the inverse functions do exist and are meromorphic functions in \(\BC^p\) with \(2p\) independent periods.  Riemann called such functions {\em Abelian functions} in honor of Abel, who had studied in the various versions of Abel's theorem sums of Abelian integrals of the sort that appear in (\ref{eqn:abelian-integral-sum}).

One measure of the multivalued nature of the Abelian integrals used in (\ref{eqn:abelian-integral-sum}) is to use what are called the {\em periods of an Abelian integral}, and these will turn out to be the periods of the corresponding Abelian functions described briefly above.  Let \(A_k(z)\) be \(k\) linearly independent (over the real numbers) Abelian integrals of the first kind on a Riemann surface \(S\) of genus \(p\), \(k=1,\dots,p\),  and let \(\g_1,\dots,\g_{2p},\) be closed curves ({\em cycles}) on \(S\) which represent cuts which reduce the Riemann surface \(S\) to a simply-connected open subset \(S'\subset S\), then let
\be
\label{eqn:abelian-periods}
\om_k^j:=\int_{\g_j} \a_k(z),
\ee
where \(\a_k(z)\) is the {\em integrand} of \(A_k(z)\) considered as a holomorphic one-form on \(S\), as was illustrated in our examples earlier.  These are the {\em periods} of these Abelian integrals \(A_k(z)\) for this choices of cycles \(\g_1,\dots,\g_{2p}.\) It follows that if \(\g\) and \(\tilde \g\) are any two paths joining an initial point \(z_0\) to a variable point \(z\) on \(S\), and if \(A_k(z)\) represents the value of the Abelian integral along the path \(\g\) and \(\tilde A_k (z)\) represents the value of the same Abelian integral along the path \(\tilde \g,\) then 
\[
\tilde A_k(z) = A_k(z) + m_1\om_k^1+\cdots + m_k^{2p}\om_k^{2p},
\]
where \(m_k^j\) are integers.  Thus, the periods \(\om_k^1,\dots,\om_k^{2p},\) represent precisely the multivalued nature of the Abelian integral \(A_k(z)\).

We now define a {\em meromorphic function} \(f(z)=f(z_1,\dots,z_p)\) on an open domain \(D\subset \BC^p, p\ge 1\), to be a holomorphic function on \(D-S\), where \(S\) is a closed lower-dimensional subset of \(D\) such that near any point \(z^0=(z_1^0,\dots,z^0_p)\) of \(D\), \(f\) can be represented as the quotient of two holomorphic functions.  The singular points \(S\) correspond to the points where the local holomorphic function in the denominator is zero.  Thus the set \(S\) is generically a \(p-1\)-dimensional locally defined complex submanifold of \(\BC^p\). For instance if we set
\[
f(z_1,z_2)=\frac{z_1-z^0_1}{z_2-z^0_2},
\]
then \(f(z)\) is meromorphic on \(\BC^2\), and it has zeros on the line \(z_1=z^0_1\), and it has poles along the line \(z_2=z^0_2\), and it has no well defined value at the singular point \((z^0_1,z^0_2)\).

Consider a meromorphic function \(f(z)\) on \(\BC^p\). Let \(\om\in \BC^p\) be a fixed complex \(p\)-tuple, \((\om_1,\dots,\om_p)\ne 0\).  Then we say the function \(f(z)\) is {\it periodic with respect to the period \(\om\) } if 
\[
f(z+\om)=f(z), \; \textrm{for all}\; z\in\BC^p.
\]
  Let us define the vectors 
\[
\om^j = (\om^j_1,\dots,\om^j_p), j= 1,\dots,2p,
\]
where the \(\om^j_k\) are defined as the periods of the Abelian integrals as in (\ref{eqn:abelian-periods}), then Riemann, Weierstrass and others showed that there exist meromorphic functions \(f(z)\) on \(\BC^p\), such that 
\[
f(z+m_1\om^1+\cdots+m_{2p}\om^{2p})=f(z), z\in \BC^p,
\]
where the \(m_j\) are arbitrary integers.  More precisely, there exist functions with \(2p\) independent periods (namely the periods defined above are linearly independent over the real numbers), and they also showed there are no functions with more than \(2p\) periods. 

Riemann used theta-functions in \cite{riemann1857d} to demonstrate the existence of Abelian functions, in a manner similar to our discussion of elliptic functions in Section \ref{sec:abel-theorem}.  Weierstrass formulated (by differentiation) (\ref{eqn:abelian-integral-sum}) as differential equations and solved these using power series methods. Both solutions gave great impetus to further research in the rich theory of Abelian functions during the latter half of the 19th century.  We recommend highly the interesting book by Markushevich \cite{markushevich1992} which gives a detailed and very well written history of the early development of both elliptic and Abelian functions.

\section{Conclusion}
\label{sec:conclusion}
In the preceding sections of this paper we have seen how the complex plane evolved into the concept of a Riemann surface and how the special class of holomorphic functions began to play an important role in analysis. There were  several other significant ideas which arose in the 19th century which played an important role in complex geometry.  The first of these was the development of projective geometry and more specifically the notion of projective space, a generalization of classical Euclidean space which evolved over numerous decades of the 19th century.  In an earlier paper  \cite{wells2013} we  traced the  development of this and other geometric ideas during this time period.  For an outstanding historic reference to this development, we recommend Felix Klein's beautiful lectures on non-Euclidean geometry from the end of the 19th century  which were  published in 1928 \cite{klein1928}. In modern complex geometry, complex \(n\)-dimensional projective space \(\BP_n(\BC)\) plays a very important role.

At the end of the 19th century all the ingredients had been developed which allowed Hermann Weyl (1885--1955)  to develop in 1913 the first theory of abstract manifolds in the important special case of abstract Riemann surfaces.  We will discuss this in more detail towards the end of this conclusion.  We now discuss some fundamental developments that preceded Weyl's work.

The first of these, the theory of transformation groups, has become more well-known under its modern appellation of {\it Lie groups}.  The first study of Lie groups arose as transformation groups of specific geometric spaces in various papers of Felix Klein (1849--1925) and Sophus Lie (1842--1899).  For instance in 1871 they wrote a joint paper \cite{klein-lie1871} which referred to two earlier papers each of them had written that dealt with  transformation groups on quite specific geometric spaces (\cite{klein1870},\cite{lie1870}. In 1872 Klein wrote his famous {\it Erlangen Program} paper \cite{klein1872} in which he outlined, among other things, the role he foresaw for transformation groups, or more generally, continuous groups and their subgroups (in particular discrete subgroups) to play in geometry. This turned out to be a very significant paper, and the study of the action of Lie groups on manifolds became  an important topic in the 20th century. In 1880 Lie published the first of his fundamental papers on what became the theory of Lie groups entitled ``Theorie der Transformationsgruppen I" \cite{lie1880}.   In this paper he classified the Lie groups acting on Euclidean spaces of one and two dimensions and developed the tools of Lie algebras as means of determining the classification.  He gives a summary in this paper of all earlier references known to him at the time concerning this generic topic.  One can leaf through the collected works of both Lie and Klein (\cite{klein1921}, \cite{lie1922}, which are all available on-line today) to get a good overview of the development of transformation groups and their role in studying geometry in a wide variety of contexts.

A very important relation between transformation groups to complex geometry came at the end of the 19th century in what became known as the {\it uniformization theorem} of Riemann surfaces, representing all compact Riemann surfaces as quotients of the three distinct simply-connected Riemann surfaces (the Riemann sphere, the complex plane, and the unit disc) by discrete groups of the biholomorphic automorphisms of these spaces.   Let us give a brief summary of this important work, which became a role model for many similar questions in higher dimensions.

The history of non-Euclidean geometry geometry has been well documented (see the classical treatment by Klein \cite{klein1928}, for instance).  In the 19th century there were various discoveries of geometries that were not Euclidean with abstract axiomatic systems which had variations on the parallel axiom and, more particularly, specific {\it models} of a given geometry. Here we will only mention that the developments of complex geometry in the 19th century led to three very specific models of the three types of geometries that have evolved. Namely, first the complex plane \(\BC\) with its Euclidean metric 
\[
ds^2 = dzd\overline{z}= dx^2 +dy^2
\]
is a model for the classical Euclidean plane {\it plane geometry}, with its geodesics being the classical straight lines, and the Euclidean translations and rotations being given by \(z \mapsto z+a\), and rotations \(z \mapsto ze^{i\theta}\).  The second case of {\it elliptic geometry} is  represented in terms of the two sphere \(S^2\), which  can be described in complex terms as the Riemann sphere, that is, the one-point compactification of the complex plane \(\overline{\BC}= \BC\cup \infty\),which is biholomorphic to one-dimensional complex projective space \(\BP_1(\BC)\). If we let the complex plane be the standard coordinate patch (the complement of the point at infinity) with the metric
\[
ds^2 = \frac{dzd\overline{z}}{(1+ |z|^2)^2}
\]
then we have the geodesics are the great circles and the transformation group of (holomorphic) isometries is the orientation-preserving rotations of the sphere \(SO(3) \cong PSU(2)\), in terms of real and complex coordinates respectively. As a model for non Euclidean geometry, it is necessary to consider the two sphere with antipodal points identified (which gives two-dimensional real projective space, so that through any two points there is only one geodesic joining them, one of the axioms of all the geometries). This projective space structure doesn't preserve the complex structure.  Finally, we have the very important case of {\it hyperbolic geometry}, which can be modeled in terms of the Poincar\'{e} disk \(\D\), which is the unit disk in the complex plane
\[
\D := \{z\in \BC: |z| < 1\},
\]
equipped with the Poincar\'{e} metric
\[
ds^2= \frac{dzd\overline{z}}{(1-|z|^2)^2} = \frac{dx^2 +dy^2}{(1-(x^2+y^2)^2}.
\]
Here the geodesics are the arc of circles in the unit disc which have endpoints on the unit circle and which are orthogonal to the unit circle at those points.  The transformation group of holomorphic isometries is \(\mathrm{SU}(1,1)\), which can be represented as the set of M\"{o}bius transformations of the form
\[
z \mapsto e^{i\theta}\left(\frac{z-a}{-\overline{a}z+1}\right), \th \in \BR, |a| < 1.
\]

As is well known, these three two-dimensional models of Euclidean and non-Euclidean geometry \(\BC, \overline{\BC},\) and \(\D\) give also the complete classification of the connected and simply-connected complex manifolds of dimension one, which was finally proved satisfactorily at the beginning of the 20th century by Henri Poincar\'{e} (1854-1912) \cite{poincare1907} 
and  Paul Koebe (1882-1945) \cite{koebe1907c}.  
This theorem is referred to in the literature as the {\it uniformization theorem}%
\footnote{See the very informative historical paper by Jeremy Gray \cite{gray1994} on the history of both the Riemann mapping theorem and the uniformization theorem.}%
.
 Moreover, any compact Riemann surface of genus 1 is equivalent to the quotient of the complex plane by a lattice (a complex torus of dimension 1) , and any compact Riemann surface of genus \(g>1\) is equivalent to  \(\D/\G\), where \(\G\) is a properly discontinuous subgroup of \(\mathrm{SU}(1,1)\), the automorphisms of the unit disc (see, e.g., the extensive survey paper by Lipman Bers \cite{bers1972}).

The final developments of the 19th century critical for complex geometry primarily concerned itself with the development of topology, both point-set topology and algebraic topology, and finally the abstraction of the notion of a manifold (as mentioned above).

First came the development of {\it set theory} by Georg Cantor (1845 --1918) in the 1870's, which led to the development of abstract topological spaces in the early 20th century.
Cantor's work turned out to be revolutionary for all of mathematics as well as leading to the famous continuum hypothesis and fundamental questions in the foundations of mathematics which we won't discuss here (see, for instance,  the collected works of Cantor \cite{cantor1932} with its interesting introduction Ernst Zermelo (1871--1953), and modern surveys of this important topic).  

This led to the development of {\it abstract topological spaces}. Maurice Fr\'{e}chet (1878--1973) was the first to formulate an abstract topological space \cite{frechet1906} (he used the notion of spaces of type \((L)\), which had axiomatically sequences of elements which either converged or didn't, and a few years later more general notions of a topological space (more general than that of Fr\'{e}chet) were formulated by Felix Hausdorff (1868--1942)  \cite{hausdorff1914} using axioms of neighborhoods as the fundamental notion, including, in addition, his Hausdorff separation axiom, and then in 1922 Casimir Kuratowsi (1896--1980) \cite{kuratowski1922} provided the most general theory of topological spaces (using axioms concerning closed sets as a basis for the theory).  Today we use axioms for open sets as the basis for the theory, which is equivalent to Kuratowski's theory.

Secondly was the development of {\it algebraic topology of manifolds}.  This started with the work of Riemann on Riemann surfaces (as discussed in Section \ref{sec:riemann-surfaces}), where he developed the notion of connectivity for Riemann surfaces.  This was generalized by Betti in 1871 \cite{Betti1871}, who generalized Riemann's connectivity for two-dimensional manifolds to what are now called Betti numbers of higher-dimensional manifolds. Finally, Poincar\'e, in a fundamental series of papers at the end of the 19th century, formulated the fundamental principles of what has become known as algebraic topology (see the collection of papers on topology in Volume VI of Poincar\'{e}'s collected works \cite{poincare1953} and in particular the translation of Poincar\'{e} 's topology papers into English by John Stillwell \cite{poincare2010} with its  lucid introduction to the whole topic).  In these papers Poincar\'e considered manifolds which were smooth submanifolds of Euclidean space of any dimension with or without boundary, and if with boundary, the boundaries were piece-wise smooth, all defined in terms of defining functions in the ambient space.  

The final development concerns the creation of the notion of an {\it abstract manifold} in a mathematically satisfactory way.  In Riemann's original paper \cite{riemann1854}, he discussed in a philosophical but not technical manner this concept.  For Riemann at the time it was simply something (not defined) which had local coordinate charts with suitable transition functions, and he worked from there. The fundamental creation of an abstract manifold (with topological, differentiable or complex structures) was taken by Hermann Weyl in the first edition of his  famous book on Riemann surfaces ({\it Die Idee der Riemannshen Fl\"ache} \cite{weyl1913}). He formulates for the first time in a rigorous manner the notion of an abstract topological manifold,  in the two-dimensional setting (he notes that he could work in more dimensions, but he was concerned with a new way of looking at Riemann surfaces).  He describes a (two-dimensional) manifold as any point set M with a set of neighborhoods satisfying suitable axioms and such that for each point \(p\in M\) there is a neighborhood \(U\) of \(p\) which is homeomorphic to a disc in \(R^2\).%
\footnote{Weyl had the key to an abstract topological space here, but he did not pursue it further.}%

The key here is that he starts with an abstract set, and this would only have been possible after Cantor created set theory in a manner that could be used in all parts of mathematics. He makes a further assumption to those made earlier that the manifold is triangulated and goes on to define a Riemann surface to be a triangulated topological manifold which has local coordinate systems which map to discs in the complex plane whose overlap transformations are holomorphic. He makes a point that algebraic topology will play an important role in the theory of manifolds (hence the use of triangulations), noting the earlier work of Riemann and Poincar\'e. He shows later in the book that the extra hypothesis of being triangulated for Riemann surfaces can be eliminated (his book came out in a second edition in 1955 \cite{weyl1955} and there is an English translation as well \cite{weyl2009}).   His second edition uses the now more common version of being a topological space (set with neighborhoods satisfying axioms) which is Hausdorff and has a countable basis for the topology.  This is the point of view that was formulated by Oswald Veblen (1880--1960) and J. H. C. Whitehead (1904--1960) who formulated precisely the notion of an abstract differentiable manifold of arbitrary dimension \cite{veblen-whitehead1931} which they explored in greater detail in their monograph \cite{veblen-whitehead1932}.  This was the definition used by Hassler Whitney in his embedding theorem paper of 1936 \cite{whitney1936}, as more pathological manifolds would not be suitable for embeddings into Euclidean space, and this has become standard in the study of differentiable and complex manifolds (see for example \cite{wells2008}).

We close by quoting from Weyl's 1913 book about his new way of looking at what has become known as the theory of manifolds.
\begin{quote}
Eine solche
strenge Darstellung, die namentlich auch bei Begr\"{u}ndung der fundamentalen,
in die Funktionentheorie hineinspielenden Begriffe und S\"{a}tze der
Analysis situs sich nicht auf anschauliche Plausibilit\"{a}t beruft, sondern
mengentheoretisch exakte Beweise gibt, liegt bis jetzt nicht vor. Die
wissenschaftliche Arbeit, die hier zu erledigen blieb, mag vielleicht als
Leistung nicht sonderlich hoch bewertet werden. Immerhin glaube ich
behaupten zu k\"{o}nnen, daß ich mit Ernst und Gewissenhaftigkeit nach
den einfachsten und Sachgemäßesten Methoden gesucht habe, die zu dem
vorgegebenen Ziele f\"{u}hren; und an manchen Stellen habe ich dabei andere
Wege einschlagen m\"{u}ssen als diejenigen, die in der Literatur seit dem
Erscheinen von C. Neumanns klassischem Buche \"{u}ber „Riemanns Theorie
der Abelschen Integrale" (1865) traditionell geworden sind.%
\footnote{Such a rigorous presentation, which, namely by the establishing of the fundamental concepts and theorems in function theory and using theorems of the analysis situs which don't just depend on intuitive plausibility, but have set-theoretic exact proofs, does not exist. The scientific work that remains to be done   in this regard, may perhaps not be be particularly highly valued.  But, nevertheless, I believe I can maintain that I have tried in a serious and conscientious manner to find the simplest and most appropriate methods, that lead to the  asserted goal; and at many points I have had to proceed in a different manner than that which has become traditional in the literature since the appearance of C. Neumann's classical book about "Riemann's theory of Abelian Integrals."}
\end{quote}

\bibliography{references}
\bibliographystyle{plain}%

\end{document}